\documentclass[12pt]{amsart}

\usepackage{nccmath}
\usepackage{empheq}

\usepackage{amssymb, amsmath, amsthm,commath,xspace}

\usepackage[utf8]{inputenc}
\usepackage[T1]{fontenc}

\usepackage{dsfont}
\usepackage[labelfont={normalsize}]{caption,subfig}
	
\usepackage{microtype}
\usepackage{times}

\usepackage{subfiles}
\usepackage[normalem]{ulem}

\usepackage[usenames,dvipsnames]{xcolor}
\usepackage{tikz} 
\usetikzlibrary{patterns}
\usetikzlibrary{patterns.meta}

\tikzdeclarepattern{
    name=hatch,
    parameters={\hatchsize,\hatchangle,\hatchlinewidth},
    bounding box={(-.1pt,-.1pt) and (\hatchsize+.1pt,\hatchsize+.1pt)},
    tile size={(\hatchsize,\hatchsize)},
    tile transformation={rotate=\hatchangle},
    defaults={
        hatch size/.store in=\hatchsize,hatch size=5pt,
        hatch angle/.store in=\hatchangle,hatch angle=0,
        hatch linewidth/.store in=\hatchlinewidth,hatch linewidth=.4pt,
    },
    code={
        \draw[line width=\hatchlinewidth] (0,0) -- (\hatchsize,\hatchsize);
    }
}

\usepackage{enumitem}
\usepackage{thmtools} 
\usepackage{thm-restate}

\usepackage{hyperref}

\usepackage[capitalize,noabbrev]{cleveref}

\usepackage{nicefrac}

\usepackage{thm-restate}

\usepackage[backend=bibtex, isbn=false, 
giveninits=true, style=alphabetic, maxalphanames=4, 
maxnames=5]{biblatex}

\DeclareFieldFormat{url}{%
    \iffieldundef{doi}{%
        \mkbibacro{URL}\addcolon\space\url{#1}%
    }{%
    }%
}

\DeclareFieldFormat{urldate}{%
    \iffieldundef{doi}{%
        \mkbibparens{\bibstring{urlseen}\space#1}%
    }{%
    }%
}

\numberwithin{equation}{section}

\usepackage[backgroundcolor=yellow]{todonotes}

\author[{\L}ukasz Ma{\'s}lanka]{{\L}ukasz Ma{\'s}lanka}
\address{Institute of Mathematics, Polish Academy of Sciences, {\'S}niadeckich~8, 00-656 Warszawa, Poland}
\email{lmaslanka@impan.pl}

\author[Piotr {\'S}niady]{Piotr {\'S}niady}
\address{Institute of Mathematics, Polish Academy of Sciences, {\'S}niadeckich~8, 00-656 Warszawa, Poland}
\email{psniady@impan.pl}

\makeatletter
\@ifclasswith{amsart}{externalize}%
{
    \usetikzlibrary{external} 
    \tikzexternalize[
    optimize=false,
    prefix=tikz/,
    mode=list and make]
    \makeatletter
    \renewcommand{\todo}[2][]{\tikzexternaldisable\@todo[#1]{#2}\tikzexternalenable}
    \makeatother
}%
{
    \newcommand{\tikzexternaldisable}{}
    \newcommand{\tikzexternalenable}{}
}
\makeatother

\newcommand{\Time}{T}

\definecolor{Set1-3-1}{RGB}{228,26,28}
\definecolor{Set1-3-A}{RGB}{228,26,28}
\definecolor{Set1-3-2}{RGB}{55,126,184}
\definecolor{Set1-3-B}{RGB}{55,126,184}
\definecolor{Set1-3-3}{RGB}{77,175,74}
\definecolor{Set1-3-C}{RGB}{77,175,74}

\newtheorem{theorem}{Theorem}[section]
\newtheorem{assumption}[theorem]{Assumption}
\newtheorem{fact}[theorem]{Fact}

\newtheorem{problem}[theorem]{Problem}
\newtheorem{lemma}[theorem]{Lemma}
\newtheorem{proposition}[theorem]{Proposition}
\newtheorem{corollary}[theorem]{Corollary}
\newtheorem{conjecture}[theorem]{Conjecture}

\theoremstyle{definition}

\theoremstyle{remark}
\newtheorem{remark}[theorem]{Remark}

\newcommand{\bigforall}{\mbox{\Large $\mathsurround0pt\forall$}}

\newcommand{\on}{\operatorname}

\newcommand{\Y}{\mathbb{Y}}
\newcommand{\N}{\mathbb{N}}
\newcommand{\Z}{\mathbb{Z}}
\newcommand{\R}{\mathbb{R}}
\newcommand{\C}{\mathbb{C}}
\newcommand{\E}{\mathbb{E}}   %

\newcommand{\F}{\mathcal{F}}

\newcommand{\ls}{\left\langle}
\newcommand{\rs}{\right\rangle}

\newcommand{\pos}{\on{pos}}
\newcommand{\sh}{\on{sh}}

\newcommand{\Ps}{\mathcal{P}}  %
\newcommand{\Pp}{\mathbb{P}_N}    %
\newcommand{\Pt}{\widetilde{\mathbb{P}}_N}   %
\newcommand{\sq}{\square}   

\newcommand{\tableaux}{\mathcal{T}}
\newcommand{\widetildetableaux}{\widetilde{\mathcal{T}}}
\newcommand{\tab}{T}
\newcommand{\multi}{M}
\newcommand{\water}{w}
\newcommand{\content}{u}   %

\newcommand{\supp}{\on{supp}}

\newcommand{\ttt}{\mathbf{T}'_N}
\newcommand{\ttilde}{\mathbf{M}'_N}

\newcommand{\tim}{t}

\newcommand{\point}{P}
\newcommand{\pat}{\mathbf{q}}

\newcommand{\Sym}{\mathfrak{S}}

\newcommand{\mysquare}{[0,1]^2}
\newcommand{\mypoint}{X}

\newcommand{\sigmapj}{\sigma_{p,j}}

\newcommand{\usc}{u}
\newcommand{\theor}{\operatorname{th}}

\newcommand*\diff{\mathop{}\!\mathrm{d}}

\newcommand{\SPar}[1]{\on{SPar}({#1})}

\newcommand{\tmeas}{\mu}
\newcommand{\cotmeas}{\nu}

\DeclareMathOperator{\id}{id}
\DeclareMathOperator{\evac}{evac}
\DeclareMathOperator{\revevac}{rev evac}

\newcommand{\jdtcomplete}{\partial^*}

\newcommand{\jdtt}{dual promotion\xspace}

\newcommand{\jdtp}{sliding path\xspace}
\newcommand{\jdtps}{sliding paths\xspace}

\newcommand{\jdtincomplete}{jeu de taquin\xspace}

\newcommand{\Longitude}{\Psi}
\newcommand{\LongitudeTilde}{\widetilde{\Psi}}
\newcommand{\longitude}{\psi}
\newcommand{\TheoreticalLongitude}{\Psi^{\operatorname{th}}_N}

\newcommand{\empiricallongitude}{G_N}
\newcommand{\dynamicempiriciallongitude}[1]{\widetilde{G}_N^{\: #1}}

\newcommand{\empiricalmeasure}{m_N}

\begin{document}


\title[Second class particles, evacuation and sliding paths]{%
Second class particles \\ and limit shapes  of  evacuation and sliding paths \\ for random tableaux}


\begin{abstract}
We investigate two closely related setups. In the first one we consider a
TASEP-style system of particles with specified initial and final
configurations. The probability of each history of the system is assumed to be
equal. We show that the rescaled trajectory of the \emph{second class particle}
converges (as the size of the system tends to infinity) to a random arc of an
ellipse.

In the second setup we consider a uniformly random Young tableau of
square shape and look for typical (in the sense of probability) \emph{sliding
    paths} and \emph{evacuation paths} in the asymptotic setting as the size of
the square tends to infinity. We show that the probability distribution of
such paths converges to a random \emph{meridian} connecting the opposite
corners of the square. We also discuss analogous results for non-square Young
tableaux.
\end{abstract}

\maketitle

\section{Introduction}
\label{sec:introduction}

This article is the full version of a 12-page extended abstract \cite{Maslanka2020}
which was published in the proceedings of the 
\emph{32nd International Conference on Formal Power Series 
and Algebraic Combinatorics, FPSAC 2020}.

\medskip

The results of the current paper concern two distinct setups which are closely connected. 
The first one, presented in \cref{sec:scp}, involves a certain interacting particle system. 
The second one, presented in \cref{sec:tableaux}, involves random Young tableaux.

\subsection{TASEP system with the uniform distribution over histories}
\label{sec:scp}

\subsubsection{The setup}
\label{sec:setup-scp}

\begin{figure}[t]
\centering
	\subfloat[]{
	\resizebox{\linewidth}{!}{
		\begin{tikzpicture}
		\draw[ultra thick,->] (-4.5,0) -- (10,0) node[anchor=west]{$x$};
		\foreach \x/\t in {-4/-N+1,-2/-2,-1/-1,0/0,1/1,2/2, 9/M-1}
		{ \draw (\x,-0.35)  node[anchor=north] {\small $ \t  $}; }
		
		\draw (5,-0.35) node[anchor=north] {\footnotesize $M-N$};

		\foreach \x in {-4,...,9}
		{ \draw[ultra thick] (\x,-0.4) -- (\x,0.4); }

		\foreach \x in {1,...,9}
		{ \draw[fill=white] (\x,0 ) circle (0.3); }

		\foreach \x in {-4,...,-1}
		{ \draw[fill=blue!50] (\x,0 ) circle (0.3); }

		\draw[preaction={fill,white},pattern=north west lines, pattern color=red] (0,0) circle (0.3);

		\end{tikzpicture}
		\label{fig:scpA}
		}
	}

	\vspace{3ex}

	\subfloat[]{
	
	\resizebox{\linewidth}{!}{
			\begin{tikzpicture}
		\draw[ultra thick,->] (-4.5,0) -- (10,0) node[anchor=west]{$x$};
		\foreach \x/\t in {-4/-N+1,-2/-2,-1/-1,0/0,1/1,2/2, 9/M-1}
		{ \draw (\x,-0.35)  node[anchor=north] {\small $ \t  $}; }
		
		\draw (5,-0.35) node[anchor=north] {\footnotesize $M-N$};
			
			\foreach \x in {-4,...,9}
			{ \draw[ultra thick] (\x,-0.4) -- (\x,0.4); }    
			
			\foreach \x in {-4,...,4}
			{ \draw[fill=white] (\x,0 ) circle (0.3); }
			
			\foreach \x in {6,...,9}
			{ \draw[fill=blue!50] (\x,0 ) circle (0.3); }
			
			\draw[preaction={fill,white},pattern=north west lines, pattern color=red] (5,0) circle (0.3);
			
			\end{tikzpicture}
		\label{fig:scpB} 
	} 
	
	}

	\caption{\protect\subref{fig:scpA} The initial configuration
	of the particle system. The first class particles are depicted as filled blue circles, the holes
	are depicted as empty circles. The striped red circle denotes the second class
	particle. 
	\protect\subref{fig:scpB} The final configuration of the particle system. 
	}

	\label{fig:scp}

\end{figure}

\begin{figure}
\subfloat[]{
\begin{tikzpicture}
\draw[ultra thick,->] (-1,0) -- (2,0);
 \foreach \x in {0,1}
{ \draw[ultra thick] (\x,-0.4) -- (\x,0.4); }

\draw[fill=blue!50] (0,0 ) circle (0.3); 
\draw[fill=white] (1,0 ) circle (0.3);
\draw[thick,red,->] (0,0.7) to[out=90,in=90] (1,0.7);

\draw (0.5,-1cm) node {$\Downarrow$};

\begin{scope}[yshift=-2cm]
\draw[ultra thick,->] (-1,0) -- (2,0);
\foreach \x in {0,1}
{ \draw[ultra thick] (\x,-0.4) -- (\x,0.4); }

\draw[fill=blue!50] (1,0 ) circle (0.3); 
\draw[fill=white] (0,0 ) circle (0.3);
\end{scope}
\end{tikzpicture}
\label{fig:transitionA}
}
\hfill
\subfloat[]{
    \begin{tikzpicture}
    \draw[ultra thick,->] (-1,0) -- (2,0);
    \foreach \x in {0,1}
    { \draw[ultra thick] (\x,-0.4) -- (\x,0.4); }
    
    \draw[preaction={fill,white},pattern=north west lines, pattern color=red] (0,0 ) circle (0.3); 
    \draw[fill=white] (1,0 ) circle (0.3);
    \draw[thick,red,->] (0,0.7) to[out=90,in=90] (1,0.7);
    
    \draw (0.5,-1cm) node {$\Downarrow$};
    
    \begin{scope}[yshift=-2cm]
    \draw[ultra thick,->] (-1,0) -- (2,0);
    \foreach \x in {0,1}
    { \draw[ultra thick] (\x,-0.4) -- (\x,0.4); }
    
    \draw[preaction={fill,white},pattern=north west lines, pattern color=red] (1,0 ) circle (0.3); 
    \draw[fill=white] (0,0 ) circle (0.3);
    \end{scope}
    \end{tikzpicture}
\label{fig:transitionB}    
}
\hfill
\subfloat[]{
    \begin{tikzpicture}
    \draw[ultra thick,->] (-1,0) -- (2,0);
    \foreach \x in {0,1}
    { \draw[ultra thick] (\x,-0.4) -- (\x,0.4); }
    
    \draw[fill=blue!50] (0,0 ) circle (0.3); 
    \draw[preaction={fill,white},pattern=north west lines, pattern color=red] (1,0 ) circle (0.3);
    \draw[thick,red,->] (0,0.7) to[out=90,in=90] (1,0.7);
    
    \draw (0.5,-1cm) node {$\Downarrow$};
    
    \begin{scope}[yshift=-2cm]
    \draw[ultra thick,->] (-1,0) -- (2,0);
    \foreach \x in {0,1}
    { \draw[ultra thick] (\x,-0.4) -- (\x,0.4); }
    
    \draw[fill=blue!50] (1,0 ) circle (0.3); 
    \draw[preaction={fill,white},pattern=north west lines, pattern color=red] (0,0 ) circle (0.3);
    \end{scope}
    \end{tikzpicture}
\label{fig:transitionC}    
}

\caption{Three allowed types of transitions: \protect\subref{fig:transitionA}
    a first class particle can jump to the right if the next node is free,
    \protect\subref{fig:transitionB} the second class particle can jump to the
    right if the next node is free, \protect\subref{fig:transitionC} the first
    class particle can jump to the right if the next node is occupied by the
    second class particle; in this case the second class particle has to yield
		and the particles exchange their places. }

\label{fig:transition}

\end{figure}

For given integers $N, M \geq 1$ we consider the particle system depicted on
\cref{fig:scp}. There are $N+M-1$ nodes, labeled by the integers from the set
\begin{equation}
    \label{eq:set-of-nodes}
\{-N+1, \dots, 0, 1, \dots, M-1\}.
\end{equation} 
In the initial configuration (which corresponds to the time
$t=1$) the $N-1$ nodes which correspond to the negative integers are occupied by
\emph{the first class particles}, the~$M-1$~nodes which correspond to the positive integers are
empty (or, equivalently, are occupied by \emph{holes}), and the node which
corresponds to zero is occupied by \textbf{\emph{the second class particle}}, see
\cref{fig:scpA}.

In each step exactly one of the following transitions occurs:
\begin{itemize}
\item any particle (first or second class) may jump right to the
next node provided that this node is empty, see
\cref{fig:transitionA,fig:transitionB}, or 
\item a first class particle may jump right to
the next node provided that this node is occupied by the second class particle.
In this case the second class particle has to yield and jumps one node to the
left, see \cref{fig:transitionC}.
\end{itemize}
The \emph{second class particle} is an analogue of a passenger with a cheap
second class ticket who has to yield the seat to any passenger with a more
expensive first class ticket.

\medskip

It is not very difficult to show that no matter which transitions occur, the
system terminates at time $t_{\max}=M N$ (that is after $M N-1$ transitions) in
the configuration shown on \cref{fig:scpB} in which no additional transition is
allowed. 

By the \emph{history} we will understand the information about the
state of the particle system over all values of the time
$t\in\{1,\dots,t_{\max}\}$,
see \cref{fig:story5} for an example.
We consider the finite set of all possible histories of the particle
system and associate to each such a history equal probability. In other words,
we consider a version of the TASEP (which is an acronym for 
Totally Asymmetric Simple Exclusion Process) 
system
\cite{Spitzer1970} with modified transition probabilities.

\begin{figure}
    \centering

    \begin{tikzpicture}[xscale=4,yscale=8]

\begin{scope}[yscale=1/32,xscale=0.17677669529]
\draw[line width=1mm,green!20] (-4,0) rectangle (8,32);    
\draw[dotted] (-4,0) grid (8,32);
\end{scope}

\draw[->] (-1.1,0) -- (1.6,0) node[anchor=north west] {$X$}; 
\draw[->] (0,-0.05) -- (0,1.1) node[anchor=south east] {$T$};

\draw[dashed] plot[smooth]   file {figures/elipsa-full.txt};

\begin{scope}[yscale=1/32,xscale=0.17677669529]
\draw[red,ultra thick] (0,1) -- (0,1) -- (-1,2) -- (-1,2) -- (-2,3) -- (-2,7) -- (-1,8) -- (-1,11) -- (0,12) -- (0,15) -- (1,16) -- (1,19) -- (2,20) -- (2,20) -- (3,21) -- (3,27) -- (4,28) -- (4,28) -- (5,29) -- (5,31)-- (4,32) ; 
\draw[blue,thick]   (-1,1) -- (-1,1) -- (0,2) -- (0,3) -- (1,4) -- (1,5) -- (2,6) -- (2,6) -- (3,7) -- (3,10) -- (4,11) -- (4,16) -- (5,17) -- (5,18) -- (6,19) -- (6,21) -- (7,22) -- (7,32);
\draw[blue,thick]   (-2,1) -- (-2,2) -- (-1,3) -- (-1,4) -- (0,5) -- (0,8) -- (1,9) -- (1,12) -- (2,13) -- (2,14) -- (3,15) -- (3,17) -- (4,18) -- (4,22) -- (5,23) -- (5,24) -- (6,25) -- (6,32);
\draw[blue,thick]   (-3,1) -- (-3,9) -- (-2,10) -- (-2,13) -- (-1,14) -- (-1,23) -- (0,24) -- (0,25) -- (1,26) -- (1,26) -- (2,27) -- (2,29) -- (3,30) -- (3,30) -- (4,31) -- (4,31) -- (5,32) -- (5,32);
\end{scope}

\draw[] (0,1) +(-1pt,0) -- +(1pt,0) node[anchor=west,fill=white,rounded corners=2pt,inner sep=1pt]{$1$};

\draw[] (0,0) +(0,0.5pt) -- +(0,-0.5pt) node[anchor=north,fill=white,rounded corners=2pt,inner sep=2pt]{$0$};

\draw[] (1,0) +(0,0.5pt) -- +(0,-0.5pt) node[anchor=north,fill=white,rounded corners=2pt,inner sep=1pt]{$1$};
\draw[] (-1,0) +(0,0.5pt) -- +(0,-0.5pt) node[anchor=north,fill=white,rounded corners=2pt,inner sep=1pt]{$-1$};
\draw[] (1.41421356237,0) +(0,0.5pt) -- +(0,-0.5pt) node[anchor=north,fill=white,rounded corners=2pt,inner sep=1pt]{$\sqrt{\theta}$};
\draw[] (-1/1.41421356237,0) +(0,0.5pt) -- +(0,-0.5pt) node[anchor=north,fill=white,rounded corners=2pt,inner sep=1pt]{$- \frac{1}{\sqrt{\theta}}$};

\end{tikzpicture}

    \caption{Sample history of the particle system for $N=4$ and $M=8$. The trajectories
        of the first class particles are shown as solid blue lines, the trajectory of the
        second class particle is shown as the thick red line.
        The thick green rectangle indicates the \emph{bounding box}. The dashed line
        indicates \emph{the arctic ellipse} which corresponds to the shape parameter $\theta=\frac{M}{N}=2$.}
    
    \label{fig:story5}    
\end{figure}

\begin{figure}
    \centering
        
    \subfile{figures/kwadrat-jdt-secondclass-N15M30.tex}

    \caption{An analogue of \cref{fig:story5} for $N=15$ and $M=30$.}
    \label{fig:story20}
    
\end{figure}

\subsubsection{Why second class particles?}

Macroscopic quantities describing an interacting particle system (such as
particle density) are usually described by nonlinear partial differential
equations (PDEs for short). For example, the famous Burgers equation
\cite{Burgers1948} describes the~density profile in the~hydrodynamical limit for
the~asymmetric simple exclusion process \cite{Benassi1987}. Weak solutions of
nonlinear PDEs can develop singularities, often referred to as \emph{shocks}.
The shocks can be found by looking for the~crossings of the~characteristic lines
of the~PDE.

Ferrari \cite{Ferrari1992} discovered that
a~second class particle, 
depending on the~place in which 
it begins its journey,
can identify microscopically the location of the shock
or describe the~behavior of the~characteristic lines of 
the~limiting hydrodynamic equation.
Ferrari and Fontes showed in \cite{Ferrari1994} that this
hydrodynamical limit converges to the~traveling wave solution
of the~inviscid Burgers equation.
This connection was later transferred to more general settings 
by Rezakhanlou \cite{Rezakhanlou1995}, Ferrari and Kipnis \cite{Ferrari1995},
Sepp\"{a}l\"{a}inen \cite{Seppaelaeinen2001} and others. 
Furthermore, the~second class particle enters naturally in the~study
of the~fluctuations of the~current of particles \cite{Ferrari1994}.

A~pictorial interpretation of TASEP as 
a~traffic model is given in \cite{Mountford2005}.
The~particles are interpreted as cars on a~single-line highway
with no possibility of passing (which corresponds to the~exclusion rule)
and the shock corresponds to the front of the traffic jam.
The motion of the shock is referred to as the \emph{propagation of the shock}
or the~\emph{rarefaction wave} (or \emph{rarefaction fan})
depending whether the shock moves to the left 
or is being resorpted. 
The second class particle identifies the shock 
and allows to qualitatively describe its motion.

Mountford and Guiol \cite{Mountford2005} studied in fact
a more advanced physical interpretation of the TASEP model
as a moving interface on the plane (space-time).
This gave them a powerful tool to analyze the shocks in
the TASEP process in terms of the last passage percolation.

A~nice heuristic (as well as rigorous) explanation of the~shocks behavior and
the~importance of the~second class particles in this context can be found
in the book of Liggett~\cite[Part~III]{Liggett1999}.

\subsubsection{The asymptotic setup}
\label{sec:setup}

We assume that $(M_i)$ and $(N_i)$ are two sequences of positive integers which tend to infinity and such that their ratio 
\[ \lim_{i\to\infty} \frac{M_i}{N_i} = \theta >0 \]
converges to some positive limit which we call \emph{the shape parameter}.
For $t\in\{1,\dots,M_i N_i\}$ we denote by $\usc_{i}(t)\in\{-N_i+1,\dots,M_i-1\}$ the position of the second class particle at time $t$.
In order to keep the notation lightweight we will sometimes omit the index $i$ and
instead of $M_i,N_i,\usc_i(t)$ we will write shortly $M,N,\usc(t)$.

\medskip

Until now we parameterized the
space using the integer parameter $x\in \{-N+1,\dots,M-1\}$
and the time using the integer parameter $t\in\{1,\dots,MN\}$, however for asymptotic
questions it is more convenient to pass to the rescaled coordinates
\begin{align*} 
      X    &= \frac{x}{\sqrt{MN}} \in \left[ -\frac{1}{\sqrt{M/N}}, \sqrt{M/N} \right], \\
    \Time &= \frac{t}{t_{\max}} \in [0,1], 
\end{align*}
see \cref{fig:story5,fig:story20}.
The rectangle 
\begin{equation}
    \label{eq:rectangle}
    B_\theta=\left\{ (X,\Time) :  X\in \left[-\frac{1}{\sqrt{\theta}}, \sqrt{\theta} \right], \quad \Time \in [0,1] \right\} 
\end{equation}
which shows the range in which the coordinates $X$ and $\Time$ vary asymptotically
will be called \emph{the bounding box}; on \crefrange{fig:story5}{fig:elipsy} it is shown as the
green rectangle.

\newcommand{\traj}{\Xi}
For each value of the parameter $s\in[-1,1]$ we define a function 
\[ \traj_s\colon [0,1] \to \left[-\frac{1}{\sqrt{\theta}}, \sqrt{\theta} \right] \]
given by
\[ \traj_s(\Time)=2 \sqrt{\Time(1-\Time)}\ s + \frac{\theta-1}{\sqrt{\theta}}\ \Time, \]
see \cref{fig:elipsy}.

\subsubsection{The main result 1: the trajectory of the second class particle}

\begin{theorem}    
    \label{thm:second-class-B}
We keep the assumptions and notations from \cref{sec:setup}.

Then there exists a sequence $(S_{i})$ of random variables  with 
the~property that the supremum distance 
\[  \sup_{\Time \in [0,1] } 
\left| \frac{1}{\sqrt{M_i N_i}}\ \usc_{i}\big( \left\lceil \Time M_i N_i \right\rceil \big) 
- \traj_{S_i}(\Time) \right|
\]
converges in probability to zero, as $i\to\infty$.

The probability distribution of $S_i$ is supported on the interval $[-1,1]$ and
for $i\to\infty$ it converges to the standard semicircular distribution with the
density
\begin{equation} 
\label{eq:SC}
f_{\operatorname{SC}}(x) = \frac{2}{\pi} \sqrt{1-x^2} \qquad \text{for $x\in [-1,1]$}.
\end{equation}
\end{theorem}

This theorem is illustrated on \cref{fig:elipsy}. 
Its proof is postponed to \cref{sec:particles}.
This result is analogous to the results of Ferrari and Kipnis \cite{Ferrari1995}, 
as well as Mountford and Guiol \cite{Mountford2005} for the usual TASEP system
 starting from a decreasing shock profile.

\begin{figure}
    \centering
    
    \subfile{figures/elipsy-graniczne.tex}

    \caption{The dashed lines are arcs of ellipses $\traj_s$ for the shape
    parameter $\theta=2$. 
		The shown values of the parameter $s\in
    F_{\operatorname{SC}}^{-1}\left\{ \nicefrac{0}{10},\nicefrac{1}{10},\dots,
    \nicefrac{10}{10} \right\}$ are the deciles of the semicircle
    distribution (above $F_{\operatorname{SC}}\colon[-1,1]\to[0,1]$ 
    denotes the cumulative distribution function of the semicircle
    distribution~\eqref{eq:SC}). 
		The shaded region forms \emph{the arctic ellipse}.
    The four black dots are the points where the arctic ellipse is tangent to
    the bounding box. The solid zigzag lines are trajectories of the second class
    particle for $N=500$ and $M=1000$. These trajectories were selected from a
    sample of $1000$ simulations; each of them corresponds to an appropriate
    empirical decile of the distribution of 
		the second class particle at time~$T=\frac{1}{2}$.
		}
    
    \label{fig:elipsy}
\end{figure}

\subsubsection{The limit trajectories}
\label{sec:limit-trajectories}

Each curve $\traj_s$ for the parameter \mbox{$s\in[-1,1]\setminus \{0\}$}
is an arc of an ellipse which fits into the bounding box $B_\theta$,
passes through the points
\begin{equation}
    \label{eq:points}
  (0,0), \qquad   \left( \frac{\theta-1}{\sqrt{\theta}} ,1 \right),
\end{equation}
and is tangent there to the bottom and the top edge of the bounding
box~$B_\theta$. In the degenerate case $s=0$ the curve $\traj_0$ is a straight
line connecting the aforementioned two points \eqref{eq:points}. The union of
the two extreme curves $\traj_{-1}$ and $\traj_1$ is the unique ellipse (which
we call \emph{the arctic ellipse}) which is inscribed into the bounding box
$B_\theta$ and is tangent to its bottom and its top side in the aforementioned
two points \eqref{eq:points}, see \cref{fig:elipsy}. In the special case
$\theta=1$ the scaling of the axes can be chosen in such a way that the arctic
ellipse becomes a circle which we call the \emph{``the arctic circle''}. Our use
of this name is not a coincidence since it turns out to be indeed related to the 
celebrated \emph{arctic circle theorem} \cite{Romik2012}.

\subsection{Random sorting networks}
\label{sec:random-sorting-networks}

\cref{thm:second-class-B} can be regarded as the solution to a toy version of
the problem of \emph{random sorting networks} considered by Angel, Holroyd,
Romik, and Vir\'{a}g \cite{Angel2007}.
More specifically, we consider the
symmetric group $\Sym_{N+M-1}$ viewed as the set of permutations of the set of
nodes \eqref{eq:set-of-nodes} and the permutation $\rho_{N,M}\in\Sym_{N+M-1}$
defined as
\[ \rho_{N,M}(i) = 
\begin{cases}
    i+M & \text{if } i\in\{-N+1,\dots,-1\},\\
    M-N & \text{if } i=0, \\
    i-N& \text{if } i\in\{1,\dots,M-1\}
\end{cases}
\]
which describes the change of the positions of the particles and holes during
the passage from the initial configuration shown on \cref{fig:scpA} to the final
configuration shown on \cref{fig:scpB}. For an integer $s\in\{-N+1,\dots,M-2\}$ denote
the adjacent transposition at location $s$ by 
$\tau_s=(s,s+1)\in\Sym_{N+M-1}$.

Any history of the particle system considered in \cref{sec:setup-scp} can be
encoded by the sequence $s_1,\dots,s_{t_{\max}-1}\in\{-N+1,\dots,M-2\}$, where
$s_t$ and $s_t+1$ are the nodes which are interchanged in $t$-th transition. It
is easy to check that
\begin{equation}
    \label{eq:factorization}
     \rho_{N,M} = \tau_{s_{t_{\max}-1}} \tau_{ s_{t_{\max}-2}}  \cdots \tau_{s_{2}}\ \tau_{s_{1}}
\end{equation}
and the corresponding sequence of partial products 
 \begin{equation}
     \label{eq:shortest-path}
     \id,\quad \tau_{s_{1}},\quad \tau_{s_{2}} \tau_{s_{1}},\quad  \dots,\quad \tau_{s_{t_{\max}-1}} \tau_{ s_{t_{\max}-2}}\cdots\tau_{s_{2}}\ \tau_{s_{1}} 
 \end{equation}
is a shortest path from the identity permutation $\id$ to $\rho_{N,M}$ in the
Cayley graph of the symmetric group $\Sym_{N+M-1}$ generated by adjacent transpositions.

Conversely, each shortest path \eqref{eq:shortest-path} in the Cayley graph
gives a valid history of the particle system. Any such a shortest path will be
called \emph{a sorting network}. The trajectory of the second class particle
\begin{multline}
    \label{eq:RSN} 
    \big( \usc(t) : t\in\{1,\dots,t_{\max} \} \big) =\\
     \big( 0,\quad \tau_{s_1}(0),\quad \tau_{s_2} \tau_{s_1}(0),
     \quad \dots,\quad \tau_{s_{t_{\max}-1}} \cdots \tau_{s_2} \tau_{s_1}(0) \big)
\end{multline}
corresponds in this language to the sequence of images of $0$ 
under the action of the entries of the sequence \eqref{eq:shortest-path}.

\smallskip

Angel, Holroyd, Romik, and Vir\'{a}g \cite{Angel2007} considered a more
difficult version of this setup in which the permutation $\rho_{N,M}$ is
replaced by the \emph{reverse permutation} 
$\rho\in\Sym_{N+M-1}$ given by 
\[ \rho(i)= M-N-i\]
and---among several other results---stated some conjectures concerning the
asymptotic behavior of the right-hand side of \eqref{eq:RSN} for a
\emph{random sorting network}, i.e., a random shortest path from the identity
permutation $\id$ to the reverse permutation $\rho$, sampled with the uniform
distribution. We focus today on \cite[Conjecture 1]{Angel2007} which is a
direct analogue of our \cref{thm:second-class-B}. A~minor difference is that
the limit curves which appear in \cref{thm:second-class-B} form a
one-parameter family of arcs of ellipses while the limit curves which appear in
\cite[Conjecture 1]{Angel2007} form a one-parameter family of \emph{sine
    curves}. (At first sight it might appear that the family of curves in
\cite{Angel2007} has two parameters, but one of these parameters can be
eliminated by the requirement about the positions of the endpoints.)

The aforementioned conjecture \cite[Conjecture 1]{Angel2007} was proved only very recently by Dauvergne and
Vir\'{a}g \cite{Dauvergne2020} who used methods quite different
from those which we use in the current paper.

\subsection{Sliding paths and evacuation paths in random tableaux}
\label{sec:tableaux}

As promised, we turn now to the second interest area of the current paper,
namely to random Young tableaux.

\subsubsection{Notations related to Young diagrams and tableaux} 
\label{sec:tableaux-begin}

We assume the reader’s basic knowledge of tableaux theory, including
partitions, (skew) Young diagrams, standard (skew) Young tableaux, RSK
algorithm, jeu de taquin, rectification, Littlewood--Richardson coefficients
and basics of the representation theory of the symmetric groups.

\medskip

We denote the set of partitions of~$n$ by~$\Y_n$.
We draw Young diagrams on the Cartesian plane
using the French convention, that is,~we draw them from the bottom to the top, see~\cref{fig:Youngdiagram}.

For any Young diagram~$\lambda$ we denote the set of standard Young tableaux
of shape~$\lambda$ by $\tableaux_\lambda$. \emph{The shape} of
a~tableau~$\tab$ will be denoted by $\sh(\tab)$ and its size by $|\sh(\tab)|$,
or shortly $|\tab|$. 

Let $\tab$ be a tableau. 
If $p$ is a number which appears exactly once in $\tab$ 
(which will always be the case in our considerations), 
we define the~\emph{position of the~box with the number $p$} as the
Cartesian coordinates of the bottom-left corner of the~unique square which
contains $p$; we denote this position by~$\pos_p(\tab)$. For example, for
$\tab$ from~\cref{fig:Youngtableau}, $\pos_5(\tab) = (2,0)$.

We will have a particular interest in Young diagrams and tableaux of square shape. 
By $\sq_N\in\Y_{N^2}$ we denote the square diagram with side of length $N$.

\begin{figure}\label{Youngdt}
    \centering
		\subfloat[]{
			\begin{tikzpicture}[scale=0.75]
                    
                    \fill[blue!5] (4,0) -- (4,1) -- (3,1) -- (3,2) -- (1,2) -- (1,4) -- (0,4) -- (0,0);
                    
                    \draw [->] (0,0) -- (5.8,0) node[anchor=north west] {$x$};
                    \draw [->] (0,0) -- (0,5.8) node[anchor=south east] {$y$};

                    \foreach \x in {1,...,5}
                    { \draw (\x,0 ) +(0,5pt) -- +(0,-5pt) node[anchor=north] {\small $ \x  $}; }
                    
                    \foreach \y in {1,...,5}
                    { \draw  (0, \y ) +(5pt,0) -- +(-5pt,0) node[anchor=east] {\small $ \y  $}; }

					\begin{scope}
					\clip[](4,0) -- (4,1) -- (3,1) -- (3,2) -- (1,2) -- (1,4) -- (0,4) -- (0,0);
					\draw[thick] (0,0) grid (10,10);
					\end{scope}
					\draw[ultra thick] (4,0) -- (4,1) -- (3,1) -- (3,2) -- (1,2) -- (1,4) -- (0,4) -- (0,0) -- cycle;
			\end{tikzpicture}  
	\label{fig:Youngdiagram} 
}
\hfill
\subfloat[]{
			\begin{tikzpicture}[scale=0.75]	
       
\fill[blue!5] (4,0) -- (4,1) -- (3,1) -- (3,2) -- (1,2) -- (1,4) -- (0,4) -- (0,0);

\draw [->] (0,0) -- (5.8,0) node[anchor=north west] {$x$};
\draw [->] (0,0) -- (0,5.8) node[anchor=south east] {$y$};

\foreach \x in {1,...,5}
{ \draw (\x,0 ) +(0,5pt) -- +(0,-5pt) node[anchor=north] {\small $ \x  $}; }

\foreach \y in {1,...,5}
{ \draw  (0, \y ) +(5pt,0) -- +(-5pt,0) node[anchor=east] {\small $ \y  $}; }

\begin{scope}
\clip[](4,0) -- (4,1) -- (3,1) -- (3,2) -- (1,2) -- (1,4) -- (0,4) -- (0,0);
\draw[thick] (0,0) grid (10,10);
\end{scope}
\draw[ultra thick] (4,0) -- (4,1) -- (3,1) -- (3,2) -- (1,2) -- (1,4) -- (0,4) -- (0,0) -- cycle;

					\draw (0.5,0.5) node {$1$};
					\draw (0.5,1.5) node {$3$};
					\draw (1.5,0.5) node {$2$};
					\draw (2.5,0.5) node {$5$};
					\draw (3.5,0.5) node {$9$};
					\draw (1.5,1.5) node {$4$};
					\draw (0.5,2.5) node {$6$};
					\draw (2.5,1.5) node {$8$};
					\draw (0.5,3.5) node {$7$};
		\end{tikzpicture}
\label{fig:Youngtableau} } 

\caption{\protect\subref{fig:Youngdiagram} The Young
diagram $\lambda = (4, 3, 1, 1)$ shown in the Cartesian coordinate
system. \protect\subref{fig:Youngtableau} Example of a standard Young tableau
of shape~$\lambda$. } 

\label{fig:Young}
\end{figure}

\subsubsection{Sliding paths and evacuation paths}
\label{sec:jdt-definition}

One of the operations heavily used in the~study of Young tableaux is \emph{\jdtincomplete} 
\cite[Section~1.2]{Fulton1997}, which acts on Young tableaux in the following way
(see~\cref{fig:jdtA,fig:jdtB}): we remove the bottom-left box of the given
tableau $\tab$ and obtain a~\emph{hole} in its place. Then we look at the two
boxes: the one to the~right and the one above the~hole, and choose the one
which contains the smaller number. We slide this smaller box into the location
of the hole, see \cref{fig:jdtZ}. As a result, the hole moves in the opposite
direction. We continue this operation as long as there is some box to the right or
above the hole. The path which was traversed by the `traveling hole' will be 
called the~\textbf{\emph{\jdtp}}, see~\cref{fig:jdtA}.
The result of \jdtincomplete
applied to a~tableau~$\tab$ will be denoted by~$j(\tab)$, 
see \cref{fig:jdtB}.

\newcommand{\letterX}{r}
\newcommand{\letterY}{s}
\begin{figure}
    \subfloat[]{
        \begin{tikzpicture}[scale=1.2]
        \clip (-0.3,-0.3) rectangle (2.3,2.3);
        \draw[black!20] (-1,-1) grid (3,3);
        \draw[pattern color=red!50,pattern=north east lines] (0,1) +(0.1,0.1) rectangle +(0.9,0.9);
        \draw[pattern color=green!50,pattern=north west lines] (1,0) +(0.1,0.1) rectangle +(0.9,0.9);
        \draw[fill=black!5] (1,1) +(0.1,0.1) rectangle +(0.9,0.9);
        \draw[fill=black!5] (1,2) +(0.1,0.1) rectangle +(0.9,0.9);
        \draw[fill=black!5] (2,1) +(0.1,0.1) rectangle +(0.9,0.9);
        \draw[fill=black!5] (2,0) +(0.1,0.1) rectangle +(0.9,0.9);
        \draw[fill=black!5] (0,2) +(0.1,0.1) rectangle +(0.9,0.9);
        \draw[fill=black!5] (-1,0) +(0.1,0.1) rectangle +(0.9,0.9);
        \draw[fill=black!5] (0,-1) +(0.1,0.1) rectangle +(0.9,0.9);
        \draw[fill=black!5] (-1,1) +(0.1,0.1) rectangle +(0.9,0.9);
        \draw[fill=black!5] (1,-1) +(0.1,0.1) rectangle +(0.9,0.9);
        \draw[fill=black!5] (-1,2) +(0.1,0.1) rectangle +(0.9,0.9);
        \draw[fill=black!5] (2,-1) +(0.1,0.1) rectangle +(0.9,0.9);
        \draw[fill=black!5] (-1,-1) +(0.1,0.1) rectangle +(0.9,0.9);
        \draw[fill=black!5] (2,2) +(0.1,0.1) rectangle +(0.9,0.9);
        \draw (0.5,1.5) node [circle,inner sep=2pt,fill=white] {$\letterX$};
        \draw (1.5,0.5) node [circle,inner sep=2pt,fill=white] {$\letterY$};
        \end{tikzpicture}
        \label{subsig:jdtA}}
    \hfill
    \subfloat[]{
        \begin{tikzpicture}[scale=1.2]
        \clip (-0.3,-0.3) rectangle (2.3,2.3);
        \draw[black!20] (-1,-1) grid (3,3);
        \draw[pattern color=red!50,pattern=north east lines] (0,0) +(0.1,0.1) rectangle +(0.9,0.9);
        \draw[pattern color=green!50,pattern=north west lines] (1,0) +(0.1,0.1) rectangle +(0.9,0.9);
        \draw[fill=black!5] (1,1) +(0.1,0.1) rectangle +(0.9,0.9);
        \draw[fill=black!5] (1,2) +(0.1,0.1) rectangle +(0.9,0.9);
        \draw[fill=black!5] (2,1) +(0.1,0.1) rectangle +(0.9,0.9);
        \draw[fill=black!5] (2,0) +(0.1,0.1) rectangle +(0.9,0.9);
        \draw[fill=black!5] (0,2) +(0.1,0.1) rectangle +(0.9,0.9);
        \draw[fill=black!5] (-1,0) +(0.1,0.1) rectangle +(0.9,0.9);
        \draw[fill=black!5] (0,-1) +(0.1,0.1) rectangle +(0.9,0.9);
        \draw[fill=black!5] (-1,1) +(0.1,0.1) rectangle +(0.9,0.9);
        \draw[fill=black!5] (1,-1) +(0.1,0.1) rectangle +(0.9,0.9);
        \draw[fill=black!5] (-1,2) +(0.1,0.1) rectangle +(0.9,0.9);
        \draw[fill=black!5] (2,-1) +(0.1,0.1) rectangle +(0.9,0.9);
        \draw[fill=black!5] (-1,-1) +(0.1,0.1) rectangle +(0.9,0.9);
        \draw[fill=black!5] (2,2) +(0.1,0.1) rectangle +(0.9,0.9);
        \draw (0.5,0.5) node [circle,inner sep=2pt,fill=white] {$\letterX$};
        \draw (1.5,0.5) node [circle,inner sep=2pt,fill=white] {$\letterY$};
        \draw[ultra thick,->] (0.5,1.5) -- (0.5,0.7);
        \end{tikzpicture}
        \label{subsig:jdtB}}
    \hfill
    \subfloat[]{
        \begin{tikzpicture}[scale=1.2]
        \clip (-0.3,-0.3) rectangle (2.3,2.3);
        \draw[black!20] (-1,-1) grid (3,3);
        \draw[pattern color=red!50,pattern=north east lines] (0,1) +(0.1,0.1) rectangle +(0.9,0.9);
        \draw[pattern color=green!50,pattern=north west lines] (0,0) +(0.1,0.1) rectangle +(0.9,0.9);
        \draw[fill=black!5] (1,1) +(0.1,0.1) rectangle +(0.9,0.9);
        \draw[fill=black!5] (1,2) +(0.1,0.1) rectangle +(0.9,0.9);
        \draw[fill=black!5] (2,1) +(0.1,0.1) rectangle +(0.9,0.9);
        \draw[fill=black!5] (2,0) +(0.1,0.1) rectangle +(0.9,0.9);
        \draw[fill=black!5] (0,2) +(0.1,0.1) rectangle +(0.9,0.9);
        \draw[fill=black!5] (-1,0) +(0.1,0.1) rectangle +(0.9,0.9);
        \draw[fill=black!5] (0,-1) +(0.1,0.1) rectangle +(0.9,0.9);
        \draw[fill=black!5] (-1,1) +(0.1,0.1) rectangle +(0.9,0.9);
        \draw[fill=black!5] (1,-1) +(0.1,0.1) rectangle +(0.9,0.9);
        \draw[fill=black!5] (-1,2) +(0.1,0.1) rectangle +(0.9,0.9);
        \draw[fill=black!5] (2,-1) +(0.1,0.1) rectangle +(0.9,0.9);
        \draw[fill=black!5] (-1,-1) +(0.1,0.1) rectangle +(0.9,0.9);
        \draw[fill=black!5] (2,2) +(0.1,0.1) rectangle +(0.9,0.9);
        \draw (0.5,1.5) node[circle,inner sep=2pt,fill=white] {$\letterX$};
        \draw (0.5,0.5) node[circle,inner sep=2pt,fill=white] {$\letterY$};
        \draw[ultra thick,->] (1.5,0.5) -- (0.7,0.5);
        \end{tikzpicture}
        \label{subsig:jdtC}}
    \caption{Elementary step of the \jdtincomplete transformation: 
        \protect\subref{subsig:jdtA} the initial configuration of boxes,
        \protect\subref{subsig:jdtB} the outcome of the slide in the case when $\letterX<\letterY$,
        \protect\subref{subsig:jdtC} the outcome of the slide in the case when $\letterY<\letterX$.
 Copyright \textcopyright2014 Society for Industrial and Applied Mathematics.  Reprinted from \cite{Sniady2014} with permission.  All rights reserved.    
}

    \label{fig:jdtZ}
\end{figure}

If $\tab$ is a standard tableau then $j(\tab)$ is no longer standard because the~numbering of its boxes
starts with $2$; however, if we decrease each entry of $j(\tab)$ by $1$ then it becomes standard. This observation allows us to define the \emph{\jdtt} 
$\jdtcomplete\colon \tableaux_\lambda\to \tableaux_\lambda$ which is a bijection on the set 
of standard Young tableaux of any fixed shape $\lambda$. 
The idea is to
put once again the~box with the~number~$|\tab|$ to the~aforementioned standardized version of
the tableau~$j(\tab)$ 
in the place where we removed a box during \jdtincomplete, see~\cref{fig:jdtC}.

\begin{figure}[t]
    \centering
    \subfloat[]{
			\begin{tikzpicture}[scale=0.75]
					\fill[blue!20] (0,0) rectangle +(1,1);
					\fill[blue!20] (0,1) rectangle +(1,1);
					\fill[blue!20] (1,1) rectangle +(1,1);
					\fill[blue!20] (2,1) rectangle +(1,1);
					\fill[blue!20] (2,2) rectangle +(1,1);
					
					\begin{scope}
					\clip[](5,0) -- (5,1) -- (4,1) -- (4,2) -- (3,2) -- (3,3) -- (1,3) -- (1,4) -- (0,4) -- (0,0);
					\draw (0,0) grid (10,10);
					\end{scope}
					\draw[ultra thick] (5,0) -- (5,1) -- (4,1) -- (4,2) -- (3,2) -- (3,3) -- (1,3) -- (1,4) -- (0,4) -- (0,0) -- cycle;
					\draw (0.5,0.5) node {$1$};
					\draw (0.5,1.5) node {$2$};
					\draw (1.5,0.5) node {$3$};
					\draw (2.5,0.5) node {$7$};
					\draw (3.5,0.5) node {$10$};
					\draw (4.5,0.5) node {$13$};
					\draw (1.5,1.5) node {$4$};
					\draw (0.5,2.5) node {$5$};
					\draw (2.5,1.5) node {$6$};
					\draw (1.5,2.5) node {$8$};
					\draw (0.5,3.5) node {$11$};
					\draw (2.5,2.5) node {$9$};
					\draw (3.5, 1.5) node {$12$};
			\end{tikzpicture}        
\label{fig:jdtA}
}
\hfill
\subfloat[]{
			\begin{tikzpicture}[scale=0.75]
					\fill[blue!20] (0,0) rectangle +(1,1);
					\fill[blue!20] (0,1) rectangle +(1,1);
					\fill[blue!20] (1,1) rectangle +(1,1);
					\fill[blue!20] (2,1) rectangle +(1,1);
					\fill[blue!10] (2,2) rectangle +(1,1);
                    
					\begin{scope}
					\clip[](5,0) -- (5,1) -- (4,1) -- (4,2) -- (2,2) -- (2,3) -- (1,3) -- (1,4) -- (0,4) -- (0,0);
					\draw (0,0) grid (10,10);
					\end{scope}
					\draw[ultra thick] (5,0) -- (5,1) -- (4,1) -- (4,2) -- (2,2) -- (2,3) -- (1,3) -- (1,4) -- (0,4) -- (0,0) -- cycle;
					\draw (0.5,0.5) node {$2$};
					\draw (0.5,1.5) node {$4$};
					\draw (1.5,0.5) node {$3$};
					\draw (2.5,0.5) node {$7$};
					\draw (3.5,0.5) node {$10$};
					\draw (4.5,0.5) node {$13$};
					\draw (1.5,1.5) node {$6$};
					\draw (0.5,2.5) node {$5$};
					\draw (2.5,1.5) node {$9$};
					\draw (1.5,2.5) node {$8$};
					\draw (0.5,3.5) node {$11$};
					\draw (3.5, 1.5) node {$12$};
			\end{tikzpicture}
\label{fig:jdtB}
}        
\hfill
\subfloat[]{ 
	\begin{tikzpicture}[scale=0.75]
	\fill[blue!20] (0,0) rectangle +(1,1);
	\fill[blue!20] (0,1) rectangle +(1,1);
	\fill[blue!20] (1,1) rectangle +(1,1);
	\fill[blue!20] (2,1) rectangle +(1,1);
    \fill[blue!20] (2,2) rectangle +(1,1);
    
	\fill[pattern=north west lines, pattern color=blue] (2,2) rectangle +(1,1);

	\begin{scope}
	\clip[](5,0) -- (5,1) -- (4,1) -- (4,2) -- (3,2) -- (3,3) -- (1,3) -- (1,4) -- (0,4) -- (0,0);
	\draw (0,0) grid (10,10);
	\end{scope}
	\draw[ultra thick] (5,0) -- (5,1) -- (4,1) -- (4,2) -- (3,2) -- (3,3) -- (1,3) -- (1,4) -- (0,4) -- (0,0) -- cycle;
	\draw (0.5,0.5) node {$1$};
	\draw (0.5,1.5) node {$3$};
	\draw (1.5,0.5) node {$2$};
	\draw (2.5,0.5) node {$6$};
	\draw (3.5,0.5) node {$9$};
	\draw (4.5,0.5) node {$12$};
	\draw (1.5,1.5) node {$5$};
	\draw (0.5,2.5) node {$4$};
	\draw (2.5,1.5) node {$8$};
	\draw (1.5,2.5) node {$7$};
	\draw (0.5,3.5) node {$10$};
	\draw (3.5, 1.5) node {$11$};
	\draw (2.5, 2.5) node[fill=blue!20,rounded corners=2pt,inner sep=1pt] {$13$};
	\end{tikzpicture}        
	\label{fig:jdtC} 
} 
	\caption{\protect\subref{fig:jdtA} A standard Young tableau
$\tab$ of shape $\lambda = (5,4,2,1)$. The highlighted boxes form the~\emph{\jdtp}.
\protect\subref{fig:jdtB} The outcome $j(\tab)$ of the \jdtincomplete
transformation. \protect\subref{fig:jdtC} The result $\jdtcomplete(\tab)$ of the \emph{\jdtt} applied to $\tab$.}
\label{fig:jdt}
\end{figure}

For a given tableau $\tab\in \tableaux_\lambda$ with $n=|\lambda|$ boxes the \jdtincomplete transformation $j$ can be 
iterated $n$ times until we end with the empty tableau. 
During each iteration the box with the biggest number $n$ either moves one node left or down, or stays put. 
Its trajectory 
\begin{equation}
\label{eq:evacuationtrajectory}
\evac(T)=
\Big( 
\pos_n(\tab),\ 
\pos_n\!\big( j(\tab) \big),
\ \dots,\
\pos_n\! \big( j^{n-1}(\tab) \big) 
\Big) 
\end{equation}
will be called the~\textbf{\emph{evacuation path}}.

\subsubsection{The main results 2 and 3: asymptotics of sliding paths and evacuation paths}

Observe that if we draw the boxes of a given square tableau $\tab\in
\tableaux_{\sq_N}$ as little squares of size $\frac{1}{N}$ then the
corresponding \jdtp is a zigzag line connecting the opposite corners of the
unit square $[0,1]^2$. \emph{Let $\tab\in \tableaux_{\sq_N}$ be a random
    standard Young tableau of square shape (sampled with the uniform probability
    distribution on $\tableaux_{\sq_N}$ which will be denoted $\Pp$).  Our goal is to find asymptotics
    of such random zigzag lines in the limit as $N\to\infty$, see
    \cref{fig:simulation-jdtZ}.} We will show that there exists a family of smooth
lines, called \emph{meridians}, which connect the opposite corners of the unit
square, with the property that the probability distribution of \emph{the
    scaled \jdtp for a random tableau converges, as~$N\to\infty$, to a random
    meridian}.

\medskip

An analogous result holds true for the scaled evacuation path for a random
square tableau: during iteratively applied \jdtincomplete operations~$j$, the
biggest box of the tableau asymptotically moves along a random meridian. 

A version of this result applies also to the other boxes of the tableau;
it follows that the time evolution of the tableau in the iterated applications of \jdtincomplete
\begin{equation}
\label{eq:jdt-tableaux-sequence}
T,\; j(T),\; j^2(T),\; \ldots, j^{N^2}(T) 
\end{equation}
converges in probability, as $N\to\infty$, to dynamics of an incompressible liquid which flows along the meridians.

\medskip

For the details of our results, see
\cref{thm:evacuation,thm:jdt} in \cref{sec:problems}.

\begin{figure}
    \centering

        \subfile{figures/kwadrat-jdt-duzo-100.tex}

    \caption{The nine zigzag lines are sample \jdtps for random square tableaux
    of size $N=100$, selected so that they cross the anti-diagonal near the
    corresponding meridians (smooth thick curves) with the longitudes
    $\longitude\in\{\nicefrac{1}{10},\nicefrac{2}{10},\dots,\nicefrac{9}{10}\}$. The gray lines are the meridians with the longitudes \mbox{$\longitude\in\{\nicefrac{2}{100},\nicefrac{4}{100},\dots,\nicefrac{98}{100}\}$}. See also the blue-to-red family of curves on \cref{fig:geographic}.}

\label{fig:simulation-jdtZ}
    
\end{figure}

\subsubsection{Not only squares}
\label{sec:tableaux-end}
\label{sec:no-only-squares}

For simplicity and concreteness we stated our main results
concerning random Young tableaux only for
large random Young tableaux of \emph{square shape}.
However, analogous results hold true also for random tableaux
of shape which is a \emph{balanced Young diagram} 
(see~\cref{sec:definitions} for the~definition 
and \cref{fig:L-shape} for a teaser).
In \cref{sec:generalization} we present a~way
in which the results obtained in this paper
can be used (or generalized) 
in order to cover the~class of balanced Young diagrams.

\begin{figure}[t]
	\centering
		\subfile{figures/L-jdt-duzo.tex}	
		\caption{
		Sample \jdtps in random Young tableaux of an $L$-shape with $3600$ boxes.}
\label{fig:L-shape}
\end{figure}

\subsection{The content of the paper}

The paper is organized as follows. 

In \cref{sec:problems} we state the main results (\cref{thm:evacuation,thm:jdt})
about the typical shapes of the evacuation paths and sliding paths 
in square Young tableaux. 

In \cref{sec:definitions} we give basic definitions on permutations and representation theory.

In \cref{sec:longitude-begins} we introduce a `surfers' language which we will use to 
describe dynamics of the box with the biggest entry (which we will call `the surfer')
and the smaller boxes (`the water'). 
In this spirit we also introduce the story of the \emph{multisurfers} which will play a crucial role 
in our proofs and considerations. 
We will use this new multisurfer story later as a point of reference for the original problem of
the (single) surfer in order to prove \cref{thm:all-the-same}
concerning the position of the surfer along its journey. 
We sketch the proof in \cref{sec:plan}. 

In \cref{sec:single-multi} we show the way in which we will embed simultaneously 
both the story of the single surfer and the story of the multisurfers into a common universe. 

In \cref{sec:udist} we provide \cref{prop:longitude-experimental} concerning the distribution
of the multisurfers on the water. 
We use here the Jucys--Murphy elements to give a direct link between the statistical properties 
of the multisurfers and 
the symmetric group characters evaluated on certain polynomials in the Jucys--Murphy elements.

In \cref{sec:longitude-ends} we use 
the aforementioned results from \cref{sec:single-multi,sec:udist} 
to prove \cref{thm:all-the-same}.

In \cref{sec:proof-of-evacuation-for-squares} we prove \cref{thm:evacuation} 
concerning typical evacuation paths. 

\cref{sec:equivalence-of-problems} is devoted to the proof of \cref{thm:equivalence}
which shows the equivalence  
between the problems of finding the sliding paths and 
the evacuation paths in random Young tableaux of given shape.

In \cref{sec:generalization} we extend our 
main results (\cref{thm:extension}) 
concerning typical evacuation and sliding paths
to some subset of $C$-balanced Young tableaux. 

In \cref{sec:particles} we provide the link between the trajectory of the second class particle
in an interacting particle system and the sliding path 
for a random Young tableau and prove \cref{thm:second-class-B}. 

\section{The limit shape of \jdtps and evacuation paths}
\label{sec:mainresults}
\label{sec:problems}

\subsection{Asymptotics of a single box in the evacuation trajectory}
\label{sec:single-box}
As we mentioned in~\cref{sec:tableaux},  
we will focus on random standard Young tableaux 
of square shape~$\sq_N$ sampled according to 
the~uniform measure~$\Pp$. 
The symbol $\tab_N$ will be reserved for such a~uniformly~random
square Young tableau of shape~$\sq_N$.

The position of each of the boxes in the evacuation
path~\eqref{eq:evacuationtrajectory} coincides with the~position of a specific
box in the standard Young tableau obtained by iterating the \jdtt $\jdtcomplete$:
\begin{equation}
    \label{eq:happy-box} 
    \pos_{N^2}\big( j^{i}(\tab_N)\big)  = \pos_{N^2-i} \big( {\jdtcomplete}{}^{i}(\tab_N)\big).
\end{equation}
Since $\jdtcomplete\colon \tableaux_{\sq_N}\to \tableaux_{\sq_N}$ is a bijection, for each
$i \geq 0$ the latter standard tableaux $\jdtcomplete{}^i(\tab_N)$ is also a
uniformly random square Young tableau. It follows that the solution to the
(much simpler) problem of understanding the asymptotics of a \emph{single}
element of the evacuation trajectory~\eqref{eq:evacuationtrajectory} follows
from the work of Pittel and Romik \cite{Pittel2007},
see also \cref{sec:random-position} below.
In the current section we will recall the details of their work 
and we will use it to state our second main
result, \cref{thm:evacuation} (which addresses the more complex problem of
understanding the \emph{whole} evacuation
trajectory~\eqref{eq:evacuationtrajectory}) and the third main result, \cref{thm:jdt}.

\subsection{The circles of latitude $g_\alpha$}
\label{sec:levelcurves}

The~\emph{Russian coordinate system} is given by the following transformation of the
Cartesian plane 
(warning: our notations differ from those of Pittel and Romik \cite{Pittel2007}
who scale the coordinates below by an additional
factor $1/\sqrt{2}$):
\[
u := x-y, \qquad v := x+y,
\]
see~\cref{fig:Cart-Russian}.

\begin{figure}[t]
	\begin{tikzpicture}[scale=1.2]
		\clip (-1.2,-2) rectangle (5,5); 

		\draw [thick,blue,->] (0,0) -- (4.8,0) node[anchor=north] {\textcolor{blue}{$x$}};
		\draw [thick,blue,->] (0,0) -- (0,4.8) node[anchor=east] {\textcolor{blue}{$y$}};

		\foreach \x in {1,...,4}
		{ \draw[thick,blue] (\x,0 ) +(0,5pt) -- +(0,-5pt) node[anchor=north] {\small $ \x  $}; }

		\foreach \y in {1,...,4}
		{ \draw[thick,blue]  (0, \y ) +(5pt,0) -- +(-5pt,0) node[anchor=east] {\small $ \y  $}; }

		\foreach \x/\xx in {-1,0/{},1,2,3}
		{ \draw[thick,red] (0.5*\x,-0.5*\x) +(2pt,2pt) -- +(-2pt,-2pt) node[anchor=north east] {\small $ \xx  $}; }

		\fill[blue!5] (4,0) -- (4,1) -- (3,1) -- (3,2) -- (1,2) -- (1,4) -- (0,4) -- (0,0);

		\begin{scope}
				\clip[](4,0) -- (4,1) -- (3,1) -- (3,2) -- (1,2) -- (1,4) -- (0,4) -- (0,0);
				\draw[thick] (0,0) grid (10,10);
		\end{scope}
		\draw[ultra thick] (4,0) -- (4,1) -- (3,1) -- (3,2) -- (1,2) -- (1,4) -- (0,4) -- (0,0) -- cycle;

		\foreach \x/\xx in {-4,-3,-2,-1,0,1,2,3}
		{  \draw[red, dotted]  (0.5*\x,-0.5*\x) -- +(5,5); }

    	\foreach \x/\xx in {1,...,5}
    {  \draw[red, dotted]  (0.5*\x,0.5*\x) +(-5,5) -- +(5,-5); }
    
    		\foreach \x/\xx in {1,...,5}
    { \draw[thick,red] (0.5*\x,0.5*\x) +(-2pt,2pt) -- +(2pt,-2pt) node[anchor=north west] {\small $    \xx  $}; }

        \draw [thick,red,->] (-1,1) -- (2,-2) node[anchor=south west] {$u=x-y$};
        \draw [thick,red,->] (-0.5,-0.5) -- (3,3) node[anchor=south west,fill=white,rounded corners=2pt,inner sep=1pt] {$v=x+y$};

	\end{tikzpicture}
		 \caption{The Cartesian and the Russian coordinate systems on the plane.   
			}
		\label{fig:Cart-Russian}
\end{figure}

For each $0 \leq \alpha \leq 1$ and $u \in \left[ -2\sqrt{\alpha(1-\alpha)}, 2\sqrt{\alpha(1-\alpha)} \right]$ define 
\[k_{\alpha, u} := \sqrt{4\alpha(1-\alpha) - u^2}\]
and for any $0 < \alpha < 1$ the function
\[
h_\alpha \colon \left[ -2\sqrt{\alpha(1-\alpha)}, 2\sqrt{\alpha(1-\alpha)} \right] \to \R
\]
given by  
\begin{equation}\label{height}
h_\alpha(u) := \left\{ \begin{array}{ll} \frac{2 u}{\pi} \arctan\left(\frac{1-2\alpha}{k_{\alpha, u}} \cdot u \right) + 
\frac{2 }{\pi} \arctan \left(\frac{k_{\alpha, u}}{1-2\alpha}\right) & \text{ if \ } 0 < \alpha < \frac{1}{2}, \\
2  - \frac{2 u}{\pi} \arctan\left(\frac{2\alpha-1}{k_{\alpha, u}} \cdot u \right) - 
\frac{2}{\pi} \arctan \left(\frac{k_{\alpha, u}}{2\alpha-1}\right) & \text{ if \ } \frac{1}{2} < \alpha < 1 , \\
1 & \text{ if \ } \alpha = \frac{1}{2}. \end{array} \right.
\end{equation}
In the expression above there may occur $k_{\alpha, u} = 0$ 
in the denominator; 
in such case (for fixed $\alpha \in (0,1)$)
we consider the appropriate limit:
if $u_0 = \pm 2\sqrt{\alpha(1-\alpha)}$ then 
\[
h_\alpha(u_0) := 
\lim_{u \to u_0} h_{\alpha}(u) 
= \begin{cases}
u_0 & \text{for $0 < \alpha < \frac{1}{2}$}, \\
2-u_0 & \text{for $\frac{1}{2} < \alpha < 1$}.
\end{cases}
\]
Additionally we define the one-point functions $h_0(0) := 0$ and $h_1(0) := 2$.

For any $\alpha \in [0,1]$ we consider 
the curve which in the Russian coordinate system is defined as
\begin{equation}
\label{eq:alpha-level-curve}
g_\alpha^{\text{Rus}} := \left\{\big(u, h_\alpha(u) \big) : |u| \leq 2\sqrt{\alpha(1-\alpha)}  \right\} \subset \R^2 
\end{equation}
and, equivalently, in the Cartesian coordinates is given by (see~\cref{fig:level-curves,fig:geographic})
\[g_\alpha := \left\{\left(\frac{u+h_\alpha(u)}{2}, \frac{h_\alpha(u) - u}{2}\right) :
|u| \leq 2\sqrt{\alpha(1-\alpha)} \right\} \subset  [0,1]^2 . \]
We call $g_\alpha$ \textbf{\emph{the circle of latitude}} with the latitude $\alpha$.

\begin{figure}
    \centering
    \begin{tikzpicture}[scale=1]
        \draw[ultra thick,->] (0,0) -- (11,0) node[anchor=north west] {$x$};
        \draw[ultra thick,->] (0,0) -- (0,11) node[anchor=south east] {$y$};
        { \draw[ultra thick] (10,0 ) +(0,5pt) -- +(0,-5pt) node[anchor=north] {\small $ 1 $}; }
        { \draw[ultra thick]  (0,10 ) +(5pt,0) -- +(-5pt,0) node[anchor=east] {\small $1$}; }
        
        \begin{scope}
            \clip[](10,0) -- (10,10) -- (0,10) -- (0,0);
            \draw[black!20] (0,0) grid (10,10);
            
            \definecolor{darkcyan}{rgb}{0.0, 0.55, 0.55}

            \draw[darkcyan!25!orange,ultra thick] plot[smooth,scale=10] file {circle025.txt};
            \draw[darkcyan!50!orange,ultra thick] (0,10) -- (10,0);
            \draw[darkcyan!75!orange,ultra thick] plot[smooth,scale=10] file {circle075.txt};
            
            \draw[line width=1mm,darkcyan!75!orange,dashed,rounded corners](10,4) -- (8,4) -- (8,5) -- (7,5) -- (7,7) -- (6,7) -- (6,8) -- (5,8) -- (5,9) -- (2,9) -- (2,10); 
            \draw[line width=1mm,darkcyan!50!orange,dashed,rounded corners](10,2) -- (8,2) -- (8,3) -- (7,3) -- (7,4) -- (5,4) -- (5,5) -- (4,5) -- (4,6) -- (2,6) -- (2,8) -- (1,8) -- (1,10); 
            \draw[line width=1mm,darkcyan!25!orange,dashed,rounded corners](9,0) -- (9,1) -- (6,1) -- (6,2) -- (3,2) -- (3,4) -- (2,4) -- (2,5) -- (1,5) -- (1,7) -- (0,7); 
            
            \draw (0.5,0.5) node {$1$};
            \draw (1.5,0.5) node {$2$};
            \draw (2.5,0.5) node {$3$};
            \draw (0.5,1.5) node {$4$};
            \draw (1.5,1.5) node {$5$};
            \draw (3.5,0.5) node {$6$};
            \draw (0.5,2.5) node {$7$};
            \draw (4.5,0.5) node {$8$};
            \draw (5.5,0.5) node {$9$};
            \draw (0.5,3.5) node {$10$};
            \draw (1.5,2.5) node {$11$};
            \draw (2.5,1.5) node {$12$};
            \draw (2.5,2.5) node {$13$};
            \draw (1.5,3.5) node {$14$};
            \draw (2.5,3.5) node {$15$};
            \draw (3.5,1.5) node {$16$};
            \draw (0.5,4.5) node {$17$};
            \draw (4.5,1.5) node {$18$};
            \draw (5.5,1.5) node {$19$};
            \draw (1.5,4.5) node {$20$};
            \draw (6.5,0.5) node {$21$};
            \draw (7.5,0.5) node {$22$};
            \draw (0.5,5.5) node {$23$};
            \draw (8.5,0.5) node {$24$};
            \draw (0.5,6.5) node {$25$};
            \draw (0.5,7.5) node {$26$};
            \draw (3.5,2.5) node {$27$};
            \draw (4.5,2.5) node {$28$};
            \draw (0.5,8.5) node {$29$};
            \draw (3.5,3.5) node {$30$};
            \draw (6.5,1.5) node {$31$};
            \draw (2.5,4.5) node {$32$};
            \draw (4.5,3.5) node {$33$};
            \draw (5.5,2.5) node {$34$};
            \draw (6.5,2.5) node {$35$};
            \draw (1.5,5.5) node {$36$};
            \draw (1.5,6.5) node {$37$};
            \draw (3.5,4.5) node {$38$};
            \draw (2.5,5.5) node {$39$};
            \draw (3.5,5.5) node {$40$};
            \draw (5.5,3.5) node {$41$};
            \draw (1.5,7.5) node {$42$};
            \draw (7.5,1.5) node {$43$};
            \draw (0.5,9.5) node {$44$};
            \draw (8.5,1.5) node {$45$};
            \draw (9.5,0.5) node {$46$};
            \draw (6.5,3.5) node {$47$};
            \draw (4.5,4.5) node {$48$};
            \draw (7.5,2.5) node {$49$};
            \draw (9.5,1.5) node {$50$};
            \draw (4.5,5.5) node {$51$};
            \draw (8.5,2.5) node {$52$};
            \draw (2.5,6.5) node {$53$};
            \draw (9.5,2.5) node {$54$};
            \draw (5.5,4.5) node {$55$};
            \draw (2.5,7.5) node {$56$};
            \draw (3.5,6.5) node {$57$};
            \draw (7.5,3.5) node {$58$};
            \draw (5.5,5.5) node {$59$};
            \draw (4.5,6.5) node {$60$};
            \draw (6.5,4.5) node {$61$};
            \draw (8.5,3.5) node {$62$};
            \draw (1.5,8.5) node {$63$};
            \draw (5.5,6.5) node {$64$};
            \draw (3.5,7.5) node {$65$};
            \draw (4.5,7.5) node {$66$};
            \draw (2.5,8.5) node {$67$};
            \draw (3.5,8.5) node {$68$};
            \draw (1.5,9.5) node {$69$};
            \draw (4.5,8.5) node {$70$};
            \draw (5.5,7.5) node {$71$};
            \draw (6.5,5.5) node {$72$};
            \draw (9.5,3.5) node {$73$};
            \draw (6.5,6.5) node {$74$};
            \draw (7.5,4.5) node {$75$};
            \draw (7.5,5.5) node {$76$};
            \draw (8.5,4.5) node {$77$};
            \draw (7.5,6.5) node {$78$};
            \draw (8.5,5.5) node {$79$};
            \draw (5.5,8.5) node {$80$};
            \draw (6.5,7.5) node {$81$};
            \draw (2.5,9.5) node {$82$};
            \draw (9.5,4.5) node {$83$};
            \draw (8.5,6.5) node {$84$};
            \draw (7.5,7.5) node {$85$};
            \draw (6.5,8.5) node {$86$};
            \draw (3.5,9.5) node {$87$};
            \draw (9.5,5.5) node {$88$};
            \draw (9.5,6.5) node {$89$};
            \draw (4.5,9.5) node {$90$};
            \draw (7.5,8.5) node {$91$};
            \draw (8.5,7.5) node {$92$};
            \draw (5.5,9.5) node {$93$};
            \draw (9.5,7.5) node {$94$};
            \draw (8.5,8.5) node {$95$};
            \draw (6.5,9.5) node {$96$};
            \draw (7.5,9.5) node {$97$};
            \draw (8.5,9.5) node {$98$};
            \draw (9.5,8.5) node {$99$};
            \draw (9.5,9.5) node {$100$};
        \end{scope}
        
        \draw[thick](0,0) -- (0,10) -- (10,10) -- (10,0) -- cycle ;
    \end{tikzpicture}
    \caption{A scaled down sample random square tableau of size $N=10$. The zigzag
        lines are the level curves for $\alpha\in\{ \nicefrac{1}{4} ,\nicefrac{2}{4},
        \nicefrac{3}{4} \}$. The smooth lines are the corresponding circles
        of~latitude~$g_\alpha$, see also the orange-to-green family of curves on
        \cref{fig:geographic}.} 
    
    \label{fig:level-curves}
\end{figure}

\medskip

Roughly speaking, for each $\alpha\in[0,1]$ the (scaled down)
\emph{$\alpha$-level curve} (which separates the boxes with entries $\leq
\alpha N^2$ from the boxes bigger than this threshold, see
\cref{fig:level-curves}) in a uniformly random square tableau $\tab_N$ 
converges in probability, as $N\to\infty$, to the circle of latitude $g_\alpha$, see
\cite[Theorem~1]{Pittel2007} 
for a precise statement.

\subsection{The random position of the box $\lfloor \alpha N^2 \rfloor$, 
the limit measure $\cotmeas_{\alpha}$}
\label{sec:random-position}

Pittel and Romik \cite[Theorem~2]{Pittel2007} also found the~explicit formula
for the~limit distribution of the~scaled down location 
\[ \frac{1}{N} \pos_{\lfloor \alpha N^2
    \rfloor} \left( \tab_N \right)
\]
of the~entry~$\lfloor \alpha N^2 \rfloor$ in a uniformly random square Young
tableau $T_N$, as \mbox{$N\to\infty$}. This limit distribution turns out to be
supported on the circle of latitude~$g_\alpha$ and thus it is uniquely determined by the
probability distribution of the~\mbox{$u$-coordinate}. The latter, denoted by
$\cotmeas_\alpha$, turns out to be the semicircular distribution on the interval
$\left[ -2\sqrt{\alpha(1-\alpha)}, 2\sqrt{\alpha(1-\alpha)} \right]$ with the
density
\begin{equation}\label{eq:density} 
f_{\cotmeas_\alpha}(u) := \frac{k_{\alpha, u}}{2 \pi \alpha (1-\alpha)} .
\end{equation}

\subsection{Geographic coordinates on the square}
\label{sec:geographic}

For any point $p = (x,y) \in \mysquare$ of the unit square there is exactly
one $\alpha=\alpha(p) \in [0,1]$ such that $p$ lies on the curve $g_\alpha$.
We say that \emph{the~latitude of $p$} is equal to~$\alpha$. With the
notations of Pittel and Romik \cite{Pittel2007} the latitude
$\alpha(x,y)=L(x,y)$ is just the limit height function of random square
standard Young tableaux.

The \emph{longitude of $p$}, 
which we denote by 
\[ 
\longitude(p)
=\cotmeas_\alpha \! \left( \left( -\infty, x-y \right] \right) 
= F_{\cotmeas_\alpha}(x-y), 
\]
is defined as the~mass of the~points on the curve $g_\alpha$ which have their
\mbox{$u$-coordinate} not greater than the~\mbox{$u$-coordinate} of the~point~$p$ or, equivalently, in terms of the~cumulative
distribution function $F_{\cotmeas_\alpha}$ of the~measure~$\cotmeas_\alpha$. 
Notice that $\psi((0,0)) = \psi((1,1)) = 1$.
The set of points of the unit square $\mysquare$ with equal longitude~$\longitude$ is a~curve called
the~\emph{meridian~$\longitude$}, see~\cref{fig:simulation-jdtZ,fig:geographic}.

For given $\alpha \in (0,1)$ and $\longitude\in [0,1]$ 
we denote by 
\[
\point_{\alpha, \longitude} = \left(x_\alpha^\longitude, y_\alpha^\longitude \right) \in \mysquare
\]
the unique point of~the~unit square~$\mysquare$ with the~latitude~$\alpha$
and the longitude~$\longitude$. 
We set additionally $P_{0, \psi} = (0,0)$ and $P_{1, \psi} = (1,1)$
for any $\psi \in [0,1]$. 
We~denote by
\[
u_{\alpha}^{\longitude} := x_\alpha^\longitude - y_\alpha^\longitude
\quad \text{ and } \quad
v_{\alpha}^{\longitude} := x_\alpha^\longitude + y_\alpha^\longitude
\] 
the~$u$- and $v$-coordinate of the point 
$\left(x_\alpha^\longitude, y_\alpha^\longitude \right) = \point_{\alpha, \longitude}$.

\begin{remark}
    We will not use the following observation, but fans of cartography may find
it interesting: for reasons which will hopefully become obvious later on,
the map
    \[ [0,1]^2 \ni (\alpha, \longitude) \mapsto \point_{\alpha, \longitude} \in [0,1]^2 \]
    is \emph{equiareal} which manifests  by equality of the areas of the curvilinear
    rectangles on \cref{fig:geographic}.
\end{remark}

\begin{remark}
The result of Pittel and Romik \cite[Theorem~2]{Pittel2007} 
is a special case of a general phenomenon
of existence of the level curves (circles of latitudes)
for random Young tableaux of specified shape.
Biane \cite{Biane1998} proved that such level curves exist 
for any balanced sequence of Young diagrams, 
see \cref{sec:generalization} 
for more information. 
\end{remark}

\begin{figure}[t]
    \centering 
     
\centering
\begin{tikzpicture}[scale=0.1]

\draw[ultra thick,->] (0,0) -- (110,0) node[anchor=north west] {$x$};
\draw[ultra thick,->] (0,0) -- (0,110) node[anchor=south east] {$y$};
{ \draw[ultra thick] (100,0 ) +(0,50pt) -- +(0,-50pt) node[anchor=north] {\small $ 1 $}; }
{ \draw[ultra thick]  (0,100 ) +(50pt,0) -- +(-50pt,0) node[anchor=east] {\small $1$}; }

\foreach \x in {2,4,...,8}    { \draw[black!25] plot[smooth,scale=100] file {rownolezniki/rownoleznik-0.0\x.txt}; }
\foreach \x in {12,14,...,98}    { \draw[black!25] plot[smooth,scale=100] file {rownolezniki/rownoleznik-0.\x.txt}; }

\foreach \x in {2,4,...,8}    { \draw[black!25] plot[smooth,scale=100] file {poludniki/poludnik-0.0\x.txt}; }
\foreach \x in {12,14,...,98} { \draw[black!25] plot[smooth,scale=100] file {poludniki/poludnik-0.\x.txt}; }

\foreach \x in {10,20,...,90}
{
    \draw[red!\x!blue,thick] plot[smooth,scale=100] file {poludniki/poludnik-0.\x.txt};
}

\definecolor{darkgreen}{rgb}{0.0, 0.2, 0.13}
\definecolor{darkcyan}{rgb}{0.0, 0.55, 0.55}

\foreach \x in {10,20,...,90} { \draw[darkcyan!\x!orange,thick] plot[smooth,scale=100] file {rownolezniki/rownoleznik-0.\x.txt}; }

\draw[thick](0,0) -- (0,100) -- (100,100) -- (100,0) -- cycle ;

\end{tikzpicture} 
    
    \caption{ The geographic coordinate system on the unit square $[0,1]^2$.
    The blue-to-red family of thick colored curves connecting the bottom left and
    the upper right corner are the meridians with the longitudes
    $\longitude\in\{\nicefrac{1}{10},\nicefrac{2}{10},\dots,\nicefrac{9}{10}\}$. 
    The gray lines are the meridians with the longitudes 
    \mbox{$\longitude\in\{\nicefrac{2}{100},\nicefrac{4}{100},\dots,\nicefrac{98}{100}\}$}, cf.~\cref{fig:simulation-jdtZ}. 
    The orange-to-green family of thick colored curves are the~circles of
    latitudes $g_\alpha$ with the latitudes $\alpha\in\{
    \nicefrac{1}{10},\nicefrac{2}{10},\dots,\nicefrac{9}{10}\}$. The thin,
    gray lines are the~circles of latitude $g_\alpha$ with
    the latitudes
    $\alpha\in\{\nicefrac{2}{100},\nicefrac{4}{100},\dots,\nicefrac{98}{100}\}$, see \cref{fig:level-curves}. 
    The shown meridians and circles of latitude split the square into a $50\times 50$ grid of curvilinear rectangles with equal areas.}
    
    \label{fig:geographic}
\end{figure}

\subsection{The second main result. Typical evacuation path}
\label{subsec:evacuation-Xt-definition}

For a given tableau $\tab_N\in \tableaux_{\sq_N}$ and $\tim\in[0,1)$ we denote
by
\begin{equation}
    \label{eq:evacuation}
\mypoint_\tim= \mypoint_\tim(\tab_N) 
= \frac{1}{N} \pos_{N^2}\left( j^{\lfloor \tim N^2 \rfloor}(\tab_N) \right) 
\in \mysquare
\end{equation}
the scaled down position of the~point from the evacuation path
$\evac(\tab_N)$, cf.~\eqref{eq:evacuationtrajectory}. Clearly, the parameter
$t$ indicates how many boxes
were removed so far and therefore we can relate
to it as the~\emph{time}.

\medskip

Our second main result states that, asymptotically, \emph{the scaled
    evacuation path in a~random square tableau is a random meridian},
see~\cref{fig:simulation-jdtZ}.

\begin{restatable}{theorem}{thm:evacuation}
    \label{thm:evacuation} 
    
For each $N\in \N$ there exists a random variable~$\Longitude_N\colon
\tableaux_{\sq_N} \to [0,1]$ such that the supremum distance
\begin{equation}
\label{eq:main-theorem}
\sup_{\tim \in [0,1]} 
			\left| \mypoint_\tim(\tab_N) -
			\point_{1-\tim, \Longitude_N(\tab_N)}  \right|
\end{equation}
converges in probability to zero, as $N\to\infty$. More explicitly: for each $\varepsilon>0$
\[ 
\lim_{N \to \infty} \Pp 
\Big\{ \tab_N \in \tableaux_{\sq_N}: 
\sup_{\tim \in [0,1]} 
\left| \mypoint_\tim(\tab_N) -
\point_{1-\tim, \Longitude_N(\tab_N)}  \right|
             > \varepsilon \Big\} = 0.
\]

The probability distribution
of the random variable~$\Longitude_N$ converges, as \mbox{$N\to\infty$}, to
the~uniform distribution on the~unit interval~$[0,1]$.

		In other words (with a small abuse of notation), 
		the random trajectory $\left(\mypoint_\tim(\tab_N) \right)_{\tim \in [0,1]}$ 
		with respect to the supremum norm converges in distribution to 
		a random meridian $\left( \point_{1-\tim, \psi} \right)_{\tim \in [0,1]}$ 
		when $N\to \infty$, 
		that is shortly, 
		\[
		\left( \mypoint_\tim(\tab_N) \right)_{\tim \in [0,1]} 
        \xrightarrow[N \to \infty]{\on{d}} 
        \left( \point_{1-\tim, \psi} \right)_{\tim \in [0,1]},
		\]
		where $\psi$ is a random variable with the uniform $U(0,1)$ distribution. 
\end{restatable}

The proof is postponed to \cref{sec:proof-of-evacuation-for-squares}
and the preparations to it will take the majority of time.

\subsection{The third main result. Typical \jdtp} 
\label{sec:jdt}

Let $\tab$ be a~standard Young tableau with $n$ boxes. Sometimes when
investigating the \jdtp, we are concerned not only about its shape, but also we
would like to be able to tell which entries were placed in the~rearranged
boxes. This motivates the~notion of \emph{the \jdtp in the lazy
    parametrization} (or shortly \emph{the lazy \jdtp}) which is a sequence of boxes
$\pat(\tab) :=(\pat_1,\dots,\pat_n) \subset \N^2$ where $\pat_i$ is the last
box along the \jdtp corresponding to $\tab$ (cf.~\cref{fig:jdtA}) which
contains a number $\leq i$, cf.~\cite[Section~3.3]{Romik2015a}.

Our third main result is stated in the~language of lazy \jdtp
and says that, asymptotically, 
\emph{the scaled \jdtp in~a~random square tableau 
is, just like in~\cref{thm:evacuation}, a~random meridian}, 
see~\cref{fig:simulation-jdtZ}.

\begin{theorem}
    \label{thm:jdt}
For each $N\in \N$ there exists a random variable~$\LongitudeTilde_N\colon
\tableaux_{\sq_N} \to [0,1]$ such that the supremum distance
\[\sup_{\tim \in [0,1]} 
\left| 
 \frac{1}{N} \pat_{\lceil \tim N^2\rceil}(\tab_N)  -
P_{\tim,\LongitudeTilde_N(\tab_N)} \right| 
\]
converges in probability to zero, as $N\to\infty$. 

The probability distribution
of the random variable $\LongitudeTilde_N$ converges, as \mbox{$N\to\infty$},
to the uniform distribution on the unit interval $[0,1]$.

		In other words (with a small abuse of notation), 
		the lazy \jdtp $\left( \pat_{\lceil \tim N^2\rceil}(\tab_N) \right)_{\tim \in [0,1]}$
		with respect to the supremum norm converges in distribution to 
		a random meridian $\left( \point_{\tim, \psi} \right)_{\tim \in [0,1]}$ 
		when $N\to \infty$, 
		that is shortly, 
		\[
		\frac{1}{N} \pat_{\lceil \tim N^2\rceil}(\tab_N) \xrightarrow[N \to \infty]{\on{d}} \point_{\tim, \psi}
		\]
		where $\psi$ is a random variable with the uniform $U(0,1)$ distribution. 
\end{theorem}

The proof is postponed to \cref{sec:equivalence-of-main-problems}
and is based on showing that the random lazy \jdtp and the (reversed) random evacuation path
have the same distribution, cf.~\cref{thm:equivalence}.
In other words, we show that the problem of finding typical \jdtps 
is equivalent to the problem of finding typical evacuation paths from \cref{thm:evacuation}.

\subsection{Conjecture on the independence of the iterated sliding paths.} 
\label{sec:conjecture}

Recall from \cref{sec:jdt-definition} that the \jdtt $\jdtcomplete$ is a
bijection defined by a two-step procedure in which first we apply the
\jdtincomplete and then we add the box which was removed from the initial tableau
$\tab$. We can iterate $\jdtcomplete$ on the random tableau $\tab_N$ and get
\emph{the sequence of iterated \jdtps}
\begin{equation}
\label{eq:iterated-jdtps}
\pat \! \left( \tab_N\right), 
\quad
\pat \! \left(\jdtcomplete(\tab_N)\right), 
\quad
\pat \! \left((\jdtcomplete)^2(\tab_N)\right), 
\quad 
\dots.
\end{equation}
Since $\jdtcomplete$ is a bijection, \cref{thm:jdt} provides asymptotics
of the distribution of each element of the sequence \eqref{eq:iterated-jdtps} separately.
The following conjecture aims to provide asymptotic information about their \emph{joint}
distribution. 

\begin{conjecture}
For each integer $k\geq 1$ 
the probability distribution
of the random vector 
\[ 
\left( \LongitudeTilde_N(\tab_N), 
\quad 
\LongitudeTilde_N\big(\jdtcomplete (\tab_N)\big), 
\quad \dots, \quad 
\LongitudeTilde_N\big((\jdtcomplete)^{k-1} (\tab_N)\big) \right)
\]
converges, as \mbox{$N\to\infty$},
to the uniform distribution on the unit cube $[0,1]^k$.
\end{conjecture}

Romik and the second named author proved an analogue of this conjecture 
in the case of the Plancherel-distributed 
random infinite tableaux, see  
\cite[the comment below Theorem~1.5]{Romik2015a}
and \cref{sec:old-new} for more details.

\subsection{Something old, something new, something borrowed: \jdtps for random infinite tableaux}
\label{sec:old-new}

The problem of the asymptotic shape of a \jdtp (analogous to the
one from \cref{thm:jdt}) was studied before by Romik and the~second named author
\cite[Section~1.3]{Romik2015a} for a certain \emph{infinite} random Young
tableau,  more specifically for the recording tableau $Q(x_1,x_2,\dots)$
obtained by applying the Robinson--Schensted correspondence to an infinite sequence
$(x_1,x_2,\dots)$ of independent, identically distributed random variables with
the uniform distribution on the unit interval~$[0,1]$. Such a choice
corresponds to sampling the random Young tableau according to, so called,
\emph{the Plancherel measure}.

In such a setting the \jdtp happens to converge almost surely to a~straight
line with a~random direction \cite[Theorem~1.1]{Romik2015a}, in other words
the~analogue of our meridians in the context of the Plancherel measure is
given by the~straight lines emanating from the origin of the coordinate
system.

\medskip

One of the main difficulties in the proof of \cref{thm:jdt} will be the
construction of the random variables $\LongitudeTilde_N$ which provide the
longitude of the meridian along which the sliding path travels.
This difficulty is absent in \cite{Romik2015a} because in that context 
the analogue of the longitude turns out to be simply equal to $x_1$, the first
entry of the random sequence to which the Robinson--Schensted correspondence is
applied.

It would be tempting to repeat the approach from \cite{Romik2015a} in our
context, for example one could proceed as follows. Let 
\[ \pi^{(N)}=\left(
\pi^{(N)}_1,\dots, \pi^{(N)}_{N^2} \right)\] 
be a uniformly random element of
the set of \emph{extremal Erdős--Szekeres permutations}, i.e., permutations with
the property that the corresponding tableaux associated via the Robinson--Schensted
correspondence have the square shape~$\Box_N$. Then the corresponding recording
tableau $\tab_N=Q(\pi^{(N)})$ is, as required, a uniformly random standard Young
tableau of square shape $\sq_N$. A naive guess would be that one possible choice for
the random variable $\LongitudeTilde_N$ is again (a rescaled version of) the
first entry of the permutation, i.e.,~$\pi^{(N)}_1$.

Uniformly random extremal Erdős--Szekeres permutations were investigated by
Romik \cite{Romik2006} who proved, among other results, that the probability distribution of
$\frac{1}{N^2} \pi^{(N)}_1$ converges to the point measure concentrated in~$\frac{1}{2}$. 
For this reason it seems that the random variables
$\pi^{(N)}_1$ do not carry any information which would be useful for our
purposes, hence the approach from \cite{Romik2015a} is not applicable here
directly, and the construction of the random variables $\LongitudeTilde_N$
must follow different ideas.

\medskip

Despite this fundamental difference, an astute reader may notice that our proof
of \cref{thm:evacuation} follows a path parallel to the one of
\cite[Theorem~5.1]{Romik2015a}. For example, the counterpart of our
\cref{prop:expvalue} (which can be viewed as a result about a certain random process
of \emph{removal} of boxes from a Young diagram) is
\cite[Theorem~4.4]{Romik2015a} (which concerns a certain random process
of \emph{addition} of boxes to a Young diagram).

\section{Preliminaries}
\label{sec:definitions}

\subsection{Permutations, Young diagrams and Young tableaux, continued}

We continue \cref{sec:tableaux-begin} where some basic definitions were
introduced.

By $\N$ we denote the~set of positive integers and we denote $\N_0 := \N \cup
\{0 \}$. For any natural number~$n$, we define the~set $[n]:=\{1, \ldots, n\}$.
From the following on  (unlike in \cref{sec:random-sorting-networks}) we will
view the symmetric group $\Sym_{n}$ as the group of permutations of 
the set $[n]$. 
We define \emph{the~length of~a permutation $\pi$} to be 
the~minimal number of factors necessary to write $\pi$
as a~product of (arbitrary, not necessarily adjacent!) transpositions, 
and denote it by~$|\pi|$. 

For any Young diagram $\lambda$, we denote \emph{the~number of standard Young tableaux} of shape
$\lambda$ by~$f^\lambda=\left| \tableaux_\lambda \right|$. 
If $\tab \in \tableaux_\lambda$ and $0 \leq p \leq |\lambda|$ we consider the
\emph{restriction of $\tab$ to its $p$ least boxes} by removing the entries
which are bigger than $p$, and denote the~obtained tableau by $\tab|_{\leq p}$
(clearly, it is also a~standard Young tableau).

For $C \geq 1$ we say that a Young diagram $\lambda$ is \emph{$C$-balanced} if $\lambda$ has
at most $C\sqrt{|\lambda|}$ rows and at most $C\sqrt{|\lambda|}$ columns.

For a tableau $T$ and an integer $p$ which appears exactly once in $\tab$
(which will always be the case in our considerations)
we define \emph{the $u$-coordinate of the box with the~entry~$p$} as
\[
\content^\tab_p := x - y \qquad \text{ for $(x,y) = \pos_p(\tab)$}.
\]
Note that in the literature, for instance
in~\cite{Ceccherini-Silberstein2010}, such a \mbox{$u$-coordinate} is called
the \emph{content}. 

We will also consider skew tableaux obtained by removing some boxes from Young
tableaux. Let $\tab \in \tableaux_\lambda$ be a standard Young tableau. If
$\tab$ with the~boxes with entries $a_1, \dots, a_i$ removed is a skew
tableau, we denote it by $\tab \setminus \{a_1, \dots, a_i\}$. Clearly, $\tab
\setminus\{p+1, \dots, |\lambda|\} = \tab|_{\leq p}$ is also a standard Young
tableau.

\subsection{Representation theory}

If $G$ is a finite group and $\rho\colon G \to \on{End} V$ is its
representation on a~finite-dimensional complex linear space~$V$, then by
\[\xi_V(g) := \on{Tr} \rho(g), \qquad g \in G,\]
we denote its \emph{character} (we write just $\xi$ if it is clear which representation we consider). 
We also consider the \emph{normalized character} $\chi_{_V}\colon G\to\C$ given by 
\[\chi_{_{V}} := \frac{1}{\on{dim} V} \cdot \xi_V = \frac{1}{\xi_V(\on{id})} \cdot \xi_V.\]
Additionally, for any element of the group algebra
\[ f= \sum_{g \in G} f_g g \in \C G \]
we denote by
$\chi_{_{V}}(f)$ the \emph{extension of the character by linearity} given by
\[\chi_{_{V}}(f) := \sum_{g \in G} f_g\, \chi_{_{V}}(g) \in\C.\]

On the vector space of functions $\mathcal{G}:=\{f \colon G \to \C \}$ 
we consider the \emph{standard scalar product} given by
\[\ls f, g \rs := \frac{1}{|G|} \sum_{h \in G} f(h)\ \overline{g(h)}, \qquad \text{for }f,g \in \mathcal{G},\]
where $\overline{c}$ denotes the complex conjugation. 

We will denote by $\hat{G}$ the family of irreducible representations
(\emph{irreps} for short) of $G$. It is known that the system of irreducible characters
$\{\xi_{V_x}\}_{x \in \hat{G}}$ is orthonormal. Therefore for the normalized
irreducible characters the following holds:
\[
\ls\chi_{_{V_1}}, \chi_{_{V_2}} \rs = 
\left\{ \begin{array}{ll} \frac{1}{\on{dim} V_1} & 
\text{if \ } \chi_{_{V_1}} = \chi_{_{V_2}}, \\[1.5ex] 
0 &\text{otherwise.} \end{array} \right.
\]

Let $V$ be a finite-dimensional representation and let
\[ 
V = \displaystyle \bigoplus_{x \in
    \hat{G}} m_x V_x\] 
be its decomposition into irreducible components, where $m_x\in\N$ denotes the
multiplicity of $x$. By \emph{a random irreducible component} of $V$ we will
understand a random element of $\hat{G}$ sampled according to the probability
measure $\Ps_V$ which is proportional to the total dimension of all copies of a
given irrep in $V$:
\begin{equation}
    \label{eq:probformula}
\Ps_V(x) := \frac{m_x \on{dim} V_x}{\on{dim} V} =
\left(\on{dim} V_x\right)^2 \cdot \ls \chi_{_{V}}, \chi_x \rs \qquad \text{for } x\in\hat{G}.
\end{equation} 

The \emph{trivial representation} will be denoted by $\on{triv}$. 
If $H$ is~a~subgroup of $G$ then 
we denote by $\rho \downarrow^G_H$ 
the~\emph{restriction of a representation $\rho$ to $H$} 
(if $G$ is fixed we just write $\rho \downarrow_H$).  

\medskip

For $\lambda\in\Y_n$ we denote by $\rho_\lambda \colon \Sym_n \to \on{End} V_\lambda$ 
\emph{the irreducible representation of the symmetric group $\Sym_n$ 
corresponding to the Young diagram $\lambda$}
and by $\chi_{_\lambda}$ its \emph{normalized character}.

\subsection{Asymptotics of characters and the approximate factorization property}
\label{sec:asymptotics-characters}

We will use two results concerning asymptotics of characters. 
The first one, due to Fer\'ay and the~second named author, 
gives an upper bound on the irreducible characters.

\begin{fact}[{\cite[Theorem~1]{Feray2011}}]
\label{thm:character-estimation}
There exists a~constant~$a>0$ 
such that for any Young diagram $\lambda$
and any permutation $\pi \in \Sym_{|\lambda|}$
\begin{equation}
\label{eq:FeraySniady}
\left|\chi_{_\lambda}(\pi) \right| 
\leq \left[ a \max\left(
\frac{r(\lambda)}{|\lambda|}, 
\frac{c(\lambda)}{|\lambda|}, 
\frac{|\pi|}{|\lambda|} 
\right) \right]^{|\pi|}
\end{equation}
where $r(\lambda)$ and $c(\lambda)$ stand, accordingly, 
for the~number of rows 
and the number of columns of $\lambda$.
\end{fact}

The second result, due to Biane (and its generalization due to the second-named author), 
shows that the~character calculated on the~product of two fixed permutations 
with disjoint supports approximately factorizes.

\begin{fact}[{\cite[Corollary~1.3]{Biane1998}, 
\cite[Section~0]{Biane2001}, \cite[Theorem~1]{Sniady2006a}}]
\label{thm:character-value}
Let $C \geq 1$ and $m \in \N$. 
There exists a~constant~$K>0$
such that for each $C$-balanced Young diagram $\lambda$ 
and all permutations $\sigma, \tau \in \Sym_{|\lambda|}$ 
with disjoint supports and satisfying $|\sigma|, |\tau| \leq m$
we have
\[
\left| \chi_{_\lambda}(\sigma \tau) 
- \chi_{_\lambda}(\sigma) \chi_{_\lambda}(\tau) \right|
\leq  \frac{K}{\left(\sqrt{|\lambda|}\right)^{|\sigma|+|\tau|+2}}.
\]
\end{fact}

\subsection{Jucys--Murphy elements}
\label{subsec:Jucys--Murphy-elements-introduction}

In the applications of \cref{thm:character-estimation} and \cref{thm:character-value}
we will encounter the expressions in the left-hand-side of \eqref{eq:jucys-estimation}.  
The lemma below gives an upper bound for their values. 

Its proof uses \emph{Jucys--Murphy elements} which often appear in 
the modern approach to the representation theory of the symmetric groups (see~\cite {Ceccherini-Silberstein2010})
and are defined as the following elements of the symmetric group algebra $\C \Sym_n$

\begin{equation}
    \label{eq:JM}
    J_k := \sum_{1\leq i<k} (i, k) = (1, k) + \dots + (k-1, k) \in \C \Sym_n \quad \text{for } 1\leq k\leq n.
\end{equation}
The elements $J_1,\dots,J_n\in \C \Sym_n$ form a commuting family in the symmetric group algebra. 
We will use them intensively in \cref{sec:udist}.

		\begin{lemma}
    \label{lem:sumoversymmetric}
		For any $c>0$ and $n\in \N$
		\begin{equation}
		\label{eq:jucys-estimation}
    \sum_{\pi \in \Sym_n} c^{|\pi|} <
    \exp \left( \frac{n^2 c}{2} \right). 
		\end{equation}
		\end{lemma}

		\begin{proof}
				For any $c \in \C$ the following simple identity in the symmetric group algebra holds true 
				\[ \sum_{\pi \in \Sym_n} c^{|\pi|} \pi = (1+c J_1) (1+c J_2) \cdots (1+c J_n). \]
				By applying the~trivial representation to both sides of the above equality we get that for $c>0$
				 \begin{multline*}     \sum_{\pi \in \Sym_n} c^{|\pi|}= 
				(1+c)(1+2c)\cdots\big(1+(n-1)c\big) < \\ 
				e^c e^{2c} \cdots e^{(n-1)c}< \exp \left( \frac{n^2 c}{2} \right).
                \qedhere
		\end{multline*}
		\end{proof}

\section{The longitude and surfing}
\label{sec:longitude-begins}

Our strategy towards the proof of \cref{thm:evacuation} is to pass to the
geographic coordinates of the point $X_t=X_t(T_N)$ from the scaled evacuation
trajectory~\eqref{eq:evacuation}. Having the choice between the latitude and
the longitude, we start with the more challenging problem of understanding how
the longitude $\longitude(X_t)$ changes over time $t$. 

Instead of considering the longitude $\longitude(X_t)$ directly, 
it will be more convenient to study the following random variable
which we call \emph{the theoretical longitude}:
\begin{equation}
    \label{eq:estimate-longitude}
     \TheoreticalLongitude(t):= F_{\cotmeas_{_{1-t}}}\big( u(X_t)\big),  
\end{equation}
where $ u(X_t)$ denotes the $u$-coordinate of $X_t=X_t(T_N)$, and $ F_{\cotmeas_{_{1-t}}}$ is the
cumulative distribution function of the limit measure $\cotmeas_{1-t}$ which was defined in \cref{sec:random-position}.
Notice that if in time $t$ the box with the biggest number is positioned 
\emph{exactly} on the circle of latitude $\alpha = 1-t$,
that is, $\alpha(X_t) = 1-t$, 
then the theoretical longitude coincides with the longitude, 
i.e., \mbox{$\TheoreticalLongitude(t) = \longitude(X_t)$}.
Heuristically one would expect that $\TheoreticalLongitude(t)\approx \longitude(X_t)$ for $N\to\infty$.

Roughly speaking, the following result states that (away from the polar regions
which correspond to $t=0$ and $t=1$) the theoretical longitude of~$X_t$
does not change too much over time.

\begin{restatable}{theorem}{thm-all-the-same}
    \label{thm:all-the-same}

    Assume that $0<\tim_1 < \tim_2 <1$. Then for each $\varepsilon>0$
    \[ 
		\lim_{N\to\infty} 
    \Pp \bigg\{ \tab_N \in\tableaux_{\sq_N} : 
    \big| \TheoreticalLongitude(\tim_2) - \TheoreticalLongitude(\tim_1) \big| > \varepsilon \bigg\} = 0.
		\]
		In other words,
		the difference $\TheoreticalLongitude(\tim_2) - \TheoreticalLongitude(\tim_1)$ 
		converges in probability to~$0$ when $N\to \infty$, 
		that is shortly, 
		\[
		\TheoreticalLongitude(\tim_2) - \TheoreticalLongitude(\tim_1) \xrightarrow[N \to \infty]{P} 0. 
		\]
\end{restatable}	

The proof is quite involved;
\crefrange{sec:longitude-begins}{sec:propproof} 
are a preparation while
the proof of \cref{thm:all-the-same} itself will be given
in~\cref{subsec:proof-theoretical-longitude-b,subsec:proof-theoretical-longitude-a}. 
The remaining part of the current section is devoted
to a rough sketch of the proof.

\medskip

Clearly, \cref{thm:all-the-same} is equivalent to the conjunction 
of the following two statements for $\varepsilon>0$:
\begin{align}
    \lim_{N\to\infty} \Pp \bigg\{ \tab_N \in\tableaux_{\sq_N} :  
    \TheoreticalLongitude(\tim_2) - \TheoreticalLongitude(\tim_1)  > \varepsilon \bigg\} &=0,
    \label{eq:sink-B}
		\\
    \lim_{N\to\infty} \Pp \bigg\{ \tab_N \in\tableaux_{\sq_N} : 
    \TheoreticalLongitude(\tim_2) - \TheoreticalLongitude(\tim_1)  < -\varepsilon \bigg\} &= 0.
    \label{eq:sink-A} 
    \end{align}
We start with the proof of the upper bound \eqref{eq:sink-B}. Then we will use
the symmetry of the problem in order to prove the lower bound
\eqref{eq:sink-A}.

\subsection{The single surfer scenario}
\label{sec:single}

We would like to present the problem of 
the~evacuation path in a different, more vivid light. 
We will speak about a~\emph{square pool of side~$N$} 
(=the square Young diagram $\sq_N$)
filled with \emph{$N^2-1$ particles of~water} 
(=the Young tableau $\tab_N$ with
the largest entry removed), 
a~passive \emph{surfer} (=the box with the biggest
entry~$N^2$) and its trajectory (or behavior) when 
the~pool is being drained (=iteratively applying \jdtincomplete). 
Our~goal~in
\cref{thm:evacuation} is to show that, when the pool is big enough, the surfer has some
typical paths along which he/she moves as the pool is being drained.

In our proof of \cref{thm:all-the-same} we start our analysis at time $\tim_1$
when \jdtincomplete was already applied 
\begin{equation} 
    \label{eq:m1}
    m_1:=\lfloor \tim_1 N^2 \rfloor
\end{equation} 
times.
Our starting point is therefore the standard tableau 
\begin{equation}
\label{eq:initial-surfer}
\tab'_N:={\jdtcomplete}^{m_1}(\tab_N)\big|_{\leq N^2-m_1}
\end{equation}
with 
$N^2-m_1$ boxes (compare with \eqref{eq:jdt-tableaux-sequence}). 
We denote 
\begin{equation}
    \label{eq:w1}
    \water_1=N^2-m_1-1. 
\end{equation}
In this way the boxes with numbers $1,\dots,\water_1$ correspond
to~the~\emph{water} and the box with the maximal number $\water_1+1$ to the
\emph{surfer}. The position of the latter box will be called \emph{the initial
    position of the surfer} 
		and we will refer to the tableau $\tab'_N$ 
		as \emph{the initial surfer configuration}. 
		By removing the box with the surfer
\begin{equation}
    \label{eq:singlewater} 
W_N':= \tab_N' \setminus \{\water_1+1\}
\end{equation}
we get a standard Young tableau which encodes \emph{the initial configuration of the~water}.

As time goes by, at the time $\tim_2$ the \jdtincomplete was already applied 
\[ m_2:=\lfloor \tim_2 N^2 \rfloor\] 
times and we investigate the tableau
\[ \tab''_N={\jdtcomplete}^{m_2}(\tab_N)\big|_{\leq N^2-m_2} = {\jdtcomplete}^{m_2-m_1}(\tab'_N)\big|_{\leq N^2-m_2}\]
with $w_2+1$ boxes, where
\begin{equation}
    \label{eq:w2}
    \water_2=N^2-m_2-1. 
\end{equation}
The boxes with the numbers $1,\dots,w_2$ correspond to the remaining particles
of water and the box with the maximal number $w_2+1$ corresponds to the surfer; the
position of the latter box will be called \emph{the final position of the
    surfer}.

\medskip

Our aim is to relate the final position of the surfer at the time $t_2$ to its
initial position at the time $t_1$, preferably in the language of the theoretical longitude.
As a point of reference we will introduce an additional \emph{multisurfer story} 
which happens in a parallel universe in which we pay attention to $k$ surfers.

\subsection{Pieri tableaux}
\label{sec:Pieri}

We consider the~partial order on the plane $\R^2$ defined by:
\[(x_1, y_1) \preceq (x_2, y_2) \iff \left(x_1 \leq x_2 \ \wedge \ y_1 \geq y_2\right).\]

Let $k$ be a fixed natural number and $\multi$ be a tableau in which the
$k$~largest entries are numbered by~consecutive integers $l+1,\dots,l+k$. We
say that the tableau $\multi$ is a \emph{$k$-Pieri tableau} if these $k$~largest
boxes are placed in \emph{the~increasing order} with respect to $\preceq$
(i.e., they are placed from north-west to south-east) or, equivalently, their
$u$-coordinates are ordered increasingly, that is:
\[ \content^{\multi}_{l+1} < \cdots < \content^{\multi}_{l+k}.\]
If the value of the number $k$ is clear from the context,
we will shortly say that $\multi$ is \emph{Pieri}.

It is easy to check that if $\multi$ has at least $k+1$ boxes then $\multi$
is~a~$k$-Pieri tableau if and only if $j(\multi)$ is a $k$-Pieri tableau.

\medskip

For standard Young tableaux we will consider the following more general
notion. For a (skew) standard Young tableau $\tab$ with $n$ boxes and positive
integers $w$ and $k$ such that $w+k\leq n$ we say that $\tab$ is a
\emph{$(w+1,w+k)$-Pieri tableau} if
\begin{equation}
    \label{eq:Pieri}
    \content_{w+1}^\tab < \dots < \content_{w+k}^\tab.
\end{equation}
The set of $(w+1,w+k)$-Pieri standard tableaux of (skew) shape $\lambda$ will be
denoted by~$\widetildetableaux_{\lambda}^{(w+1,w+k)}$.

\medskip

\subsection{The multisurfer scenario}
\label{sec:multi}

For \emph{the multisurfer scenario} let $k=k(N)$ be a sequence of positive
integers such that
\begin{equation}
\label{eq:minimal-assumption-k}
     \lim_{N \to \infty} k=\infty\quad \text{ and } \quad \lim_{N\to\infty} \frac{k^2}{N} = 0. 
\end{equation}
For a fixed value of $N$ we consider the $N\times N$ square pool filled with $N^2-k$
\emph{particles of water} on which $k$~surfers (=$k$~boxes with the biggest
entries) are positioned in the~increasing order.  Formally speaking, by
\[
\widetildetableaux_{\sq_N}=\widetildetableaux_{\sq_N}^{(N^2-k+1,N^2)}
\] 
we denote the set of standard Young tableaux of the square shape~$\sq_N$ which are
$k$-Pieri, and~by~$\Pt$~the~uniform distribution
on~the~set~$\widetildetableaux_{\sq_N}$. 
We assume that $\multi_N$ is a~random
tableau sampled with the~uniform probability distribution~$\Pt$ on the set~$\widetildetableaux_{\sq_N}$
In this scenario, in order to refer to the~$k$~surfers, we will use the name
\emph{multisurfers}.

We start our analysis when \jdtincomplete was already applied $m_1+1-k$ times
with $m_1$ given again by \eqref{eq:m1} (notice that $m_1+1-k \geq 0$ for $N$
big enough). Our~starting point is therefore the tableau
\begin{equation}
    \label{eq:initial-multisurfers}
 \multi'_N:={\jdtcomplete}^{m_1+1-k}(\multi_N)\big|_{\leq N^2-m_1-1+k}
\end{equation}
with $N^2-m_1-1+k=\water_1+k$ boxes. We will refer to this tableau as 
\emph{the initial multisurfer configuration}. 
In this way, just as in the~single surfer scenario,
the~boxes with the numbers $1,\dots,\water_1$ correspond to the \emph{water}
(in particular, there is the same number of water particles as in the single
surfer scenario). On the other hand, the boxes with the~numbers
$\water_1+1,\dots,\water_1+k$ correspond to the \emph{multisurfers}. By
removing the~multisurfers
\begin{equation}
    \label{eq:multiwater} 
\widetilde{W}_N':= \multi_N' \setminus
\{\water_1+1,\dots,\water_1+k\}
\end{equation}
we get a standard Young tableau which encodes the initial configuration of the~water.

As time goes by, at the time $\tim_2$ the \jdtincomplete was already applied $m_2+1-k$ 
times and we investigate the tableau
\begin{equation}
    \label{eq:multi-final} 
    \multi''_N={\jdtcomplete}^{m_2+1-k}(\multi_N)\big|_{\leq w_2+k} = {\jdtcomplete}^{m_2-m_1}(\multi'_N)\big|_{\leq w_2+k} 
\end{equation}
which consists of $w_2 + k$ boxes which correspond to $w_2$ particles of water and $k$ multisurfers; we will refer
to this tableau as \emph{the final multisurfer configuration}.

\subsection{Sketch of the proof of the upper bound \eqref{eq:sink-B}}
\label{sec:plan}

\subsubsection{The collective behavior of the multisurfers is not very random}
\label{subsec:collective-behaviour}
    
Since the number of the multisurfers is small in comparison to the number of the rows/columns,
as a first-order approximation we may treat the set of positions of
$k$~multisurfers at any fixed time as a collection of $k$ independent copies
of the position of a single surfer. We can expect therefore that the law of
large numbers is applicable and, as $k\to\infty$, \emph{the multisurfer empirical
measure} (which is a random measure which encodes the scaled down 
$u$-coordinates of the multisurfers) 
converges in probability to the probability distribution of
the position of the single surfer. 
In other words: \emph{the collective behavior of the multisurfers is much less random 
than the behaviour of the single surfer}.
This phenomenon is beneficial and will allow us to use the multisurfers as a
moving frame of reference for tracing the position of the single surfer over
time.

The above naive first-order approximation clearly cannot be true if \mbox{$k=k(N)$}
grows too fast with the size $N$ of the square. Nevertheless, in
\cref{prop:longitude-experimental} we will show that if $k=k(N)$ grows at the
right speed then a version of the law of large numbers indeed holds true. The
proof of \cref{prop:longitude-experimental} is quite technically involved and
the whole \cref{sec:udist} is devoted to its proof.

\subsubsection{The multisurfers provide information about the surfer}
\label{subsec:multisurfers-with-info}

Let us fix some common initial configuration of the water for both the surfer
and the multisurfers. In another first-order approximation let us assume for a
moment that the density of the multisurfers is small enough that during the
time interval $[t_1,t_2]$ all neighboring pairs of multisurfers are separated
so that the multisurfers do not touch each other. If this is indeed the case
and there are no interferences between the multisurfers then the time
evolution of each multisurfer clearly coincides with the time evolution of the
single surfer who would have the same initial position; by reversing the
optics this means that we have a very direct information about some specific
single surfer trajectories (namely, the ones which start from the positions of
the multisurfers) in terms of the dynamics of the multisurfers, which we
understand pretty well thanks to the aforementioned 
\cref{prop:longitude-experimental}.

It is very convenient that \emph{for a fixed initial configuration of the
    water, the trajectory of the single surfer depends in a monotonic way on the
    initial position of the surfer}. For this reason it is possible to get some
partial information also about the single surfer trajectories
starting from the points \emph{between} the initial positions of two
multisurfers. If the number $k=k(N)$ of the multisurfers tends to infinity as
$N\to\infty$ such neighboring multisurfers should not be too far (in
comparison to the size $N$ of the square) which is enough to prove
\cref{thm:all-the-same}.

\subsubsection{The single surfer and the multisurfer scenario on the same water configuration}
\label{subsec:sketch-single-multi}

Above we used the idea of considering the single surfer scenario \emph{and}
the multisurfer scenario on the same configuration of water. This idea sounds
self-contradictory because each of these two scenarios gives rise to a
\emph{different} probability distribution on the set of configurations of the
water. In~\cref{sec:single-multi} we shall explain how to overcome this
difficulty and to (asymptotically) couple the surfer and the multisurfers on a single
probability space. 
The resulting object can be visualized as water on which in
two parallel universes there is (i) a single surfer, and (ii) $k$
multisurfers. The single surfer and the multisurfers are like ghosts 
to one another
and do not interact.
Furthermore, \emph{as long as the multisurfers do not touch each other, 
the relative position (with respect to the partial order $\prec$ on the plane) of the surfer and the ghosts 
of the multisurfers 
does not change over time: 
overtaking of the surfer by the multisurfers is not allowed}.

\subsubsection{Overtaking is allowed in one direction only}
\label{subsec:overtaking}

The above discussion was based on a simplistic assumption that the
multisurfers do not touch each other. Regretfully, in the real world this is not
the case; multisurfers might influence each other and hence the multisurfer
trajectories might differ from the single surfer trajectories on the same
configuration of the water. 

On the bright side, the assumption that the
multisurfers are ordered as in the definition of the Pieri tableaux implies that
the movement of each multisurfer depends only on (1) the configuration of the water,
and (2) on these multisurfers which are to the north-west;
the other multisurfers which are south-east
have no influence on its dynamics. 
Furthermore, the impact of the multisurfers is unidirectional: 
the presence of the north-west multisurfer-neighbor 
can only push the multisurfer in the south-east direction. 
For the `coupling' of the stories of the single surfer and the multisurfers
on the common water this means that 
the ghosts of the multisurfers \emph{are} allowed to overtake the surfer, but only in one direction. 
More precisely, the number of the multisurfers which are north-west to the single surfer
can only decrease over time 
(see \cref{lem:ghosts} for the formal description of the above heuristics).

\subsubsection{The order in which the proof of the upper bound \eqref{eq:sink-B} will be conducted}

The rigorous proof of \eqref{eq:sink-B} will be conducted in the following order.
First in \cref{sec:single-multi} we will develop the content of \cref{subsec:sketch-single-multi}.
Then in \cref{sec:udist} we will deal with 
the collective behavior of the multisurfers described in \cref{subsec:collective-behaviour}.
In \cref{subsec:overtaking-one-direction} 
we will work out the dynamics contained in \cref{subsec:overtaking}.
Finally, in the rest of \cref{sec:connection} 
we will formally justify the heuristics from \cref{subsec:multisurfers-with-info}
and complete the proof.

\section{Single surfer versus the multisurfer scenario}
\label{sec:single-multi}

In \cref{sec:single,sec:multi} we considered two random tableaux: $\tab_N'$
and $\multi_N'$ defined on~two~different probability spaces (which correspond
to the single surfer and the multisurfer scenario respectively). In order to
proceed with the ideas sketched in \cref{subsec:sketch-single-multi} we need to
define some random tableaux $\ttt$ and $\ttilde$ on a \emph{common} probability space with
(almost) the same distributions as $\tab_N'$ and $\multi_N'$ and such that the
corresponding configurations of water coincide. 
A solution to this problem is provided by \cref{prop:single-and-multi} below which is also the main result of the current section.
We will compare the distributions with respect to \emph{the total variation distance}.

\medskip

Suppose that $X$ and $Y$ are random variables (possibly defined on different probability spaces)
taking values in some finite set $S$
with probability distributions $\mathbb{P}_X$ and $\mathbb{P}_Y$ respectively. 
We define the~\emph{total variation distance} 
between the~random variables $X$ and $Y$ \cite[Section~3.6.1]{Durrett2010}
(or, alternatively, between the probability distributions
$\mathbb{P}_X$ and $\mathbb{P}_Y$) as 
\begin{multline} 
\label{eq:TVD1}
\delta( X, Y) = \delta( \mathbb{P}_X, \mathbb{P}_Y):= \\ \frac{1}{2} \sum_{s \in S} \left| 
\mathbb{P}_X (s) - \mathbb{P}_Y (s) 
\right| 
= \max_{Z \subset S} \left| \mathbb{P}_X (Z)  
- \mathbb{P}_Y (Z)  \right|.
\end{multline}
Sometimes we will also denote this quantity by $\delta(X,\mathbb{P}_Y)$, etc.

\newcommand{\mydelta}{\Delta}

\begin{proposition}
    \label{prop:single-and-multi} 
For each $C \geq 1$ and $\mydelta \in (0,1)$ 
there exists a~constant~\mbox{$d>0$} with the following property.

Let $\lambda$ be a $C$-balanced Young diagram
and let $k, w,a$ be positive integers such that 
$w < (1-\Delta)\ |\lambda|$ and
$w \leq a \leq |\lambda| - k$.

Then there exists a pair of random tableaux $\mathbf{T}$ and $\mathbf{M}$ which
are defined on the same probability space with the following properties:
\begin{enumerate}[label=(\alph*)]
    \item $\mathbf{T}$ is a uniformly random element of $\tableaux_{\lambda}$;
    \item \label{item:uretilde}
    $\mathbf{M}$ is a random element of $\widetildetableaux_\lambda^{(a+1,a+k)}$;
    \item
    \label{item:TVD} the total variation distance between the distribution of $\mathbf{M}$ 
		and the~uniform distribution on $\widetildetableaux_\lambda^{(a+1,a+k)}$ fulfills the bound
\begin{equation}\label{eq:TVD}
         \delta\left( \mathbf{M}, \mathbb{P}_{\widetildetableaux_\lambda^{(a+1,a+k)}} \right) < d \frac{k^2}{\sqrt{|\lambda|-w}}; 
\end{equation}
    
    \item \label{item:M=T} $\mathbf{T}\big|_{\leq w}=\mathbf{M}\big|_{\leq w}$ holds true almost surely.
\end{enumerate}

\end{proposition}

The proof is postponed to \cref{sec:proof-single-and-multi}.
In the next two subsections we prepare to it by proving general lemmas 
on the probabilities of some particular events 
in the multisurfers scenario (\cref{subsec:probability-Pieri})
and comparing the distributions of water beneath surfer(s) in both scenarios (\cref{sec:proof:lem:probdiv}).
It turns out that these distributions are similar in terms of the total variation distance 
which asymptotically converges to 0,
see \cref{lem:probdivB}.

\subsection{Probability that a random tableau is Pieri}
\label{subsec:probability-Pieri}

\begin{lemma}
    \label{prop:jdtinv}
    
    Let $\lambda/\mu$ be a skew Young diagram with $n$ boxes and let \mbox{$1\leq k \leq n$}.
    Then the cardinality of the set
\begin{equation}
    \label{eq:mysets}
     \widetildetableaux_{\lambda/\mu}^{(w+1,w+k)} \qquad \text{for $w\in\{0,\dots,n-k\}$} 
\end{equation}
    does not depend on the choice of $w$.
\end{lemma}

\begin{proof}
    There is a simple bijection between the set 
$\widetildetableaux_{\lambda/\mu}^{(w+1,w+k)}$ and the set of semistandard
skew tableaux of shape $\lambda/\mu$ and of weight $(1^w,k,1^{n-w-k})$ which is
defined as follows. For a tableau
$T\in\widetildetableaux_{\lambda/\mu}^{(w+1,w+k)}$ we replace each entry
from the set $\{w+1,\dots,w+k\}$ by the same number $w+1$, and we replace each
entry $i\in\{w+k+1,\dots,n\}$ by $i+1-k$. It follows therefore that
the cardinality of \eqref{eq:mysets} is equal to the coefficient
\[ \left[ x_1 \cdots x_w x_{w+1}^k x_{w+2} \cdots x_{n+1-k} \right] s_{\lambda/\mu}\]
in the expansion of the skew Schur function in the basis of monomials. Since
the skew Schur function is a symmetric polynomial, the proof is completed.
\end{proof}

Some notions defined for ordinary Young diagrams have their
natural counterparts in the skew setup, given as follows. 
For $C\geq 1$ we say that a skew Young diagram $\lambda/\mu$ 
is \emph{$C$-balanced} if $\lambda/\mu$ has
at most $C\sqrt{|\lambda/\mu|}$ rows and at most $C\sqrt{|\lambda/\mu|}$
columns. Moreover we denote by
$\mathbb{P}_{\lambda/\mu}$ the uniform measure on the set of standard Young
tableaux of skew shape $\lambda/\mu$. 

\medskip

We calculate now the~probability of choosing a Pieri tableau 
from the set of standard Young tableaux of a specified skew shape.

\begin{lemma}
    \label{lem:Pieriprob}
    For each $C\geq 1$ there exists a constant $c>0$ with the following property.
    Let $\lambda/\mu$ be a $C$-balanced skew Young diagram with $n$ boxes.
    Let $k$ and $w$ be integers such that $1\leq k<\sqrt[4]{n}$ and $0\leq w\leq n-k$. Then
     \begin{equation}
         \label{eq:pieri-proba}
         \left| k!\; \mathbb{P}_{\lambda/\mu} \! \left(\widetildetableaux_{\lambda/\mu}^{(w+1,w+k)}\right) - 1 \right| <  c \frac{k^2}{\sqrt{n}}. 
     \end{equation}
\end{lemma}
\begin{proof}
        Let $T$ be a uniformly random standard tableau of skew shape $\lambda/\mu$. It
is easy to check that $T$ is $(w+1,w+k)$-Pieri if and only if the rectified
tableau $\on{rect} T$ is $(w+1,w+k)$-Pieri 
(see \cref{sec:Pieri} for the definition of Pieri tableaux). 
In the following we will describe
the probability distribution of $\on{rect} T$.

    Recall that the plactic skew Schur polynomial $S_{\lambda/\mu}$ is the
formal sum of the elements in the plactic monoid which correspond to all
semistandard tableaux of shape $\lambda/\mu$. 
The relations in the plactic monoid \cite[Section~2, Corollary~1]{Fulton1997} 
		allow us to identify a
skew tableau with its rectification and to express the skew plactic Schur
polynomial as a linear combination of (non-skew) plactic Schur
polynomials \cite[Section~5.1, Corollary~4]{Fulton1997}:
\begin{equation}
    \label{eq:LR}
     \frac{S_{\lambda/\mu}}{f^{\lambda/\mu}} = \sum_{\nu} \frac{c^{\lambda}_{\mu \nu} f^\nu}{f^{\lambda/\mu}} \cdot \frac{S_{\nu}}{f^\nu},
\end{equation}
where $c^{\lambda}_{\mu \nu}$ is the Littlewood--Richardson coefficient. If we
restrict our attention only to the summands which correspond to (skew) \emph{standard}
tableaux, the~left-hand side can be identified with the
probability distribution of $\on{rect} T$. The right-hand can be interpreted
as a linear combination (over some Young diagrams $\nu$) of the uniform
measure on the set $\tableaux_{\nu}$.

In this way we proved that a random tableau with the same distribution as
$\on{rect} T$ can be generated by the following two step procedure. Firstly,
we select a random Young diagram $\nu$ with the probability distribution
\[ \mathbb{P}(\nu)=\frac{c^{\lambda}_{\mu\nu} f^\nu}{f^{\lambda/\mu}}. \]
Secondly, we select a uniformly random standard tableau with the shape $\nu$.

In particular, the probability that $\on{rect} T$ is a $(w+1,w+k)$-Pieri tableau
is a weighted arithmetic mean (over certain diagrams $\nu$) of the probability
that a uniformly random element of $\tableaux_\nu$ is a $(w+1,w+k)$-Pieri
tableau. It follows that it is enough to prove a version of the inequality
\eqref{eq:pieri-proba} in which the skew diagram $\lambda/\mu$ is replaced by
any diagram $\nu$ which contributes to~\eqref{eq:LR}; we will do it in the
following.

\medskip

We start with an observation that in the process of rectification the number of
rows and the number of columns of a tableau cannot increase. Since
$\lambda/\mu$ is $C$-balanced, 
it follows that any Young diagram $\nu$ which contributes to~\eqref{eq:LR} is also $C$-balanced.

By~\cref{prop:jdtinv} it is enough to consider the case $w=0$. Let $T$ be a
uniformly random standard tableau of shape $\nu$ and let $\xi=\on{sh} T|_{\leq
    k}$ be the positions of the first $k$ boxes of $T$. The tableau $T$ is
$(1,k)$-Pieri if and only if $\xi=(k)$ is the one-row diagram. The remaining
difficulty is therefore to identify the probability distribution of $\xi$.

The link between the combinatorics of standard Young tableaux and the
irreducible representations of the symmetric groups (in particular, the
branching rule) implies that the probability distribution of $\xi$ coincides with
the measure $ \Ps_{\left. V_{\nu} \right\downarrow_{\Sym_k}}$ 
on the irreducible components of $\left. V_{\nu} \right\downarrow_{\Sym_k}$ 
which was defined in \eqref{eq:probformula}. 
\cref{eq:probformula} gives therefore an exact formula
    \begin{multline}
        \label{eq:pieripieri}
        k! \, \mathbb{P} \! \left( T\in\widetildetableaux_{\nu}^{(1,k)} \right) 
        = k! \, \Ps_{\left. V_{\nu} \right\downarrow_{\Sym_k}}(\on{triv}) 
        = k! \ls \left. \chi_{{\nu}} \right\downarrow_{\Sym_k}, \chi_{\on{triv}} \rs = \\
         \sum_{\pi \in \Sym_k} \chi_{{\nu}}(\pi)\  \overline{\chi_{\on{triv}}(\pi)} 
        = 1 + \sum_{\substack{\pi \in \Sym_k \\ \pi \neq \on{id}}} 
        \chi_{\nu}(\pi).
    \end{multline}
    In the following we will find an asymptotic bound for the second summand on
    the right-hand side. 
    
		\medskip
		
    By \cref{thm:character-estimation} there exists a~universal constant~$a>0$
    (which depends only on $C$) such that for any $\pi \in \Sym_{k}$
    \begin{equation}
        \label{eq:FeraySniady-kwadrat}
        \left|\chi_{{\nu}}(\pi) \right|  \leq 
        \left( \frac{a}{\sqrt{n}} \right)^{|\pi|}.
    \end{equation}
    It follows that the second summand on the right-hand side of \eqref{eq:pieripieri} is bounded by
    \[
    \Bigg|\sum_{\substack{\pi \in \Sym_k \\ \pi \neq \on{id}}} \chi_{{\nu}}(\pi) \Bigg|  
    \leq 
    \sum_{\pi \in \Sym_k} \left(\frac{a}{\sqrt{n}}\right)^{|\pi|}
    - 1
    \leq e^{\frac{a k^2}{\sqrt{n}}}  - 1 = O\!\left(\frac{k^2}{\sqrt{n}}\right), 
    \]
    where we used \cref{lem:sumoversymmetric} and the assumption that $\frac{k^2}{\sqrt{n}}=O(1)$.
\end{proof}

\subsection{Comparison of distributions of water beneath surfer(s) in both scenarios}
\label{sec:proof:lem:probdiv}

The following lemma shows that in the asymptotic setting 
the probability distributions of water beneath surfer(s) in
the single surfer and the multisurfer scenarios are nearly equal.

\begin{lemma}
    \label{lem:probdivB} 
	For each $C\geq 1$ and $\mydelta\in(0,1)$ there exists a constant $d>0$ with the following property.
    
    Let $\lambda$ be a $C$-balanced Young diagram
		and $k, w$ and $a$ be positive integers such that 
        $w<(1-\mydelta) \ |\lambda|$ and $w \leq a \leq |\lambda|-k$.
		Denote by $T$ a~uniformly random element of $\tableaux_{\lambda}$ and
		by $M$ a~uniformly random element of
		$\widetildetableaux_{\lambda}^{(a+1,a+k)}$. 
		The total variation distance between the distributions of the restricted tableaux 
		$T\big|_{\leq w}$ and $M\big|_{\leq w}$ 
		fulfills the bound		
	\begin{equation}
			\label{eq:TVD-probdivB}
				 \delta\left( T\big|_{\leq w}, M\big|_{\leq w} \right) <  d \frac{k^2}{\sqrt{|\lambda|-w}}. 
	\end{equation}
 \end{lemma}
\begin{proof}
We start with an observation that the total variation distance is trivially
bounded from above by $1$. It follows that (provided $d\geq 1$) it is enough to
consider the case when $k^4\leq |\lambda|-w$.

\medskip

Notice that the~probability distribution
$\mathbb{P}_{\widetildetableaux_{\lambda}^{(a+1,a+k)}} (\cdot)$
coincides with the conditional probability 
$\mathbb{P}_{\tableaux_{\lambda}} \! 
\left( \cdot \: \big| \: \widetildetableaux_{\lambda}^{(a+1,a+k)} \right)$.
Therefore, for any $\mu \in \Y_w$ such that $\mu\subseteq \lambda$ and any $S
\in \tableaux_\mu$ we have, by the Bayes rule, that 
\begin{multline*}
\mathbb{P}_{\widetildetableaux_{\lambda}^{(a+1,a+k)}} \!  \left( M\big|_{\leq w} = S
\right) = \\ 
\frac{
    \mathbb{P}_{\tableaux_{\lambda}} \!  \left( M \in
    \widetildetableaux_{\lambda}^{(a+1,a+k)} \;\, \Big| \;\,  M\big|_{\leq w} = S
     \right) 
 }%
{\mathbb{P}_{\tableaux_{\lambda}} \!  \left(
    \widetildetableaux_{\lambda}^{(a+1,a+k)} \right)} 
\cdot
\mathbb{P}_{\tableaux_{\lambda}} \! \left( M\big|_{\leq w} =
S \right).
\end{multline*}
By elementary algebra this equality can be rewritten as
\begin{multline}
    \label{eq:lolcats2}
\mathbb{P}_{\widetildetableaux_{\lambda}^{(a+1,a+k)}} \! \left( M\big|_{\leq w} = S
\right)  
- \mathbb{P}_{\tableaux_{\lambda}} \! \left( M\big|_{\leq w} =
S \right)
 = \\
\shoveleft{ =\mathbb{P}_{\widetildetableaux_{\lambda}^{(a+1,a+k)}} \! \left( M\big|_{\leq w} = S
\right)  \!  \left[ 1- k! \ \mathbb{P}_{\tableaux_{\lambda}} \! \left(
\widetildetableaux_{\lambda}^{(a+1,a+k)} \right) \right] 
+   }  \\
+ \mathbb{P}_{\tableaux_{\lambda}} \! \left( M\big|_{\leq w} =
S \right) \!
 \left[k! \ \mathbb{P}_{\tableaux_{\lambda}} \! \left( M \in
 \widetildetableaux_{\lambda}^{(a+1,a+k)} \; \Big| \;  M\big|_{\leq w} = S
 \right)  -1 \right].
\end{multline}
Our strategy is to find an upper bound for the absolute value of the right-hand side.

\medskip

The conditional probability in the second summand on the right-hand side, i.e.,
\begin{equation}
    \label{eq:first-factor}
     \mathbb{P}_{\tableaux_{\lambda}} \! \left( M \in
    \widetildetableaux_{\lambda}^{(a+1,a+k)} \;\, \Big| \;\,  M\big|_{\leq w} = S
    \right), 
\end{equation} is equal
to the conditional probability that the restricted tableau $M |_{> w}$ is an
$(a+1,a+k)$-Pieri tableau. In order to calculate this conditional probability
we notice that the conditional probability distribution of the restricted
tableau $M |_{> w}$ (under the condition $M\big|_{\leq w} = S$) is the uniform
measure on the set of tableaux of shape $\lambda/\mu$ such that their
entries form the multiset \mbox{$(w+1,\dots,n)$}.
In other words, the probability
distribution of the random tableau \mbox{$\left( M |_{> w}\right)-w$} 
(which is obtained by decreasing each entry of $M |_{> w}$ by $w$) is given by
$\mathbb{P}_{\tableaux_{\lambda/\mu}}$.
In this way we proved that \eqref{eq:first-factor} is equal to
\begin{equation}
    \label{eq:trololo}
      \mathbb{P}_{\tableaux_{\lambda/\mu}}\! \left(  \widetildetableaux_{\lambda/\mu }^{(a+1-w,a+k-w)} \right).  
\end{equation}

By comparing the number of rows and columns, as well as the number of boxes of
the skew diagram $\lambda/\mu$ with their counterparts for  $\lambda$ it follows that $\lambda/\mu$ is $C'$-balanced with 
\[    C'=  C \sqrt{\frac{|\lambda|}{|\lambda|-w}} < \frac{C}{\sqrt{\mydelta}}. \] 
A fortiori $\lambda$ and $\lambda/\mu$ are $C''$-balanced, 
where $C''$ is the right-hand side of the above inequality.

\medskip

We apply \cref{lem:Pieriprob} twice: for both expressions in the square brackets
on the right-hand side of \eqref{eq:lolcats2}.  It follows that there exists a universal
constant $c>0$ (which depends only on $C''$) such that
\begin{multline*}
    \left| \mathbb{P}_{\widetildetableaux_{\lambda}^{(a+1,a+k)}} \! \left( M\big|_{\leq w} = S
    \right)  
    -
    \mathbb{P}_{\tableaux_{\lambda}} \! \left( M\big|_{\leq w} =
    S \right) \right|
    \leq  \\
    \leq 
    \mathbb{P}_{\widetildetableaux_{\lambda}^{(a+1,a+k)}} \! \left( M\big|_{\leq w} = S
    \right)  \frac{c k^2}{\sqrt{|\lambda|}}   
    +      
    \mathbb{P}_{\tableaux_{\lambda}} \! \left( M\big|_{\leq w} =
    S \right) \frac{c k^2}{\sqrt{|\lambda|-w}} 
   .
\end{multline*}
By summing over all choices of $S$ we get~\eqref{eq:TVD-probdivB} for $d:=\max(1,2c)$,
as required.
\end{proof}

\subsection{Proof of \cref{prop:single-and-multi}}
\label{sec:proof-single-and-multi}

\begin{proof}[Proof of \cref{prop:single-and-multi}]

We will sample the random tableaux $\mathbf{T}$ and $\mathbf{M}$ by the
following two-step procedure. Firstly, we sample $\mathbf{T}$ with the uniform
probability measure on $\tableaux_\lambda$, that is
$\mathbb{P} \! \left( \mathbf{T} = T \right) := \mathbb{P}_{\tableaux_{\lambda}}( T )$ 
for any $T \in \tableaux_\lambda$.
After the tableau $\mathbf{T}$ was
selected, we sample the tableau $\mathbf{M}$ with the conditional probability
\begin{multline*} \mathbb{P}\left( \cdot \;\big|\; \mathbf{T}=T \right) :=\\
\mathbb{P}_{\widetildetableaux_{\lambda}^{(a+1,a+k)}}\! \left( \cdot \;\Big|\; \! \left\{M\in \widetildetableaux_{\lambda}^{(a+1,a+k)} : M|_{\leq w}=\mathbf{T}|_{\leq w} \right\}\right).
\end{multline*}
In this way the condition \ref{item:M=T} is fulfilled trivially.

Therefore the probability distribution of $\mathbf{M}$ is given by
\begin{multline} 
    \label{eq:chainrule-A}
 \mathbb{P}\left( \mathbf{M}=M \right)= \\
\shoveleft{\mathbb{P}_{\tableaux_{\lambda}} \! \left\{ T\in \tableaux_{\lambda} : T|_{\leq w}= M|_{\leq w} \right\} \times} \\
{\mathbb{P}_{\widetildetableaux_{\lambda}^{(a+1,a+k)}} \! \left( M \;\Big|\; \! \left\{T\in \widetildetableaux_{\lambda}^{(a+1,a+k)} : T|_{\leq w}=M|_{\leq w} \right\}\right)}.
\end{multline}
The probability measure $\mathbb{P}_{\widetildetableaux_{\lambda}^{(a+1,a+k)}}$ can be written in an analogous way as
\begin{multline} 
        \label{eq:chainrule-B}
    \mathbb{P}_{\widetildetableaux_{\lambda}^{(a+1,a+k)}}(M)= \\
    \shoveleft{\mathbb{P}_{\widetildetableaux_{\lambda}^{(a+1,a+k)}} \! \left\{ T\in \widetildetableaux_{\lambda}^{(a+1,a+k)} : T|_{\leq w}= M|_{\leq w} \right\} \times} \\
    {\mathbb{P}_{\widetildetableaux_{\lambda}^{(a+1,a+k)}}\! \left( M \;\big|\; \! \left\{T\in \widetildetableaux_{\lambda}^{(a+1,a+k)} : T|_{\leq w}=M|_{\leq w} \right\}\right)}.
\end{multline}
It follows that the process of sampling $\mathbf{M}$ as well as 
the process of sampling the uniformly random
element of $\widetildetableaux_{\lambda}^{(a+1,a+k)}$ can be viewed as 
a two-step procedure: we first sample the positions of the boxes $1,\dots,w$ 
and
in the second step the remaining boxes. 
Notice that in both sampling procedures in the second step 
we sample the remaining boxes with the same conditional distribution. 
It follows that the total variation
distance between the measures \eqref{eq:chainrule-A} and
\eqref{eq:chainrule-B} is bounded from above by 
the total variation distance~\eqref{eq:TVD-probdivB} from \cref{lem:probdivB} 
which completes the proof.
\end{proof}

\section{The distribution of the~$u$-coordinates of the~multisurfers}
\label{sec:udist}

Our main result in this section is \cref{prop:longitude-experimental} which
shows that the multisurfer empirical measure
(i.e., the distribution of the~$u$-coordinates of the~multisurfers)
after draining a $(1-\alpha)$-fraction of the water converges to the limit measure
$\cotmeas_{\alpha}$ (the limit distribution of the $u$-coordinate of
the~single surfer on the level curve $g_{\alpha}$, cf.~\cref{sec:random-position}). Also \cref{prop:expvalue}
might be interesting from the viewpoint of algebraic combinatorics as it
provides a direct link between the statistical properties of uniformly random
$(a,b)$-Pieri tableaux of some fixed shape and the celebrated Jucys--Murphy
elements.

The following assumption on the `amount of water' and the number of multisurfers
will be central in the forthcoming results (\cref{prop:longitude-experimental,thm:moments,prop:dynamicalempiricallongitude}). 

\begin{assumption}
\label{ass:multisurfers}
Let $c > 0$ be the constant in \cref{lem:Pieriprob} obtained for $C=1$
for which \eqref{eq:pieri-proba} holds. 
We assume that $k=k(N)$ and $w=w(N)$ 
are sequences of positive integers 
which fulfill 
\[
\lim_{N \to \infty} k(N) = \infty
\quad \text{and} \quad 
k <  \sqrt{\frac{N}{2 c}} 
\quad \text{and} \quad 
w+k < N^2. 
\]
\end{assumption}

\subsection{Counting multisurfers gives the longitude}
\label{multiplesketch-0}

Let $w=w(N)$ and $k=k(N)$
be sequences of nonnegative integers such
that $0 < w+k < N^2$.
Let $\multi_N$ be a uniformly random tableau from $\widetildetableaux_{\sq_N}^{(w+1,w+k)}$.
We use a~shorthand notation $\content_n := \content_{n}^{\multi_N}$
for the $u$-coordinate of the box with the number $n$ in the tableau $\multi_N$. 
For $\content\in\R$ we define the random variable
$\empiricallongitude(\content)$ 
to~be~the~fraction of the multisurfers 
which have their scaled \mbox{$u$-coordinate} smaller than~$u$,~that is
\begin{equation} 
\label{eq:psi-via-multi}	
	\empiricallongitude(u) := \frac{1}{k} \max\left\{ p\in\{1,\dots,k\}  : 
	\frac{1}{N} \content_{\water+p} \leq u \right\}.  
\end{equation}
Clearly, $\empiricallongitude$ is the cumulative distribution function of
the~random measure~$\empiricalmeasure$ on~$\R$
\[\empiricalmeasure := \frac{1}{k} \sum_{1\leq p\leq k} \delta_{N^{-1} \, \content_{w+p}}\]
where $\delta_x$ denotes the \emph{delta measure concentrated at $x$}. 

\begin{theorem}
\label{prop:longitude-experimental} 
Let $\alpha \in (0,1)$. Let $w = w(N)$ and $k=k(N)$ fulfill 
\cref{ass:multisurfers} and 
\[
\lim_{N \to \infty} \frac{w}{N^2} = \alpha.
\] 
Then for each $\varepsilon>0$
\[ 
\mathbb{P}_{\widetildetableaux_{\sq_N}^{(w+1,w+k)}}
\bigg\{ 
\multi_N
 :
\sup_{u \in\R}
\left| F_{\cotmeas_{\alpha}}(x) - \empiricallongitude(x) \right| > \varepsilon  
\bigg\} 
 =  O\!\left(\frac{1}{k} + \frac{k^2}{N}\right)
\]
with the constant in the $O$-notation depending only on $\alpha$ and $\varepsilon$.

		In particular, if $k(N) \to \infty$ and $k(N) = o(\sqrt{N})$, i.e.,
		\eqref{eq:minimal-assumption-k} is satisfied,
		then the random sequence of 
		cumulative distribution functions 
		$\empiricallongitude$ 
		with respect to the supremum norm 
		converges in probability to
		the cumulative distribution function 
		$F_{\cotmeas_{\alpha}}$
		when $N\to \infty$, 
		that is shortly, 
		\[
		\sup_{x\in\R} \left| F_{\cotmeas_{\alpha}}(x) - \empiricallongitude(x) \right|
		\xrightarrow{P} 0. 
		\]
\end{theorem}

In order to prove this result we will compare the (random) moments of the
empirical measure $\empiricalmeasure$ and the moments of the measure
$\cotmeas_\alpha$ which gives the asymptotics of the $u$-coordinate of a single
box, cf.~\cref{sec:random-position}. In \cref{thm:moments} below we shall
calculate the~moments of the~empirical measure $\empiricalmeasure$. In
\cref{sec:proof-prop:longitude-experimental} we will complete the proof of
\cref{prop:longitude-experimental}.

\subsection{Moments of the empirical measure $\empiricallongitude$}
\label{sec:empirical-moments}

For each $\beta \in \N$ we define 
the~\emph{$\beta$-th moment of the random measure $\empiricalmeasure$} as
\[
M_\beta := M_{\beta}(w, k) := \int_{\R} z^\beta \dif  \empiricalmeasure(z) 
= \frac{1}{k} N^{-\beta} \sum_{1\leq p\leq k} u_{w+p}^\beta.
\]
Notice that
$M_\beta$ is also a random variable. The following result expresses the first
two moments of the random variable $ M_\beta(w,k)$ (which is related to the
problem of multisurfers) in terms of the first moment of the random variable
$M_{\beta}(w,1)$ (which is related to the much simpler problem of a single
surfer).
This proposition is crucial to the proof of \cref{prop:longitude-experimental}.

\begin{proposition}
    \label{thm:moments}
		Let $w = w(N)$ and $k=k(N)$ fulfill \cref{ass:multisurfers}. 
For each $\beta \in \N$, we have
\begin{align}
\label{EM}
\E_{\Pt}  M_\beta(w, k) &= \E_{\Pp}  M_\beta(w, 1) + O\!\left(\frac{k^2}{N}\right); 
 \\
\label{VarM}
\on{Var}_{\Pt} \! M_\beta(w, k) &= 
 O\!\left(\frac{1}{k} + \frac{k^2}{N}\right)
\end{align}
with the constants in the $O$-notation depending only on $\beta$.
\end{proposition}

\begin{remark}
The counterpart of the above proposition 
in the paper of Romik and the second named author is \cite[Theorem~4.6]{Romik2015a}.
There the error terms for the expected value and variance are much smaller, accordingly, 
$O\!\left(\frac{k}{N}\right)$ and $O\!\left(\frac{1}{k} + \frac{k}{N}\right)$.
There are two reasons for which our error terms are much bigger, of the form, accordingly, 
$O\!\left(\frac{k^2}{N}\right)$ and $O\!\left(\frac{1}{k} + \frac{k^2}{N}\right)$.
\begin{itemize}
\item In the proof of \cref{prop:expvalue} we will view the probability
distribution of the multisurfers as a \emph{conditional distribution} of the
boxes with certain numbers in a uniformly random standard skew tableau with a
specified shape \emph{under the condition} that these boxes are suitably ordered (i.e., the tableau is Pieri). Unfortunately, the probability of the latter event depends heavily on the shape of the diagram 
and we do not have a very good control over the error term, see \cref{lem:Pieriprob}.

Using our terminology in their context, the placement of the multisurfers by
Romik and the second named author can also be seen as a conditional process: one
first adds $k$ boxes to a given Young diagram~$\lambda$ by $k$ independent steps of the Plancherel growth
process, and then \emph{conditions} that these boxes are suitably ordered. In this case,
however, the conditioning does not create additional difficulties because the
probability of the event that the newly created boxes are Pieri is equal to
$\frac{1}{k!}$ and does not depend on the shape of $\lambda$.

\item The error $O\!\left(\frac{k^2}{N}\right)$ also appears during the application of \cref{prop:cosetcharacter}.
		Romik and the second named author make use of \cite[Theorem~4.4]{Romik2015a} 
		which is a counterpart of ours \cref{prop:expvalue}.
		They deal with the character of 
		the left-regular representation which obviously is not troublesome
		and need not be estimated.
\end{itemize}
\end{remark}

The proof of \cref{thm:moments} is quite long. 
In~\cref{sub:connection,sub:cosets,sub:Jucys-Murphy-product} 
we gather some tools helpful in proving \cref{thm:moments}. 
In particular, the~goal of~\cref{sub:connection} is to provide 
a~connection between the statistical properties of the multisurfers 
and the~representation theory (see~\cref{prop:expvalue}).
\cref{sub:cosets} gives background for calculating 
the~character $\chi_{_{\sq_N}}$ of the cosets appearing in~\eqref{eq:exp-a}. 
\cref{sub:Jucys-Murphy-product} 
is devoted mostly to an~analysis of 
the permutations arising from the~powers of 
Jucys--Murphy elements.   
Eventually, in~\cref{sec:expected-value} we give the~proof of~\eqref{EM}
and in~\cref{sec:variation} the~proof of~\eqref{VarM}
which completes the proof of \cref{thm:moments}.

\subsection{Multisurfers and the~representation theory}
\label{sub:connection}

The following result, \cref{prop:expvalue}, provides a~link 
between the~statistical properties of
the~multisurfers and~the~representation theory of the~symmetric groups.

\medskip

Let $w,k,n$ be positive integers such that $w+k\leq n$. 
With a small abuse of notation we denote by $\Sym_w$ the group of
permutations of~the~set $\{1,\dots,w\}$ and by $\Sym_k$ 
the group of permutations of~the~set $\{w+1,\dots,w+k\}$. 
In this way $\Sym_w$ and $\Sym_k$ are 
commuting subgroups of $\Sym_{w+k}\subset \Sym_{n}$.
We define the element of the symmetric group
algebra
\[ p_{\Sym_k} = \frac{1}{k!} \sum_{\sigma \in \Sym_k} \sigma \in
\C \Sym_{n}.\] 

Recall from \cref{subsec:Jucys--Murphy-elements-introduction} that the Jucys--Murphy elements $J_1,\dots,J_n$
form a commuting family in the symmetric group algebra $\C \Sym_n$
and are given by 
\[
    J_k := \sum_{1\leq i<k} (i, k) = (1, k) + \dots + (k-1, k) \in \C \Sym_n \quad \text{for } 1\leq k\leq n.
\]

\newcommand{\diagramA}{\lambda}
\newcommand{\diagramB}{\mu}
\newcommand{\diagramC}{\nu}

\begin{proposition}
    \label{prop:expvalue} 
    Let $w,k,n$ be positive integers such that $w+k\leq n$. 
		Let $W(x_1, \dots, x_k)$ be a symmetric polynomial in $k$ variables. 
    Let $\diagramA\in\Y_n$ be a Young diagram and $\tab$ be a random
    element (sampled with the uniform distribution) of the set
    $\widetildetableaux_{\diagramA}^{(w+1,w+k)}$ 
		of $(w+1,w+k)$-Pieri tableaux of shape~$\diagramA$. 
		Then
    \begin{align} 
        \E\ W \! \left(\content^\tab_{w+1},\ \dots,\ \content^\tab_{w+k}\right) 
        &= \frac{ \chi_{_{\diagramA}}\Big(W\!\left( J_{w+1}, \dots, J_{w+k}\right) \cdot p_{\Sym_k} \Big)}
        {\chi_{_{\diagramA}} \! \left(p_{\Sym_k} \right) } 
				\label{eq:exp-a}
        \\
        &= \frac{\chi_{_{\diagramA}}\Big(W\!\left( J_{w+1}, \dots, J_{w+k}\right) \cdot p_{\Sym_k} \Big)}
				{\mathbb{P}_{\tableaux_\diagramA}\! \left(\widetildetableaux_{\diagramA}^{(w+1,w+k)}\right)}.
				\nonumber
    \end{align}
\end{proposition}

The proof is postponed until the end of the current section until we gather the necessary tools.

\medskip 

We start with the following fundamental property of Jucys--Murphy elements.
\begin{fact}[{\cite{Jucys1974}}]
    \label{lem:Jucys}
Let $\lambda \in \Y_n$ be a Young diagram, and let $\content_1, \dots, \content_n$ be 
the $u$-coordinates of its boxes (listed in an~arbitrary order). 
Let $W(x_1, \dots, x_n)$ be a symmetric polynomial in 
$n$ variables. Then:
\begin{itemize}
\item $W(J_1, \dots, J_n) \in \C \Sym_n$ belongs to the center of the group 
algebra.
\item The operator $\rho_\lambda\big(W(J_1,\dots,J_n)\big)$ is a multiple 
of the identity operator, so it can be identified with a complex number. 
The value of this number is equal to
\begin{equation}\label{jmlemma}
\chi_{_{\lambda}} \Big(W(J_1, \dots, J_n) \Big) = W(\content_1, \dots, \content_n).
\end{equation}
\end{itemize}
\end{fact}

\newcommand{\niagramA}{\mu}
\newcommand{\niagramB}{\nu}

\begin{lemma}
    \label{lem:act}
Let $\niagramA \in \Y_{w+k}$ be a Young diagram 
and let $W(x_1, \dots, x_k)$ be a symmetric polynomial in 
$k$ variables. 
Then the operator 
\[ \rho_\niagramA\big(W(J_{w+1},\dots,J_{w+k})\big)\] 
acts on each irreducible component $V_\niagramB$ of the~restriction~$V_\niagramA \big\downarrow^{\Sym_{w+k}}_{\Sym_w}$  
as a multiple 
of the identity operator 
and can be identified with the~complex number
\begin{equation}\label{eq:jmk}
W\!\left(\content_{w+1}^{\niagramA / \niagramB}, \dots,
\content_{w+k}^{\niagramA / \niagramB} \right).
\end{equation}
Above, for a diagram $\niagramB$ with $w$ boxes such that $\niagramB \subset \niagramA$ we denote by
$\content_{w+1}^{\niagramA / \niagramB}, \dots, \content_{w+k}^{\niagramA
    / \niagramB}$ the~$u$-coordinates of the boxes of the~skew diagram $\niagramA / \niagramB$ (listed in an~arbitrary order).
\end{lemma}

\begin{proof}
 
We will show that the lemma holds 
in the particular case when 
$W$ is the~\emph{power-sum symmetric function}, that is 
\[W(x_1, \dots, x_k) := p_\beta(x_1, \dots, x_k) := \sum_{1\leq i\leq k} x_i^\beta\] 
for some $\beta \in \N_0$. Since the power-sum symmetric functions generate
the~algebra of the~symmetric polynomials and the representation $\rho_\niagramA$ is an algebra
homomorphism, in this way we~will prove that the lemma holds true in general.

Clearly,
\begin{align*}
W \! \left(J_{w+1}, \dots, J_{w+k} \right) &= p_\beta\left(J_{w+1}, \dots, J_{w+k} \right) = \\
&= p_\beta \left(J_1, \dots, J_{w+k} \right) - p_\beta \left(J_1, \dots, J_w \right).
\end{align*}
By~\cref{lem:Jucys}, the operator $p_\beta \left(J_1, \dots, J_{w+k} \right)$ 
acts on the component $V_\niagramA$ 
as multiplication by the~factor 
\begin{equation}
\label{eq:star}
p_\beta \left(\content^\niagramA_1, \dots, \content^\niagramA_{w+k} \right) = 
\sum_{1\leq i\leq w+k} \left( \content^\niagramA_i \right)^\beta. 
\end{equation}
Again by~\cref{lem:Jucys}, 
for any $\niagramB \in \Y_w$,
the~operator $p_\beta \left(J_1, \dots, J_{w} \right)$ 
acts on the component $V_\niagramB$  
as multiplication by the~factor 
\begin{equation}
\label{eq:starstar}
p_\beta \left(\content^\niagramB_1, \dots, \content^\niagramB_{w} \right) = \sum_{1\leq i\leq w} \left( \content^\niagramB_i \right)^\beta. 
\end{equation} 

Let $\niagramB \subset \niagramA$. The multiset of the $u$-coordinates of the~boxes
of~$\niagramA$ is the union of (i) the multiset of the $u$-coordinates of the~boxes
of~$\niagramB$, and (ii) the multiset of the $u$-coordinates of the~boxes of~$\niagramA
/ \niagramB$. Therefore, by subtracting \eqref{eq:starstar} from
\eqref{eq:star}, we get that the~operator 
$p_\beta \left(J_{w+1}, \dots, J_{w+k} \right)$ 
acts on the component $V_\niagramB$ of the~restriction $V_\niagramA
\downarrow^{\Sym_{w+k}}_{\Sym_w}$ as multiplication by the~scalar
\[
p_\beta \left(\content^{\niagramA / \niagramB}_{w+1}, \dots,
\content^{\niagramA / \niagramB}_{w+k} \right) = \sum_{1\leq i\leq k}
\left( \content^{\niagramA / \mu}_{w+i} \right)^\beta,
\]
as required.
\end{proof}

\begin{proof}[Proof of \cref{prop:expvalue}]
    Observe that any tableau $\tab \in
\widetildetableaux_{\diagramA}^{(w+1,w+k)}$ can be split into the
following three parts: (i)~a~standard tableau $P$ with entries from
$\{1,\dots,w\}$; we denote its shape by $\diagramC$,\; (ii)~a~skew tableau
$Q$ which is $k$-Pieri with entries in $\{w+1,\dots,w+k\}$; we denote its
shape by $\diagramB/\diagramC$, and\; (iii)~a~skew tableau $R$
with~the~entries~$>w+k$ with shape $\diagramA/ \diagramB$.
    
    For fixed partitions $\diagramB$ and $\diagramC$ it is easy to count the
number of tableaux~$P$ which contribute to~(i) and the number of tableaux
$R$ which contribute to~(iii): their cardinalities are by definition given
by $f^\diagramC$ and $f^{\diagramA/ \diagramB}$ respectively.
    \newline
		The number of $k$-Pieri tableaux $Q$ which contribute to (ii) is slightly
more challenging: it is equal to $1$ if $\diagramB/ \diagramC$ has at most
one box in each column and is equal to zero otherwise. A combinatorial
interpretation of the Littlewood--Richardson coefficient
$c^{\diagramB}_{\diagramC,(k)}$ for a single-row partition $(k)$ (or,
nomen omen, the Pieri rule) implies that this coefficient coincides with the
latter cardinality.

    In this way we proved that the left-hand side of~\eqref{eq:exp-a} is given by
    \begin{multline}
        \label{eq:ciagE}
        \E\ W\!\left(\content^\tab_{w+1},\ \dots,\ \content^\tab_{w+k}\right) = \\
        \frac{1}{\left|\widetildetableaux_{\diagramA}^{(w+1,w+k)}\right|}
				\sum_{\tab \in \widetildetableaux_{\diagramA}^{(w+1,w+k)}}  
				W\!\left(\content^\tab_{w+1},\ \dots,\ \content^\tab_{w+k}\right) = 
        \\
        \frac{1}{\left|\widetildetableaux_{\diagramA}^{(w+1,w+k)}\right|} \sum_{\substack{ \diagramB\in\Y_{w+k} \\ \diagramB \subset \diagramA}}\  
        \sum_{\substack{\diagramC \in \Y_w \\ \diagramC\subset\diagramB}} f^{\diagramA / 
            \diagramB} \, f^\diagramC \, c_{\diagramC, (k)}^\diagramB \, W\!\left(\content^{\diagramB/ 
            \diagramC}_{w+1}, \dots, \content^{\diagramB/ \diagramC}_{w+k} \right).
    \end{multline}

    \bigskip 
    
    We will now investigate the numerator on the~right hand side
of~\eqref{eq:exp-a}. 
	Each multiplicity $f^{\diagramA/\diagramB}$ in the decomposition of the
	restriction of $V_\lambda$ into irreducible components
			\[
	V_{\diagramA}\big\downarrow^{\Sym_{N^2}}_{\Sym_{w+k}} = 
	\bigoplus_{\substack{\diagramB\in \Y_w \\ \diagramB\subset \diagramA}} 
	f^{\diagramA / \diagramB} V_\diagramB
	\]
	is equal to the number of skew standard Young tableaux of shape $\diagramA/\diagramB$. 
Since \mbox{$W\!\left(J_{w+1},\ \dots,\ J_{w+k}\right)
    \cdot p_{\Sym_k} \in \C \Sym_{w+k}$}, we get that
    \begin{multline}
        \label{eq:numerator}
            C_{\diagramA} := \chi_{_{\diagramA}}
						\Big(W \left(J_{w+1},\ \dots,\ J_{w+k}\right) \cdot p_{\Sym_k} \Big) \\
            = \frac{1}{\on{dim} V_{\diagramA}} \on{Tr}_{{V_{\diagramA}}} 
						\rho_{{\diagramA}}  \Big(W \left(J_{w+1}, \dots, J_{w+k}\right) 
						\cdot p_{\Sym_k} \Big) = \\
            = \frac{1}{f^{\diagramA}} \sum_{\substack{\diagramB \in\Y_{w+k} \\ \diagramB \subset \diagramA}} 
            f^{\diagramA / \diagramB} \on{Tr}_{V_{\diagramB}} 
						\rho_{{\diagramA}}  \Big(W \left(J_{w+1}, \dots, J_{w+k}\right) 
						\cdot p_{\Sym_k} \Big).
\end{multline}

		The multiplicities in the decomposition of the restricted representation into irreducible components  
    \[
    V_{\diagramB}\big\downarrow^{\Sym_{w+k}}_{\Sym_{w} \times \Sym_{k}} 
    = \bigoplus_{\substack{ \diagramC \in\Y_{w} \\ \xi\in\Y_{k}}} c_{\diagramC, \xi}^\diagramB V_\diagramC \otimes V_\xi
    \]
    are given by Littlewood--Richardson coefficients.
    Therefore, since $p_{\Sym_k}$ is a~projection to the~trivial representation, 
    its image is given by
    \begin{equation}
        \label{eq:psk}
        p_{\Sym_k} \! \left( V_{\diagramB} \right) 
        = p_{\Sym_k} \! \left( V_{\diagramB}\big\downarrow^{\Sym_{w+k}}_{\Sym_{w} \times \Sym_{k}} \right)
        = \bigoplus_{\diagramC\in\Y_{w}} c_{\diagramC, (k)}^\diagramB V_\diagramC \otimes V_{(k)} = \bigoplus_{\diagramC\in\Y_{w}} c_{\diagramC, (k)}^\diagramB V_\diagramC.
    \end{equation}
    
    By combining \eqref{eq:numerator}, \eqref{eq:psk}, and \cref{lem:act} we
get the following closed formula for the numerator on the right-hand side
of \eqref{eq:exp-a}
    \begin{equation}
        \label{eq:Csq}
        C_{\diagramA} = 
        \frac{1}{f^{\diagramA}} \sum_{\substack{\diagramB \in\Y_{w+k} \\ \diagramB \subset \diagramA}} \sum_{\diagramC\in\Y_{w} } 
        f^{\diagramA / \diagramB} \, c_{\diagramC, (k)}^\diagramB \, f^\diagramC \, W\!\left(\content_{w+1}^{\diagramB / \diagramC}, \dots, \content_{w+k}^{\diagramB / \diagramC}\right) .
    \end{equation}
    
    \medskip
    
    On the other hand, by~\eqref{eq:Csq} evaluated for the constant polynomial $W
    \equiv 1$, we get the following formula for the denominator on the right-hand
    side of~\eqref{eq:exp-a}
    \begin{equation}
        \label{eq:chp}
        \chi_{_{\diagramA}} \! \left(p_{\Sym_k} \right) = 
        \frac{1}{f^{\diagramA}} \sum_{\substack{\diagramB \in\Y_{w+k} \\ \diagramB \subset \diagramA}} \sum_{\diagramC\in\Y_{w} } 
        f^{\diagramA / \diagramB} c_{\diagramC, (k)}^\diagramB f^\diagramC  
        = \frac{1}{f^{\diagramA}} \cdot \left|\widetildetableaux_{\diagramA}^{(w+1,w+k)}\right|
    \end{equation}
    where the~last equality comes from~\eqref{eq:ciagE} evaluated for $W \equiv 1$. 
		Observe also that 
		\begin{equation}
		\label{eq:chi-prob}
		\mathbb{P}_{\tableaux_\diagramA}\! \left(\widetildetableaux_{\diagramA}^{(w+1,w+k)}\right)
		=	\frac{f^{\diagramA}}{\left|\widetildetableaux_{\diagramA}^{(w+1,w+k)}\right|}.
		\end{equation}
		
    \medskip
    
    Equations \eqref{eq:ciagE}, \eqref{eq:Csq}, \eqref{eq:chp} and \eqref{eq:chi-prob} complete the proof of \eqref{eq:exp-a}.		
\end{proof}

\subsection{Character on a coset} 
\label{sub:cosets} 

Let us call \emph{small} the~elements of the~set $\{1, \dots, \water\}$ and
\emph{big} the~elements of the~set $\{\water+1, \dots, \water+k\}$. 
As before, we view the symmetric group $\Sym_{k}$ as 
the subgroup of $\Sym_{\water+k}$ which consists of the permutations 
which can permute only the big elements.

For a Young diagram $\lambda\in\Y_n$ and a permutation $\pi\in\Sym_{w+k}$ we
define \emph{the value of the character $\chi_{_{\lambda}}$ on a left coset $\pi
    \Sym_k\in \Sym_{w+k}/\Sym_k$} as an appropriate sum over the coset, that is
\[ 
\chi_{_{\lambda}}\!\left( \pi \Sym_k\right) 
:=
\sum_{\sigma\in \pi \Sym_k} \chi_{_{\lambda}}(\sigma).
\]
This definition is motivated by \cref{prop:expvalue} because expressions of a
similar flavor (up to the factor $\frac{1}{k!}$) appear on the right-hand
side of \eqref{eq:exp-a}. 
Our goal in this section is to understand the asymptotics of such characters on cosets.

\medskip

For a left coset $\pi \Sym_k \in \Sym_{\water+k} / \Sym_k$ we define its \emph{length} as
\begin{equation}
    \label{eq:coset-length}
\left\| \pi \Sym_k \right\| 
:=  \water - \#(\text{cycles of $\pi$ which consist of only small elements}).
\end{equation}
It is easy to check that if $\pi_1 \Sym_k= \pi_2 \Sym_k$ then the cycles of $\pi_1$
which consist of only small elements coincide with the analogous cycles of
$\pi_2$; it follows that the above definition does not depend on the choice of
the representative~$\pi$ of the coset.

Remind that for a permutation $\pi$ we denote by $|\pi|$ 
\emph{the length of $\pi$}, i.e., the minimal number of transpositions
required to write $\pi$ as their product. 

\newcommand{\optimal}{minimal\space}

\begin{lemma}\label{lem:optimalpi0}
    For each left coset $\pi \Sym_k$ 
    \begin{equation} 
        \label{eq:minum-cosets}
        \left\| \pi \Sym_k \right\|  =
     \min\left\{ |\sigma| : \sigma\in \pi \Sym_k \right\}.
\end{equation}
    There exists a \emph{unique} permutation $\pi_0\in \pi \Sym_k$ for which the minimum
on the right-hand side is achieved; the permutation $\pi_0=\pi_0(\pi \Sym_k)$
with this property will be called \emph{\optimal} for the coset $\pi \Sym_k$.

This \optimal permutation has the following additional properties:
     
		\begin{enumerate}[label=(\alph*)]
		\item \label{lem:optimalpi0a}
        for each $\sigma \in \Sym_k$, 
		\begin{equation}\label{eq:pi0sigma}
		\left| \pi_0 \sigma  \right| =  \| \pi \Sym_k \| + \left| \sigma\right|. 
		\end{equation}
		\item \label{lem:optimalpi0b} 
		If $\pi$ is such that each of its cycles permutes at most one big element,
		then $\pi=\pi_0$ is the \optimal element.
    \end{enumerate}
\end{lemma}

\begin{proof}
We begin with the proof of the first assertion of the lemma. 
The permutation $\pi$ can be written as a product of disjoint cycles $\pi = \pi_1
\cdots \pi_\ell \pi_{\ell+1} \cdots \pi_L$, where $\pi_1,\dots,\pi_\ell$ are the
cycles which permute at least one big element, 
and $\pi_{\ell+1},\dots,\pi_L$ permute only small elements.

\medskip

Fix $i\in \{1, \dots, \ell\}$. Then $\pi_i$ is a~cycle of the~form
\[ \pi_i= \left( 
p_{1,1}, \dots, p_{1,r_1}, q_{1},\;\;\;
p_{2,1}, \dots, p_{2,r_2}, q_{2},\;\;\; \dots,\;\;\;
p_{n,1}, \dots, p_{n,r_n}, q_{n} \right),\]
for some $n \in \N$ and $r_1,\dots,r_n \in \N_0$, 
some big numbers $q_1,\dots,q_n$ and some small numbers $p_{r,s}$. Define
\[
\tau_i := \left(q_{n}, q_{n - 1}, \dots, q_{1} \right)\in\Sym_k
\]
as the cycle permuting (in the reverse direction) the big elements of the cycle~$\pi_i$.
 Clearly, 
\[
\pi_i \tau_i = \left(p_{1,1}, \dots, p_{1,r_1}, q_1 \right) 
\cdots 
\left(p_{n,1}, \dots, p_{n,r_n}, q_n \right)
\]
gives a product of disjoint cycles, 
each permuting exactly one big element.

\medskip

Then 
\begin{equation}
    \label{eq:cycle-decomposition} 
    \pi_0 :=\left( \pi_1
\tau_1\right) \cdots \left( \pi_\ell \tau_\ell\right) \pi_{\ell+1} \cdots
\pi_L\in \pi \Sym_k
\end{equation}
provides a decomposition into disjoint cycles which has the property that each
cycle permutes at most one big element. Hence for any $\sigma \in \Sym_k$
which permutes only big elements, the decomposition into disjoint cycles of
the product $\pi_0 \sigma\in\pi \Sym_k$ is obtained by merging appropriate
cycles of $\pi_0$, it follows therefore that
\begin{equation}
    \label{eq:pi0-sigma-length}
\left| \pi_0 \sigma  \right| =  \left|\pi_0 \right| + \left| \sigma\right|.
\end{equation} 
As a consequence, the minimum on the right-hand side of
\eqref{eq:minum-cosets} is achieved on $\pi_0$ and it is the unique
permutation with this property, as required.

\medskip

We shall now show that the equality~\eqref{eq:minum-cosets} holds true. 
It is enough to prove it in the special case when 
as the coset representative we take $\pi_0$
given by~\eqref{eq:cycle-decomposition}. In this special case
\eqref{eq:minum-cosets} is equivalent to
\begin{equation}
    \label{eq:pi0-length}
\left\| \pi_0 \Sym_k \right\|=
\left| \pi_0 \right| .
\end{equation}

The explicit decomposition into disjoint cycles \eqref{eq:cycle-decomposition} implies
that the left-hand side is equal to 
\[ w- (L-\ell).\]

On the other hand, \eqref{eq:cycle-decomposition} implies that $\pi_0$ has
exactly $L-\ell$ cycles which permute only small elements and $k$ additional
cycles, one for each big element. It follows that the  right-hand side of 
\eqref{eq:pi0-length} is equal to 
\[ |\pi_0|= (w+k) - \big( L-\ell + k \big)=w- (L-\ell)
 \] 
which concludes the proof of \eqref{eq:minum-cosets}.

\medskip

Property \ref{lem:optimalpi0a} is now a direct consequence of
\eqref{eq:pi0-sigma-length} and \eqref{eq:pi0-length}.

\smallskip

For the proof of property \ref{lem:optimalpi0b}, if $\pi$ is such that each of
its cycles permutes at most one big element then our construction gives
$\pi=\pi_0$ so $\pi$ is the \optimal representative of the coset.
\end{proof}

The next proposition gives an~insight to the irreducible characters corresponding
to square diagrams evaluated on left cosets.

\begin{proposition}
    \label{prop:cosetcharacter} 
For each positive integer $L$ there exists a constant $C_L$ with the following property. 

Let positive integers $w$ and $k$ be arbitrary and 
let $\pi_0\in\pi \Sym_k$
be the \optimal representative of a left coset $\pi \Sym_k\in \Sym_{\water+k} / \Sym_k$. 
If $\|\pi \Sym_k\|\leq L$ and  $N^2\geq w+k$ and $N>k^2$ 
then
\begin{align}
\label{eq:oszacpi0} 
\left| \chi_{_{\sq_N}}\!\left( \pi \Sym_k\right) 
- \chi_{_{\sq_N}}\!\left( \pi_0\right)  \right| & <  
\frac{C_L k^2}{N^{\|\pi \Sym_k\|+1}}, \\[2ex]
\label{eq:oszacpi1} 
\left| \chi_{_{\sq_N}}\!\left( \pi \Sym_k\right)   \right| & <  
\frac{2 C_L}{N^{\|\pi \Sym_k\|}}.
\end{align}
\end{proposition}

\begin{proof}
Let $d=d(L) \geq 1$ be a constant (which depends only on $L$) big enough 
so that it guarantees that 
\[ L+k    \leq d N. \]

By~\cref{thm:character-estimation}, \cref{lem:optimalpi0}\ref{lem:optimalpi0a} 
and \cref{lem:sumoversymmetric} we have
\begin{multline*}
\left| \chi_{_{\sq_N}}\!\left( \pi \Sym_k\right) - \chi_{_{\sq_N}}\!\left( \pi_0 \right) \right|
\leq \sum_{\substack{\sigma\in \Sym_k \\ \sigma\neq\operatorname{id}}}
\left| \chi_{_{\sq_N}}\!\left( \pi_0 \sigma \right)\right|
\leq 
\sum_{\substack{\sigma \in \Sym_k \\ \sigma \neq \on{id}}}
\left( \frac{a d}{N}  \right)^{|\pi_0 \sigma|} \leq \\
\left(\frac{a d}{N}\right)^{ \| \pi \Sym_k \|} 
\left[
\exp\left( \frac{a d k^2}{2 N} \right) - 1 \right]= \\
\frac{(a d)^{ \| \pi \Sym_k \|+1}}{2} 
\frac{k^2}{N^{ \| \pi \Sym_k \|+1}}  
\frac{
\exp\left( \frac{a d k^2}{2 N} \right) - 1 }{\frac{a d k^2}{2 N}} . 
\end{multline*}
Since $\frac{e^x - 1}{x}$ is a bounded function on the interval $\left(0,
\frac{a d}{2} \right]$, we get \eqref{eq:oszacpi0} as required.

By \eqref{eq:oszacpi0} and~\cref{thm:character-estimation} we get
\[
\left| \chi_{_{\sq_N}}\!\left( \pi \Sym_k\right)   \right|
< \frac{C_L k^2}{N^{\|\pi \Sym_k\|+1}}
+ \left( \frac{a d}{N} \right)^{|\pi_0|}
= \frac{1}{N^{\|\pi \Sym_k\|}} \left( C_L \frac{k^2}{N} + (a d)^{\|\pi \Sym_{k}\|} \right).
\]
It follows that we can increase the value of the constant $C_L$ in such a way that both \eqref{eq:oszacpi0} and \eqref{eq:oszacpi1} are fulfilled.
\end{proof}

\subsection{Products of Jucys--Murphy elements}
\label{sub:Jucys-Murphy-product}

\subsubsection{Set partitions}
\label{subsec:setpartition}

For calculations of the~moments \eqref{EM} and~\eqref{VarM}
we will need to better understand 
the~sum of products $\sum_{p=1}^k J_{\water+p}^\beta$.
We will use similar concepts and notions 
to the~ones from \cite[Section~4.9]{Romik2015a}.
Notice that 
\begin{equation}\label{eq:Jucys-product}
J_{\water+p}^\beta 
= \sum_{1 \leq j_1, \dots, j_\beta \leq \water+p-1} 
(j_1, \water+p) \cdots (j_\beta, \water+p).
\end{equation}
We denote the~summands contributing to the right-hand side of~\eqref{eq:Jucys-product} in the~following way. 
For any $p$ and a~sequence $j = (j_1, \dots, j_\beta) \in [\water+p-1]^\beta$ 
define 
\begin{equation}
    \label{eq:sigma}
\sigmapj := (j_1, \water+p) \cdots (j_\beta, \water+p). 
\end{equation}
We also  denote
\begin{align*}
Z_\Sigma \! \left(j \right) := \left\{ r \in \{1, \dots, \beta\}: j_r \leq \water \right\};\\
Z_\Pi \! \left(j \right) := \left\{ r \in \{1, \dots, \beta\}: j_r > \water \right\}. 
\end{align*}
The~sets $Z_\Sigma(j)$ and $Z_\Pi(j)$ indicate which
elements of the~sequence $j$ are, respectively, \emph{small} and \emph{big}. 
Notice that $Z_\Sigma \cup Z_\Pi = \{1, \dots, \beta\}$
and, since $Z_\Sigma(j)$ and $Z_\Pi(j)$ are disjoint,
\[
\left|Z_\Sigma(j) \right| + \left|Z_\Pi(j) \right| = \beta.
\]
  
We consider the~equivalence relation $\sim$ on~$Z_\Sigma(j)$ 
(respectively, on $Z_\Pi(j)$) given by
\[
m \sim n \iff j_m = j_n. 
\]
We denote by $\Sigma(j)$ and $\Pi(j)$ 
the~sets of the~equivalence classes of 
the~relation~$\sim$ on, respectively, $Z_\Sigma(j)$ and~$Z_\Pi(j)$. 
Then the numbers of the equivalence classes 
\begin{align*}
&\left|\Sigma(j) \right| = \left| \{j_1, \dots, j_\beta\} \cap  \{1, \dots, \water\} \right|, \\
&\left|\Pi(j) \right| =  \left| \{j_1, \dots, j_\beta\} \cap \left\{\water+1, \dots, \water+k \right\} \right|
\end{align*}
indicate how many different small / big elements appear in the~sequence~$j$.

\medskip

We say that $A$ is a set \emph{a set-partition of $X$}, and denote it by $A \in \SPar{X}$,
if  $A$ is a collection of disjoint non-empty subsets of $X$ 
such that $\bigcup A = X$. 

We call $(A, B)$
\emph{a pair of complementary set-partitions of $X$}
if $A \cup B$ is a set-partition of $X$
such that 
 \[
\left( \bigcup A \right)  \cap \left( \bigcup B \right) = \emptyset.
\]
This terminology may be a bit misleading since 
neither $A$ nor $B$ need to be a~set-partition of $X$.
Note also that we allow the situations when $A = \emptyset$ or $B= \emptyset$. 
For example, consider a sequence $j \in [\water+p-1]^\beta$, 
then the pair $\big(\Sigma(j), \Pi(j)\big)$ 
is a pair of complementary set-partitions of $[\beta]$. 

For a pair $(\Sigma,\Pi)$ of complementary set-partitions of $[\beta]$ 
we will say that
\emph{the sequence $j \in [\water+p-1]^\beta$ is of type $(\Sigma,\Pi)$} 
if $\Sigma=\Sigma(j)$ and $\Pi=\Pi(j)$. 
If this is the case, we will use a shorthand notation $j\in(\Sigma,\Pi)$.

\begin{lemma}
    \label{lem:permcosets}
Let $p_1, p_2 \in \{1, \dots, k\}$ and $s \in [\water+p_1-1]^\beta$, 
and $t \in [\water+p_2-1]^\beta$.
Suppose that 
$s$ and $t$ are of the same type.
Denote $\sigma_1 := \sigma_{p_1,s}$ and $\sigma_2 := \sigma_{p_2,t}$. 
Then
		\begin{enumerate}[label=(\alph*)]
	\item \label{lem:permcosets-a} There exists a permutation $g \in \Sym_\water \times \Sym_k$ 
	such that $g(\water+p_1) = \water+p_2$ 
	and $g(s_m) = t_m$ for $m\in\{1, \dots, \beta\}$;
	\item \label{lem:permcosets-b} Permutations $\sigma_1$ and $\sigma_2$ are conjugate by $g$,
	that is $\sigma_1 = g^{-1} \sigma_2 g$;
	\item \label{lem:permcosets-c} $\| \sigma_1 \Sym_k \| = \| \sigma_2 \Sym_k \|$.
	\end{enumerate}
\end{lemma}

\begin{proof}
Denote
\begin{align*}
A_{\leq \water} :=  [\water] \setminus \{s_m: m=1, \dots, \beta\}; \\
B_{\leq \water} :=  [\water] \setminus \{t_m: m=1, \dots, \beta\}.
\end{align*}
Since $\Sigma(s) = \Sigma(t)$, we have $\left| A_{\leq \water} \right| = \left| B_{\leq \water} \right|$,
so there exists a~bijection $\delta_1\colon A_{\leq \water} \to B_{\leq \water}$.
For the~same reason, there exists a~bijection $\delta_2\colon A_{>\water} \to B_{> \water}$ between 
analogously defined sets 
\begin{align*}
&A_{>\water} := \{\water+1, \dots, \water+k\} \setminus \left( \{s_m: m=1, \dots, \beta\} \cup \{\water+p_1\} \right); \\
&B_{>\water} :=  \{\water+1, \dots, \water+k\} \setminus \left( \{t_m: m=1, \dots, \beta\} \cup \{\water+p_2\} \right).
\end{align*}
Then $g\colon \Sym_{\water+k} \to \Sym_{\water+k}$ given by
\[
g(x) := \left\{ \begin{array}{ll}
\water+p_2 & \text{if $x = \water+p_1$}, \\
t_m & \text{if $x = s_m$ for some $m \in \{1, \dots, \beta\}$}, \\
\delta_1(x) & \text{if $x \in A_{\leq \water}$}, \\
\delta_2(x) & \text{if $x \in A_{> \water} $}.
\end{array} \right.
\]
clearly fulfills the properties required in \ref{lem:permcosets-a} and it is easy
to check that \ref{lem:permcosets-b} indeed holds true.

In order to prove \ref{lem:permcosets-c} it is enough to notice that \ref{lem:permcosets-b} implies that
$\sigma_1$ and $\sigma_2$ have the~same number of cycles
permuting only small elements.
\end{proof}

\subsubsection{Inequalities concerning the character of $J_{\water+p}^\beta$}

In the proof of \cref{thm:moments} 
we will deal with the numbers of the form $O(N^{x})$
where the exponent is given by the right-hand-side of \eqref{eq:exponents-estimation}.
Our next aim is to show that these exponents are always nonpositive and 
only in some special cases are equal to $0$. 

\begin{lemma}
\label{lem:Jucys-product}
Let $x_1,x_2,\ldots\in \{1,\dots,w+k\}$ and
$y_1,y_2,\ldots\in \{w+1,\dots,w+k\}$ be infinite sequences.
Define a~function $f: \N_0 \to \Z$ by 
\begin{multline}
\label{eq:exponents-estimation}
f(\ell) := 2 \#\left( \{x_1,\dots,x_\ell\} \cap \left\{1, \dots, \water\right\} \right) + \\
\#\left( \{x_1,\dots,x_\ell\} \cap \left\{\water+1, \dots, \water+k \right\} \right) - \\
\| (x_1, y_1) \cdots (x_\ell, y_\ell) \Sym_k\| 
- \ell.
\end{multline}
Then $f \leq 0$ and $f$ is weakly decreasing. 

\medskip

Moreover, suppose that $\ell\in \N$ is such that 
$x_{\ell+1} \in \{x_1, \dots, x_\ell\}$
and 
$x_{\ell+1}$ and $y_{\ell+1}$ belong to
different cycles in the~cycle decomposition of the~product
$(x_1, y_1)\cdots (x_\ell, y_\ell)$.
Then $f(\ell+1) < f(\ell)$ and, 
in particular, $f(n) < 0$ for all $n \geq \ell+1$.
\end{lemma}

\begin{proof}
Let $\ell$ be a non-negative integer; we will show that $f(\ell+1)\leq f(\ell)$.

\medskip

Consider first the case $x_{\ell+1} > \water$. Then
\[
(x_1, y_1) \cdots (x_{\ell+1}, y_{\ell+1}) \Sym_k
=  (x_1, y_1) \cdots (x_{\ell}, y_{\ell}) \Sym_k
\]
so 
\begin{equation}
    \label{eq:p-value}
f(\ell+1) = 
\begin{cases}
f(\ell)-1 & \text{if $x_{\ell+1} \in \{x_i: i \leq \ell\}$}, \\ 
f(\ell) & \text{otherwise},
\end{cases}
\end{equation}
thus $f(\ell+1)\leq f(\ell)$, as required.

\medskip

Assume now that $x_{\ell+1} \leq \water$. 
Our strategy is to compare the cycle decomposition of the products
\begin{equation} 
    \label{eq:products}
    (x_1, y_1) \cdots (x_{\ell}, y_{\ell}) \quad \text{ and } \quad  
(x_1, y_1) \cdots (x_{\ell+1}, y_{\ell+1}) 
\end{equation}
and to deduce in this way (via the~definition \eqref{eq:coset-length}) the
relationship between the coset lengths
\[
\| (x_1, y_1) \cdots (x_{\ell+1}, y_{\ell+1}) \Sym_k\| \quad \text{ and }
 \quad  \| (x_1, y_1) \cdots (x_\ell, y_\ell) \Sym_k\|.\]
 Consider the following two cases.
\begin{itemize}[itemsep=1ex]
\item Suppose $x_{\ell+1} \notin \{x_i: i\leq \ell\}$. 
Then the cycle decomposition of the permutation 
on the right-hand side of \eqref{eq:products} arises 
from its counterpart on the left-hand side 
by merging the fixpoint $x_{\ell+1}$
with the cycle which contains $y_{\ell+1}$. 
It follows that
\[
\| (x_1, y_1) \cdots (x_{\ell+1}, y_{\ell+1}) \Sym_k\|
= \| (x_1, y_1) \cdots (x_\ell, y_\ell) \Sym_k\| + 1
\]
hence $f(\ell+1) = f(\ell)$, as required.

\item Suppose $x_{\ell+1} \in \{x_i: i\leq \ell\}$. The cycle decomposition of
the right-hand side of \eqref{eq:products} is obtained from its counterpart on
the left-hand side either by merging two cycles or by splitting one cycle into
two cycles. Each of these two operations can change the number of cycles which
permute only small elements by at most $1$. It follows that
\[
\| (x_1, y_1) \cdots (x_{\ell+1}, y_{\ell+1}) \Sym_k\| \geq  
\| (x_1, y_1) \cdots (x_\ell, y_\ell) \Sym_k\|- 1 
\]
so $f(\ell+1)\leq f(\ell)$, as required.
\end{itemize}
This completes the proof that $f$ is weakly decreasing.

Since $f(0)=0$ it follows that $f\leq 0$ and the proof of the first part of the lemma is complete.

\medskip

We will prove the second part of the lemma by revisiting the above proof.

If $x_{\ell+1} > \water$ then by \eqref{eq:p-value}, $f(\ell) = f(\ell-1) - 1$, as required.

On the other hand, when $x_{\ell+1} \leq \water$ the cycle decomposition of the
right-hand side of \eqref{eq:products} is obtained from its counterpart on the
left-hand side by merging two cycles: the one containing $x_{\ell+1}$ with the
one containing $y_{\ell+1}$. Therefore the number of cycles which consist of
only small elements at the right-hand side of \eqref{eq:products}  is bounded
from above by its counterpart for the left-hand side of \eqref{eq:products}.
Hence
\[
\| (x_1, y_1) \cdots (x_{\ell+1}, y_{\ell+1}) \Sym_k\| \geq  
\| (x_1, y_1) \cdots (x_{\ell}, y_{\ell}) \Sym_k\|
\]
and so $f(\ell+1) \leq f(\ell) - 1$, as required.
\end{proof}

\begin{corollary}
\label{cor:Jucys-product}
Let $x_1,\ldots,x_\ell \in \{1,\dots,w+k\}$ and
$y_1,\ldots,y_\ell\in \{w+1,\dots,w+k\}$. Denote
\begin{align*}
|\Sigma| &:= \# \left(\{x_1,\dots,x_\ell\} \cap \left\{1, \dots, \water\right\} \right); \\
|\Pi| &:= \# \left(  \{x_1,\dots,x_\ell\} \cap \left\{\water+1, \dots, \water+k \right\} \right); \\
\| \sigma S_k \| &:= \|(x_1,y_1) \cdots (x_\ell,y_\ell) S_k\|.
\end{align*}
Then for $N\geq 1$
\[ 
N^{2\, |\Sigma|}\; k^{|\Pi|}\; N^{-\| \sigma S_k \| - \ell} 
\leq \left(\frac{k}{N} \right)^{|\Pi|}.
\]
\end{corollary}

\subsection{The mean value of $M_\beta$ -- the~proof of~\eqref{EM} }
\label{sec:expected-value}

\begin{proof}[Proof of \eqref{EM}]
Our goal is to calculate the~expected value of the~moment $ M_\beta(w,k)$
(recall~\cref{sec:empirical-moments}).
By~\cref{prop:expvalue},
\begin{multline}
    \label{eq:firstmoment}
\E_{\Pt} M_\beta\left(w, k \right) 
= \frac{1}{\Pp \! \left(\widetildetableaux_{\sq_N}^{(w+1,w+k)} \right)} \cdot \frac{1}{k} N^{-\beta} \chi_{_{\sq_N}}\! 
\left(\mathcal{J}  p_{\Sym_k} \right) = \\
= \frac{1}{k!\ \Pp \! \left( \widetildetableaux_{\sq_N}^{(w+1,w+k)} \right)} \cdot \frac{1}{k} N^{-\beta} \chi_{_{\sq_N}}\! 
\left(\mathcal{J}  \Sym_k \right)
\end{multline}	
where (recall the~definition~\eqref{eq:sigma} of $\sigmapj$)
\[ \mathcal{J}:=\sum_{p=1}^k J_{\water+p}^\beta =
\sum_{p=1}^k \sum_{j \in [\water+p-1]^\beta}  \sigmapj \in \C\Sym_{\water+k}.
\]
Since $\sq_N$ is $C$-balanced with $C=1$, by \cref{lem:Pieriprob}
the denominator on the right-hand side of \eqref{eq:firstmoment}
fulfills
\[ 
k!\ \Pp \! \left( \widetildetableaux_{\sq_N}^{(w+1,w+k)} \right) 
= 1 + O\!\left(\frac{k^2}{N}\right) 
\] 
with the constant in the $O$-notation equal to $c$ from \cref{lem:Pieriprob}.  
By \cref{ass:multisurfers}
the right hand side
is separated from $0$ and therefore
\begin{equation}
\label{eq:moments-denominator}
\frac{1}{k!\ \Pp \! \left( \widetildetableaux_{\sq_N}^{(w+1,w+k)} \right)} = 1 + O\!\left(\frac{k^2}{N}\right).
\end{equation}
\cref{prop:EM-asymptotics-c} from \cref{prop:EM-asymptotics} below 
provides the necessary asymptotics of the numerator and completes the proof
of \eqref{EM}.
\end{proof}

We consider an analogue of $\mathcal{J}$ in which each summand is replaced
by the minimal element of the appropriate coset
\[ 
\mathcal{J}_0:= \pi_0\left( \mathcal{J} \right) = 
\sum_{p=1}^k \sum_{j \in [\water+p-1]^\beta}  \pi_0\big( \sigmapj \big) \in \C\Sym_{\water+k}.
\]
The following proposition provides the missing element of the above proof of~\eqref{EM}.
\begin{proposition}
    \label{prop:EM-asymptotics} 
		Let $\beta \in \N$ be fixed. Let $k, N \in \N$ be such that
		\mbox{$N^2 > \water + k$} and $k^2 < N$.
		Then		
		\begin{align}
    \chi_{_{\sq_N}}\!\left( \mathcal{J}_0 \right) &= k N^\beta \! \left[ 
           \E_{\Pp} M_\beta (w,1)
           + O \left(\frac{k}{N} \right) 
           \right]; \label{prop:EM-asymptotics-a}    \\
		 \chi_{_{\sq_N}}\!  
       \left(\mathcal{J} \Sym_k \right) 
       &= \chi_{_{\sq_N}}\!\left( \mathcal{J}_0 \right)
       + O \! \left(k N^\beta  \frac{k^2}{N} \right);       
			\label{prop:EM-asymptotics-b} \\
		  \chi_{_{\sq_N}}\! 
     \left(\mathcal{J}  \Sym_k \right) &= k N^\beta \! \left[ 
    \E_{\Pp} M_\beta (w,1)
    + O \left(\frac{k^2}{N} \right) 
    \right]
			\label{prop:EM-asymptotics-c}
		\end{align}
with the constants in the $O$-notation depending only on $\beta$.    
\end{proposition}
The remaining part of this section is devoted to its proof.

\subsubsection{Decomposition of $\chi_{_{\sq_N}}\!\left( \mathcal{J}_0 \right)$}

Denote
\begin{align}
\nonumber
A & :=  \sum_{p=1}^k \sum_{j \in [\water]^\beta} \chi_{_{\sq_N}}\! \left(\pi_0\left(\sigmapj \right) \right); \\
\intertext{and for a pair of complementary set-partitions $(\Sigma, \Pi)$ of $[\beta]$ 
let us also denote}
\label{eq:definicja-bsigmapi}
B^{(\Sigma,\Pi)} & :=  \sum_{p=1}^k \sum_{j \in (\Sigma,\Pi)} \chi_{_{\sq_N}}\! \left(\pi_0\left(\sigmapj \right) \right).
\end{align}
With these notations
\begin{equation}
\label{eq:prop-a}
\chi_{_{\sq_N}}\!\left( \mathcal{J}_0 \right) = 
A + \sum_{\substack{(\Sigma,\Pi) \\ \Pi \neq \emptyset}} B^{(\Sigma,\Pi)}.
\end{equation}
Notice that the number of summands is finite, depends only on $\beta$ and does not depend on $N$ or $k$.
Therefore, it is enough to find the asymptotics for each individual summand on the right-hand side.

\subsubsection{Asymptotics of $A$}
\label{subsec:A-asymptotics}
By~\cref{lem:optimalpi0}\ref{lem:optimalpi0b}, $\pi_0 \left(\sigmapj \right)
= \sigmapj$ for all $j \in [w]^\beta$. 
Since the~character is constant on each conjugacy class,
the contribution of the summands in $A$ to the character is the~same for each value of $p$
and 
\[
A = \sum_{p=1}^k \chi_{_{\sq_N}}\! \left( \sum_{j \in [\water]^\beta} \sigmapj \right) 
= k \cdot \chi_{_{\sq_N}}\! \left(J_{\water+1}^\beta \right).
\]
We apply \cref{prop:expvalue} for $k=1$; in this special case $p_{\Sym_1}=\on{id}$ and 
$\chi_{_{\sq_N}}(p_{\Sym_1}) = 1$, thus
\begin{equation}\label{eq:EM-A}
A = k \: \E_{\Pp} u_{\water+1}^\beta 
= k N^\beta \: \E_{\Pp} M_\beta (w,1) .
\end{equation}

\subsubsection{Asymptotics of $B^{(\Sigma,\Pi)}$}
\label{subsec:B-asymptotics}

Let $j\in[w+p-1]^\beta$. By \cref{lem:optimalpi0} the coset length 
\[ \| \sigmapj \Sym_k \| \leq \| \sigmapj \| \leq \beta \]
is uniformly bounded from above by the number of factors in \eqref{eq:sigma}.
By~\cref{thm:character-estimation} and
\cref{lem:optimalpi0} it follows that there exists a universal constant $C_\beta$ such that 
\[ 
\left| \chi_{_{\sq_N}}\! \left(\pi_0\left(\sigmapj \right) \right) \right| 
\leq C_\beta N^{- \| \sigmapj \Sym_k \|} 
\]
holds true for each $j\in[w+p-1]^\beta$.

Let us fix a pair $(\Sigma,\Pi)$ of complementary set-partitions of $[\beta]$.
The number of the summands on the right-hand side of \eqref{eq:definicja-bsigmapi} 
is equal to the following sum of falling factorials
\[ \sum_{p=1}^k (w)_{|\Sigma|} 
\cdot (p)_{|\Pi|} \leq  N^{2 |\Sigma|}\ k^{|\Pi|+1}.\]
By combining these observations with \cref{cor:Jucys-product} we conclude that
\begin{equation}
    \label{eq:class-estimate}
    \left| B^{(\Sigma, \Pi)} \right|
    \leq    C_\beta N^{- \| \sigmapj \Sym_k \|}  \cdot N^{2 |\Sigma|}\ k^{|\Pi|+1} \leq
    C_\beta \, k N^\beta \! \left(\frac{k}{N} \right)^{\! |\Pi|}.
\end{equation}
The~last two arguments imply also that for any $p \in \{1, \dots, k\}$
\begin{equation}
\label{eq:EM-rest}
\sum_{j \in (\Sigma,\Pi)} \frac{k^2}{N^{\| \sigmapj \Sym_k\|+1}}
\leq \frac{k^2}{N} \, N^\beta \! \left(\frac{k}{N} \right)^{\! |\Pi|}.
\end{equation}

\subsubsection{Proof of \cref{prop:EM-asymptotics}}

\begin{proof}[Proof of \cref{prop:EM-asymptotics}]
The asymptotics of the summands which contribute to the right-hand side of
\eqref{eq:prop-a} is provided by the~equality \eqref{eq:EM-A} and the estimate
\eqref{eq:class-estimate} 
(for pairs $(\Sigma, \Pi)$ of complementary set-partitions of $[\beta]$ 
with $\Pi \neq \emptyset$). 
In this way the proof of \eqref{prop:EM-asymptotics-a} is complete.

\medskip

By~\cref{prop:cosetcharacter} there exists $H_\beta > 0$ such that 
\begin{equation}
    \label{eq:EM-est}
    \left|
    \chi_{_{\sq_N}}\!  
    \left(\mathcal{J} \Sym_k \right) 
    - \chi_{_{\sq_N}}\!\left( \mathcal{J}_0 \right)
    \right| 
    \leq H_\beta  \sum_{p=1}^k \sum_{j \in [\water+p-1]^\beta}  \frac{k^2}{N^{\| \sigmapj \Sym_k\|+1}} .
\end{equation}
The summands on the right-hand side can be grouped according to the type
$(\Sigma,\Pi)$ of the sequence $j$. For a fixed value of $\beta$ there are only finitely many possible types,
and the total contribution of the sequences $j$ of a specific type
$(\Sigma,\Pi)$ is bounded from above by \eqref{eq:EM-rest} which completes the~proof of
\eqref{prop:EM-asymptotics-b}.

\medskip

\cref{prop:EM-asymptotics-c}
is a direct consequence of 
\eqref{prop:EM-asymptotics-a} and \eqref{prop:EM-asymptotics-b}.
\end{proof}

\subsection{The variance of $M_\beta$ -- the~proof of~\eqref{VarM} }
\label{sec:variation}

We will mimic the concepts from \cref{sec:expected-value},
however, the calculations will be more involved.

\begin{proof}[Proof of \eqref{VarM}]
We first calculate the~second moment of $M_\beta$.
By~\cref{prop:expvalue} and then \eqref{eq:moments-denominator} 
\begin{multline}\label{eq:secondmoment}
\E_{\Pt}  M_\beta(w,k)^2 
= \frac{1}{k!\ \Pp \! \left(\widetildetableaux_{\sq_N}^{(w+1,w+k)} \right)} 
\cdot \frac{1}{k^2} N^{-2\beta} 
\chi_{_{\sq_N}}\! \left(\mathcal{J}^2 \Sym_k \right) 
\\
= \left( 1 + O\left(\frac{k^2}{N} \right) \right) \frac{1}{k^2} N^{-2\beta} 
\chi_{_{\sq_N}}\! \left(\mathcal{J}^2 \Sym_k \right)
\end{multline}
where (recall the~definition~\eqref{eq:sigma} of $\sigmapj$)
\[ \mathcal{J}^2:= \left( \sum_{p=1}^k J_{\water+p}^\beta \right)^2 =
\sum_{p_1, p_2=1}^k \sum_{\substack{s \in [\water+p_1-1]^\beta \\ t \in [\water+p_2-1]^\beta}} 
\sigma_{p_1,s} \sigma_{p_2,t} \in \C\Sym_{\water+k}.
\]

\cref{prop:VarM-asymptotics-c} from \cref{prop:VarM-asymptotics} below 
provides the necessary asymptotics of the numerator in \eqref{eq:secondmoment} 
and gives us 
\[
\E_{\Pt}  M_\beta(w,k)^2 
= \frac{1}{k} \E_{\Pp} M_{2\beta}(w, 1)
+ \left(1 - \frac{1}{k} \right) \left( \E_{\Pp} M_{\beta}(w,1) \right)^2 
+ O \! \left( \frac{k^2}{N} \right)
\]
with the constant in the $O$-notation depending only on $\beta$.
By~\eqref{EM} we finally get 
\begin{multline*}
\on{Var}_{\Pt} \! M_\beta(w, k) 
= \E_{\Pt} M_{\beta}(w,k)^2 - \left( \E_{\Pt} M_{\beta}(w,k) \right)^2
= \\
\frac{1}{k} \left[ \E_{\Pp} M_{2\beta}\left(w,1\right) 
- \left(\E_{\Pp} M_{\beta}\left(w,1\right) \right)^2 \right] 
+ O\!\left(\frac{k^2}{N}\right) 
= O\!\left(\frac{1}{k} + \frac{k^2}{N}\right)
\end{multline*}
with the constant in the $O$-notation depending only on $\beta$.
This completes the proof of \eqref{VarM} (and \cref{thm:moments}).
\end{proof} 

We consider an analogue of $\mathcal{J}^2$ in which each summand is replaced
by the minimal element of the appropriate coset
\[ 
\mathcal{J}^2_*:= \pi_0\left( \mathcal{J}^2 \right) = 
\sum_{p_1, p_2=1}^k \sum_{\substack{s \in [\water+p_1-1]^\beta \\ t \in [\water+p_2-1]^\beta}} 
\pi_0\left( \sigma_{p_1,s} \sigma_{p_2,t} \right) \in \C\Sym_{\water+k}.
\]

The following proposition provides the missing component of the above proof of~\eqref{VarM}.
\begin{proposition}
    \label{prop:VarM-asymptotics} 
		Let $\beta \in \N$ and $k, N \in \N$ be such that
		$N^2 \geq \water + k$ and $k^2 < N$.
		Then

			\begin{multline}
			\label{prop:VarM-asymptotics-a}
				\chi_{_{\sq_N}}\! \left( \mathcal{J}^2_* \right) = k^2 N^{2\beta} \Bigg[ 
				\frac{1}{k} \E_{\Pp} M_{2\beta} \left(w,1 \right) + \\ 
				\left(1-\frac{1}{k} \right) \left( \E_{\Pp} M_{\beta}(w,1) \right)^2 
				+ O \left(\frac{k}{N} \right) 
				\Bigg]; 
			\end{multline}
			
			\begin{multline}		
			\label{prop:VarM-asymptotics-b} 
				\chi_{_{\sq_N}}\!  
				 \left(\mathcal{J}^2 \Sym_k \right) 
				 = \chi_{_{\sq_N}}\!\left( \mathcal{J}^2_* \right)
				 + O \left( k^2 N^{2 \beta} \frac{k^2}{N} \right);   
					\hfill
			\end{multline}
			
			\begin{multline}		
			\label{prop:VarM-asymptotics-c}
				\chi_{_{\sq_N}}\! 
			 \left(\mathcal{J}^2  \Sym_k \right) = k^2 N^{2\beta} \Bigg[ 
				\frac{1}{k} \E_{\Pp} M_{2\beta} \left(w,1 \right) + \\
				\left(1-\frac{1}{k} \right) \left( \E_{\Pp} M_{\beta}(w,1) \right)^2 
				+ O \left(\frac{k^2}{N} \right) 
			\Bigg]		
			\end{multline}	
	with the constants in the $O$-notation depending only on $\beta$.
\end{proposition}

The remaining part of this section is devoted to its proof.

\subsubsection{Decomposition of $\chi_{_{\sq_N}}\!\left( \mathcal{J}^2_* \right)$}
For any $p_1, p_2 \in \{1, \dots, k\}$ denote 
\[
P_{p_1, p_2}(\beta) := [\water+p_1-1]^\beta \times [\water+p_2-1]^\beta
\]
and let 
\begin{align}
A	&:=	\sum_{p=1}^k \:
\sum_{s, t \in [w]^{\beta}} 
\chi_{_{\sq_N}}\! \left(\pi_0\left(\sigma_{p,s} \sigma_{p,t}  \right) \right); \nonumber
\end{align}

Let us define \emph{the~concatenation of sequences $a = (a_1, \dots, a_x)$ and $b = (b_1, \dots, b_y)$ }
as the~sequence $a \sqcup b := (a_1, \dots, a_x, b_1, \dots, b_y)$. 
For any pair of complementary set-partitions $(\Sigma, \Pi)$
of the~set $[2\beta]$ 
let us denote 
\begin{equation}
\label{eq:definicja-Bsigpi-VarM}
B^{(\Sigma, \Pi)} :=  \sum_{p_1, p_2=1}^k 
\sum_{\substack{ (s,t) \in P_{p_1, p_2}(\beta), \\  s \sqcup t \, \in \, (\Sigma, \Pi)}} 
\chi_{_{\sq_N}}\! \left(\pi_0\left(\sigma_{p_1,s} \sigma_{p_2,t}  \right) \right)
\end{equation}
and if $\Sigma$ is a set-partition of $[2\beta]$ denote
\begin{equation}
\label{eq:definicja-Csig-VarM}
C^{\Sigma}  := 
 \sum_{ \substack{ p_1, p_2 \in \{1, \dots, k\}, \\ p_1 \neq p_2 } } 
\sum_{\substack{ s, t \in [w]^{\beta}, \\ s \sqcup t \in (\Sigma, \emptyset) } }
\chi_{_{\sq_N}}\! \left(\pi_0\left(\sigma_{p_1,s} \sigma_{p_2,t}  \right) \right). 
\end{equation}
With these notations
\begin{equation}
\label{eq:prop-VarM}
\chi_{_{\sq_N}}\! \left( \mathcal{J}^2_* \right) = 
A + \sum_{\substack{(\Sigma, \Pi) \\ \Pi \neq \emptyset}} B^{(\Sigma, \Pi)}
+ \sum_{\substack{\Sigma}} C^{\Sigma}.
\end{equation}
Notice that the number of summands on the right hand side of \eqref{eq:prop-VarM}
is finite, depends only on $\beta$ and does not depend on $N$ or $k$.
Therefore, it is enough to find the asymptotics for each individual summand on the right-hand side.

\subsubsection{Asymptotics of $A$}
We proceed in the same way as in \cref{subsec:A-asymptotics}
to get an exact formula
\begin{equation}
\label{eq:VarM-A}
A = k N^{2\beta} \: \E_{\Pp} M_{2\beta}(w, 1).
\end{equation}

\subsubsection{Asymptotics of $B^{(\Sigma, \Pi)}$}
\label{subsec:relation-R}
We follow the lines from \cref{subsec:B-asymptotics}.

Let us fix some pair $(\Sigma, \Pi)$ of complementary set-partitions of $[2\beta]$. 
By \cref{thm:character-estimation} and
\cref{lem:optimalpi0} there exists $C_\beta > 0$ such that 
each summand corresponding to $(s,t) \in P_{p_1, p_2}(\beta)$ such that 
$s \sqcup t \in (\Sigma, \Pi)$ 
fulfills the asymptotic bound
\[ 
\left| \chi_{_{\sq_N}}\! \left(\pi_0\left(\sigma_{p_1,s} \sigma_{p_2,t}  \right) \right) \right| 
\leq C_\beta N^{- \| \sigma_{p_1,s} \sigma_{p_2,t}  \Sym_k \|}. 
\]
On the other hand, 
the number of the summands on the right-hand side 
of~\eqref{eq:definicja-Bsigpi-VarM} is bounded from above by  
\[ 
\sum_{p_1, p_2=1}^k w_{|\Sigma|} 
\cdot \max\{p_1, p_2\}_{|\Pi|} 
\leq  N^{2|\Sigma|}\ k^{|\Pi| + 2}.
\]
By combining these observations with \cref{cor:Jucys-product} used 
for the~concatenated sequence $s \sqcup t$ 
we conclude that
\begin{equation}
    \label{eq:B-class-estimate-VarM}
    \left| B^{(\Sigma, \Pi)} \right|
    \leq  
    C_\beta \: k^2 N^{2\beta} \left(\frac{k}{N} \right)^{|\Pi|}.
\end{equation}
The~last two arguments imply also that 
\begin{equation}
\label{eq:VarM-rest}
\sum_{p_1, p_2=1}^k 
\sum_{\substack{ (s,t) \in P_{p_1, p_2}(\beta), \\  s \sqcup t \, \in \, (\Sigma, \Pi)}} 
\frac{k^2}{N^{\| \sigma_{p_1,s} \sigma_{p_2,t}  \Sym_k\|+1}}
\leq \frac{k^2}{N} \cdot k^2 N^{2\beta} \left(\frac{k}{N} \right)^{|\Pi|}.
\end{equation}

These estimations show that if a pair $(\Sigma, \Pi)$ of complementary set-partitions of $[2\beta]$
is such that $\Pi \neq \emptyset$ then 
the contribution of $B^{(\Sigma, \Pi)}$ to \eqref{eq:prop-VarM} is of relatively small order.

\subsubsection{Asymptotics of $C^\Sigma$. Connected set-partitions of $[2\beta]$}
\label{subsec:C-calculations}
We will need to be much more subtle in calculating (and estimating) 
the summand $C^\Sigma$. 
We will treat the summand $C^\Sigma$ in two different ways
depending on the structure of $\Sigma$. 

Let $\Sigma$ be a set-partition of $[2\beta]$. 
We say that \emph{$\Sigma$ is connected} if there exists a block $\pi \in \Sigma$
which contains simultaneously some element of the set $\{1, \dots, \beta\}$ 
and some element of the set $\{\beta+1, \dots, 2\beta\}$,
that is formally, $\pi \cap [\beta] \neq \emptyset$ 
and $\pi \cap \{\beta+1, \dots, 2\beta\} \neq \emptyset$.
Each such a block $\pi \in \Sigma$ will be called  \emph{a link}.  
If $\Sigma$ does not have any links then we say that \emph{$\Sigma$ is disconnected}.

\subsubsection{Asymptotics of $C^{\Sigma}$ for connected $\Sigma$}
\label{subsec:asymptotics-E}
Let us fix a connected set-partition $\Sigma$ of $[2\beta]$. 
We set  
\[
n_{\Sigma} := \min\left\{ i \in \{\beta+1, \dots, 2\beta\}: \exists_{\pi \in \Sigma}\ 
\left( i \in \pi \text{ and $\pi$ is a link in $\Sigma$} \right) \right\}.
\]
In other words, $n_\Sigma$ indicates the least number $i > \beta$ 
belonging to some link $\pi$ of $\Sigma$.

Notice that for distinct $p_1, p_2$ and
a pair of sequences $(s,t)$ such that $s \sqcup t \in (\Sigma, \emptyset)$
the assumptions of the~second part of \cref{lem:Jucys-product}
are fulfilled for 
the~sequences $(x_i) = s \sqcup t$ and
$(y_i) = (p_1, \dots, p_1, p_2, \dots, p_2)$,
and $\ell := \beta + n_{\Sigma} - 1$.
Therefore an inequality
\[
2\left| \Sigma \right| 
- \| \sigma_{p_1,s} \sigma_{p_2,t} \Sym_k \| 
- 2\beta \leq -1
\]
holds for any $(s,t)$ such that $s \sqcup t \in (\Sigma, \emptyset)$
with $\Sigma$ connected.

We now follow the lines in \cref{subsec:relation-R} to
get the upper bound 
\begin{equation}
    \label{eq:C-class-connected-VarM}
    \left| C^{\Sigma} \right|
    \leq  
    C_\beta \: (k^2-k) N^{2\beta - 1}
\end{equation}
which holds for each connected set-partition $\Sigma$ of the set $[2\beta]$.
This shows that 
the contribution of $C^{\Sigma}$ with connected $\Sigma$
to \eqref{eq:prop-VarM} 
is of relatively small order. 

\subsubsection{Asymptotics of $C^\Sigma$ for disconnected $\Sigma$}
\label{subsec:asymptotics-D}

We will show that 
\begin{equation}
\label{eq:C-class-disconnected-VarM}
\sum_{\substack{\Sigma: \\ \Sigma \text{ is disconnected}}} C^{\Sigma} = (k^2-k) N^{2\beta} \! \left( \big( \E_{\Pp} M_{\beta}(w, 1) \big)^2 
+ O \! \left(N^{-2} \right) \right)
\end{equation}
with the constant in the $O$-notation depending only on $\beta$.

Recall that the $C^{\Sigma}$, defined in~\eqref{eq:definicja-Csig-VarM}, 
is the sum of characters $\chi_{_{\sq_N}}$ evaluated on the minimal permutations
$\pi_0\left(\sigma_{p_1,s} \sigma_{p_2,t}  \right)$.
When $\Sigma$ is disconnected and $p_1 \neq p_2$, 
in the permutation $\sigma_{p_1,s} \sigma_{p_2,t}$ 
each cycle permutes at most one big element,
so by \cref{lem:optimalpi0}\ref{lem:optimalpi0b} 
$\sigma_{p_1,s} \sigma_{p_2,t}$ is the minimal permutation. 
Hence whenever $\Sigma$ is disconnected 
\[
C^{\Sigma}  = 
 \sum_{ \substack{ p_1, p_2 \in \{1, \dots, k\}, \\ p_1 \neq p_2 } } 
\sum_{\substack{ s, t \in [w]^{\beta}, \\ s \sqcup t \in (\Sigma, \emptyset) } }
\chi_{_{\sq_N}}\! \left( \sigma_{p_1,s} \sigma_{p_2,t}  \right).  
\]

For any set partition $\Sigma$ of $[2\beta]$ let us define 
\[
\widetilde{C}^{\Sigma} :=  \sum_{ \substack{ p_1, p_2 \in \{1, \dots, k\}, \\ p_1 \neq p_2 } } 
\sum_{ \substack{  s, t \in [w]^{\beta}, \\ s \sqcup t \in (\Sigma, \emptyset)  } }
\chi_{_{\sq_N}}\! \left( \sigma_{p_1,s} \right) 
\chi_{_{\sq_N}}\! \left( \sigma_{p_2,t} \right).
\]
Then by \cref{prop:expvalue} applied twice for $k=1$ 
\begin{multline}
\label{eq:C-tilde-VarM}
\sum_{\Sigma} \widetilde{C}^{\Sigma} = \sum_{ \substack{ p_1, p_2 \in \{1, \dots, k\}, \\ p_1 \neq p_2 } } 	 
\left( \sum_{s \in [\water]^\beta} \chi_{_{\sq_N}}\! \left( \sigma_{p_1,s} \right)  \right) 
\left( \sum_{t \in [\water]^\beta} \chi_{_{\sq_N}}\! \left( \sigma_{p_2,t} \right) \right)
= \\
(k^2-k) N^{2\beta} \big( \E_{\Pp} M_{\beta}(w, 1) \big)^2.
\end{multline}
Our aim is to show that \eqref{eq:C-tilde-VarM} is a good approximation for 
the left-hand-side of \eqref{eq:C-class-disconnected-VarM}.

\smallskip

By~\cref{thm:character-value} 
there exists a constant $K_\beta > 0$ 
which depends only on $\beta$ (in particular it does not depend on $N$ or $k$)
such that for any distinct $p_1, p_2$
and a pair of sequences $(s,t)$ 
such that $s \sqcup t \in (\Sigma, \emptyset)$
with $\Sigma$ disconnected
\[
\left| \chi_{_{\sq_N}}\! \left(\sigma_{p_1,s} \sigma_{p_2,t} \right) 
- \chi_{_{\sq_N}}\! \left( \sigma_{p_1,s} \right) 
\chi_{_{\sq_N}}\! \left( \sigma_{p_2,t} \right) \right| 
\leq  K_\beta  \left( N^{-\left|\sigma_{p_1,s}\right| - \left|\sigma_{p_2,t} \right| - 2} \right).
\] 

Let us denote for any set-partition $\Sigma$ of $[2\beta]$ 
\begin{align*}
R^{\Sigma} &:=
 \sum_{ \substack{ p_1, p_2 \in \{1, \dots, k\}, \\ p_1 \neq p_2 } } 
\sum_{ \substack{ s, t \in [w]^{\beta}, \\ s \sqcup t \in (\Sigma, \emptyset) }} 
N^{-\left|\sigma_{p_1,s}\right| - \left|\sigma_{p_2,t} \right| - 2}.
\end{align*}
With the introduced notation the following inequality holds  
\begin{equation}
\label{eq:C-disconnected-value-VarM}
\left|\sum_{\Sigma} \widetilde{C}^{\Sigma} - \sum_{\substack{\Sigma: \\ \Sigma \text{ is disconnected}}} C^{\Sigma} \right|
\leq 
\sum_{\substack{\Sigma: \\ \Sigma \text{ is disconnected}}}
K_\beta  R^{\Sigma} + 
\sum_{\substack{\Sigma: \\ \Sigma \text{ is connected}}}
\left| \widetilde{C}^{\Sigma}\right|.
\end{equation}
We now investigate the asymptotics of the sums on the right-hand-side of~\eqref{eq:C-disconnected-value-VarM}.

\subsubsection{Asymptotics of the right-hand-side of \eqref{eq:C-disconnected-value-VarM}}
\label{subsec:asymptotics-DNS}
Observe that if $\Sigma$ is connected 
and $(s,t)$ is a pair of sequence such that $s \sqcup t \in (\Sigma, \emptyset)$ 
then 
\[
|\Sigma| = 
|\Sigma(s \sqcup t)| \leq 
|\Sigma(s)| + |\Sigma(t)| - 1.
\] 
We follow the lines from \cref{subsec:B-asymptotics} 
to get the upper bound 
\begin{equation}
    \label{eq:C-tilde-connected-VarM}
    \left| \widetilde{C}^{\Sigma} \right|  
    \leq  
    C_\beta^2 \: (k^2-k) N^{2\beta - 2}
\end{equation}
for any connected $\Sigma$ with 
universal constant $C_\beta > 0$ depending only on $\beta$. 

On the other hand, 
by \cref{lem:optimalpi0}\ref{lem:optimalpi0b} and \eqref{eq:EM-rest} 
for any set-partition $\Sigma$ of $[2\beta]$
\begin{equation}
\label{eq:C-Sigma-Rest-VarM}
R^{\Sigma} \leq (k^2-k) 
\left[ 
\sum_{s \in [w]^\beta} N^{-\left|\sigma_{p_1,s}\right| -1} 
\right]^2 
\leq (k^2-k) N^{2\beta - 2}.
\end{equation}

Inserting approximations \eqref{eq:C-tilde-connected-VarM} for connected $\Sigma$
and \eqref{eq:C-Sigma-Rest-VarM} for disconnected $\Sigma$ 
into \eqref{eq:C-disconnected-value-VarM}
and taking into account \eqref{eq:C-tilde-VarM}
proves \eqref{eq:C-class-disconnected-VarM}.

\subsubsection{Asymptotics of the sum $\displaystyle \sum_\Sigma C^{\Sigma}$}
\label{subsec:asymptotics-C}

By~\eqref{eq:C-class-connected-VarM} and \eqref{eq:C-class-disconnected-VarM}
we get 
\begin{equation}
\label{eq:C-value}
\sum_{\substack{\Sigma}} C^{\Sigma} =  
(k^2-k) N^{2\beta} \! \left( \big( \E_{\Pp} M_{\beta}(w, 1) \big)^2 
+ O \! \left(N^{-1} \right) \right).
\end{equation}

\subsubsection{Finishing the~proof of \cref{prop:VarM-asymptotics}}

\begin{proof}[Proof of \cref{prop:VarM-asymptotics}]
Inserting the~equality \eqref{eq:VarM-A}
and approximations \eqref{eq:B-class-estimate-VarM}
(for complementary set-partitions $(\Sigma, \Pi)$ with $\Pi \neq \emptyset$) 
and \eqref{eq:C-value}
into \eqref{eq:prop-VarM}
we conclude that \eqref{prop:VarM-asymptotics-a} 
holds true. 

By~\cref{prop:cosetcharacter} there exists $C_\beta > 0$ such that 
\begin{equation}
    \label{eq:VarM-est}
    \left|
    \chi_{_{\sq_N}}\!  
    \left(\mathcal{J}^2 \Sym_k \right) 
    - \chi_{_{\sq_N}}\!\left( \mathcal{J}^2_* \right)
    \right| 
    \leq C_\beta  \sum_{p_1,p_2=1}^k 
		\sum_{ (s,t) \in P_{p_1, p_2}(\beta) }  
		\frac{k^2}{N^{\| \sigma_{p_1,s} \sigma_{p_2,t} \Sym_k\|+1}} .
\end{equation}
The summands on the right-hand side can be grouped according to 
the types $(\Sigma,\Pi)$ of the sequences $s \sqcup t$. 
By \eqref{eq:VarM-rest} 
the right-hand side of \eqref{eq:VarM-est} 
estimates by $\frac{k^2}{N} \, k^2 N^{2\beta}$
which ends the~proof of \eqref{prop:VarM-asymptotics-b}.

\cref{prop:VarM-asymptotics-c}
is a direct consequence of 
\eqref{prop:VarM-asymptotics-a} and \eqref{prop:VarM-asymptotics-b}.
\end{proof}

\smallskip 

This finishes the proof of \cref{thm:moments}

\subsection{Proof of \cref{prop:longitude-experimental} }
\label{sec:proof-prop:longitude-experimental}

The proof is based on \cite[Section~4.10]{Romik2015a}.
We start with a~general fact.

\begin{lemma}
\label{prop:converge-in-moments}
Let $\varepsilon>0$ and $\mu$ be 
a~compactly supported probability measure on~$\R$. 
Let $x \in \R$ be a~continuity point 
of the~cumulative distribution function $F_\mu$ of~$\mu$. 
Then there exist $\delta >0 $ and an integer $A>0$
with the following property: 
\newline 
if $m$ is a probability measure on $\R$ 
such that its moments (up to order $A$)
are $\delta$-close to the moments of~$\mu$ then 
\[
\left| F_{\mu}(x) - F_m(x) \right| \leq \varepsilon
\]
	where $F_m$ is the~cumulative distribution function of $m$.
\end{lemma}

\begin{proof}
If this would not be the case, there would exist a sequence of probability measures which
converges to $\mu$ in moments, but would not converge to~$\mu$ in the weak
topology of probability measures. This is not possible, since $\mu$ is compactly
supported and therefore uniquely determined by its moments
\cite[Section~3.3.5]{Durrett2010}.
\end{proof}

We now prove \cref{prop:longitude-experimental}.

\begin{proof}[Proof of \cref{prop:longitude-experimental}]
We mimic the~proof of \cite[Theorem~4.1]{Romik2015a}.
Pittel and Romik 
\cite[Theorem~2]{Pittel2007}
found explicitly 
the~limit distribution $\cotmeas_\alpha$ 
which describes the $u$-coordinate of the (scaled) position 
of the~box with the entry~$\lfloor \alpha N^2 \rfloor$ 
in a uniformly random tableau $\tab_N \in \tableaux_{\sq_N}$
as the semicircle distribution~\eqref{eq:density}
(recall \cref{sec:random-position}).
In our setting this result describes 
the $u$-coordinate of the surfer after 
draining $1-\alpha$ fraction of water. 
Since `the amount of remaining water $w$' 
is such that $\frac{w}{N^2} \to \alpha$,
the $\beta$-th moment $\gamma_\beta$ of 
the distribution $\cotmeas_\alpha$ equals 
\[
\gamma_\beta := \int x^{\beta} \dif \cotmeas_{\alpha}(x)
= \lim_{N \to \infty} \E_{\Pp} M_\beta(w,1).
\]

By \cref{thm:moments} we get using Chebyshev's inequality that
for any $\varepsilon>0$ and $\beta \in \N$
\begin{equation*}
\mathbb{P}_{\widetildetableaux_{\sq_N}^{(w+1,w+k)}} \left( 
\left| M_\beta(w,k) - \gamma_\beta \right| > \varepsilon \right)
= O\left(\frac{1}{k} + \frac{k^2}{N} \right).
\end{equation*}

The cumulative distribution function $F_{\cotmeas_\alpha}$ 
of the~semicircle measure $\cotmeas_\alpha$ 
is continuous and therefore any $x \in \R$ 
is its continuity point.
Let $x \in \R$
and let $A$ and $\delta$ 
be the constants given by \cref{prop:converge-in-moments}.
Then by the union bound
\begin{multline}
\label{eq:pointwise-empirical-longitude}
\mathbb{P}_{\widetildetableaux_{\sq_N}^{(w+1,w+k)}}
\bigg\{ 
\left| F_{\cotmeas_{\alpha}}(x) - \empiricallongitude(x) \right| > \varepsilon  
\bigg\} 
\leq \\
\sum_{\beta =1}^A \mathbb{P}_{\widetildetableaux_{\sq_N}^{(w+1,w+k)}} \left( 
\left| M_\beta(w,k) - \gamma_\beta \right| > \delta \right)
= O \left(\frac{1}{k} + \frac{k^2}{N} \right)
\end{multline}
with the constant in the $O$-notation depending only on $\varepsilon$.

\medskip

For the proof of \cref{prop:longitude-experimental} we need to obtain the uniform
version of~\eqref{eq:pointwise-empirical-longitude} in which the event on the
left hand side is taken with supremum over $u \in \R$. This can be easily done by
choosing a finite set $X \subseteq \R$ with the property that its image
$F_{\cotmeas_\alpha}$ is an $\epsilon$-net for the interval $[0,1]$; such a set exists
because $F_{\nu_\alpha}$ is continuous.
The pointwise result \eqref{eq:pointwise-empirical-longitude} implies that the following 
estimate for the supremum over the finite set $X$ holds true:
\begin{equation}
    \label{eq:semi-supremum}
\mathbb{P}_{\widetildetableaux_{\sq_N}^{(w+1,w+k)}}
\bigg\{ 
\multi_N
:
\sup_{x\in X}
\left| F_{\cotmeas_{\alpha}}(x) - \empiricallongitude(x) \right| > \varepsilon  
\bigg\} 
=  O\!\left(\frac{1}{k} + \frac{k^2}{N}\right).
\end{equation}

Let $x_1<\dots<x_\ell$ be the elements of $X$. The assumption about the set $X$ implies that 
\begin{equation}\label{eq:nett} 
F_{\cotmeas_{\alpha}}(x_1)< \epsilon, \qquad    
       F_{\cotmeas_{\alpha}}(x_{i+1}) < F_{\cotmeas_{\alpha}}(x_i) + 2 \epsilon, 
       \qquad F_{\cotmeas_{\alpha}}(x_\ell) > 1- \epsilon.
\end{equation}
The elements of $X$
divide the real line into $\ell+1$ intervals: 
\[ (-\infty,x_1], \ [x_1,x_2], \ \dots ,\ [x_{\ell-1}, x_\ell], \ [x_\ell, \infty). \]
By considering each interval separately, using monotonicity of 
the cumulative distribution function $F_{\cotmeas_\alpha}$
and the monotonicity of~$\empiricallongitude$, as well as \eqref{eq:nett} it follows that
\[ \sup_{x\in \R}
\left| F_{\cotmeas_{\alpha}}(x) - \empiricallongitude(x) \right| <
2\epsilon+
\sup_{x\in X}
\left| F_{\cotmeas_{\alpha}}(x) - \empiricallongitude(x) \right|.\]
In this way \eqref{eq:semi-supremum} completes the proof.
\end{proof}

\section{Proof of \cref{thm:all-the-same}}
\label{sec:connection}
\label{sec:longitude-ends}

The current section is devoted to the proof of \cref{thm:all-the-same}.

\subsection{Overtaking only in one direction}
\label{subsec:overtaking-one-direction}

We start with a precise statement 
of the heuristic ideas from \cref{subsec:overtaking}.

\begin{lemma}
    \label{lem:ghosts}
Fix $k,n \in \N$. 
Let tableaux 
${T \in \tableaux_\mu}$ 
of shape $\mu$ with $n+1$ boxes
and $M \in \widetildetableaux_\nu^{(n+1, n+k)}$ 
of shape $\nu$ with $n+k$ boxes
be such that 
\[
T|_{\leq n} = M|_{\leq n}. 
\]
If $1\leq p\leq k$ is such that 
\[
\pos_{T}(n+1) \preceq \pos_{M}(n+p) \preceq \dots \preceq \pos_{M}(n+k)
\]
then \jdtincomplete preserves the latter relations, that is
\[
\pos_{j(T)}(n+1) \preceq \pos_{j(M)}(n+p) \preceq \dots \preceq \pos_{j(M)}(n+k).
\]
\end{lemma}

\begin{proof}
Clearly, $j(T)\big|_{\leq n}=j(M)\big|_{\leq n}$.
Notice also that 
the~boxes in jeu de taquin paths for tableaux $T$ and $M$
match at least to the boxes $\leq n$. 
As mentioned in \cref{sec:Pieri},
$j(\multi)$ is a $k$-Pieri tableau, so 
\[
\pos_{j(M)}(n+p) \preceq \dots \preceq \pos_{j(M)}(n+k)
\] 
and therefore it remains to prove that 
\begin{equation}
\label{eq:sliding-ghosts}
\pos_{j(T)}(n+1) \preceq \pos_{j(M)}(n+p). 
\end{equation}
We consider the following two cases: 
\begin{itemize}
\item[1.] The box $n+1$ in $T$ slid during \jdtincomplete, 
i.e.,~\mbox{$\pos_{T}(n+1) \neq$} $\pos_{j(T)}(n+1)$.
Consider two subcases:
\begin{itemize}
\item[1a)] The box $n+1$ in~$T$ slid to the left;
in this case it does not matter
how (and if) the~box~$n+p$ in~$M$ slid
and \eqref{eq:sliding-ghosts} holds. 
\item[1b)] The box $n+1$ in~$T$ slid to the bottom. 
Then we use the assumption that $M$ is $k$-Pieri
to show that if
$\pos_{M}(n+p) = \pos_{T}(n+1)$ then 
$n+p$ in~$M$ also slid to the bottom
and \eqref{eq:sliding-ghosts} holds. 
On the other hand if $\pos_{T}(n+1) \prec \pos_{M}(n+p)$
then the box $n+p$ in~$M$ 
must be strictly to the right of $n+1$ in~$T$ 
and it does not matter if it slides or not,
so \eqref{eq:sliding-ghosts} also holds. 
\end{itemize}

\item[2.] The box $n+1$ in $T$ did not slide, 
i.e.,~~$\pos_{T}(n+1) = \pos_{j(T)}(n+1)$.
In this case, the JDT path in $T$ ends on some box $\leq n$ 
which is strictly left-top or strictly right-bottom 
to the $\pos_{T}(n+1)$.
Hence, if $\pos_{T}(n+1) = \pos_{M}(n+p)$
then the~box $n+p$ in $M$ does not slide.
Otherwise, (by the initial relation) 
it must be strictly to the right (and weakly to the bottom) 
of $\pos_{T}(n+1)$
and it does not matter if it slides or not.
All in all, \eqref{eq:sliding-ghosts} holds. 
\qedhere
\end{itemize}
\end{proof}

\subsection{Relative position of the surfer}
\label{sec:propproof}

We recommend the reader to recall the notions in \cref{sec:longitude-begins}
and heuristics for the proof of \cref{thm:all-the-same} in \cref{sec:plan}.

Let $0 < t_1 < t_2 < 1$ and denote
\[
w := \lceil (1-t_1) N^2 \rceil - 1.
\]
Let $k$ be a positive integer such that 
$w + k \leq N^2$. 
By \cref{prop:single-and-multi}
used for $C=1$, $\Delta=t_1$, $a = w$ and $\lambda = \sq_N$ 
there exists a pair of random tableaux 
$\mathbf{T}$, $\mathbf{M}$ defined on the same probability space
with the following properties:
\begin{enumerate}[label=(A\arabic*)]
    \item \label{item:41a}
		$\mathbf{T}$ is a uniformly random element of $\tableaux_{\sq_N}$;
		\item \label{item:41b}
		$\mathbf{M}$ is a random element of $\widetildetableaux_{\sq_N}^{(w+1,w+k)}$
		sampled according to the distribution which fulfills the following
    total variation distance bound
				\[
         \delta\left( \mathbf{M}, \mathbb{P}_{\widetildetableaux_{\sq_N}^{(w+1,w+k)}} \right) \leq d \frac{k^2}{\sqrt{N^2-w}}
				\]
			for some universal constant $d>0$ which depends only on $t_1$; 
    \item \label{item:41c}
		$\mathbf{T}\big|_{\leq w}=\mathbf{M}\big|_{\leq w}$ holds true almost surely.
\end{enumerate}
Since the~dual promotion $\jdtcomplete$ is a bijection, 
by \ref{item:41a} the tableau $\mathbf{T} |_{\leq w+1}$ 
has the same distribution as the initial configuration
of the single surfer~$\tab_N'$ (see \eqref{eq:initial-surfer})
and by~\ref{item:41b} 
the tableau $\mathbf{M}\big|_{\leq w+k}$ 
has approximately the same distribution 
(up to the total variation distance in \ref{item:41b})
as the initial configuration of the multisurfers~$\multi_N'$ 
(see \eqref{eq:initial-multisurfers}). 
Moreover, by \ref{item:41c} the initial configurations of water 
are the same for both stories. 
Therefore we will refer to the~box $\water+1$ in $\mathbf{T}$
as the \emph{surfer} and to the~boxes $\water+1, \dots, \water+k$ in $\mathbf{M}$
as the \emph{multisurfers}.

Recall that the definition \eqref{eq:psi-via-multi} of the random variable
$G_N(u)$ depends implicitly on the choice of the tableau $M_N$. For an~integer
$0\leq q \leq \water$ we denote by $G^q_N(u)$ this random variable obtained by
substituting $M_N$ with $j^q(\mathbf{M})$. We underline here that $\multi_N$ and
$j^q(\mathbf{M})$ need not have the same distribution.
We also define $\dynamicempiriciallongitude{q}$ to be the~fraction of 
the~multisurfers which are to the left of the surfer, more precisely
\begin{equation}
\label{eq:dynamicempiricallongitude}
\dynamicempiriciallongitude{q} := 
G_N^q\left( \frac{1}{N} \content^{j^q(\mathbf{T})}_{\water+1} \right) =
\frac{1}{k} \max\left\{ p\in\{1,\dots,k\}  : 
\content^{j^q(\mathbf{M})}_{\water+p} \leq \content^{j^q(\mathbf{T})}_{\water+1} \right\}.
\end{equation} 
The following proposition gives
a relation between $\dynamicempiriciallongitude{q}$ and 
the theoretical longitude of the surfer
on the common probability space of the surfer and multisurfers 
defined in the beginning of this section.
\begin{proposition}
\label{prop:dynamicalempiricallongitude}
Let $s \geq 0$ and $t > 0$ be such that $s+t < 1$.
Let $w(N) = \lceil (1-t)N^2 \rceil - 1$ for $N \in \N$
and let $k = k(N)$ fulfill \cref{ass:multisurfers}.
Let $q = q(N)$ be a sequence of non-negative integers
such that 
\[
0 \leq q(N) \leq \lfloor (1-t) N^2 \rfloor
\quad \text{and} \quad
\lim_{N \to \infty} \frac{q}{N^2} = s.
\]
Then for any $\varepsilon > 0$
\[
\mathbb{P}  \left( (\mathbf{T}, \mathbf{M}): 
\left| \dynamicempiriciallongitude{q} 
- F_{\cotmeas_{_{1-t-s}}} 
\! \left( \frac{1}{N} u_{\water+1}^{j^q(\mathbf{T})} \right) \right|
> \varepsilon \right) = 
O\! \left(\frac{1}{k} + \frac{k^2}{N} \right)
\]
and the constant in the $O$-notation 
depends only on $s$, $t+s$ and~$\varepsilon$.
\end{proposition}

\begin{proof}
By the discussion below \cref{eq:dynamicempiricallongitude}
\begin{multline*}
\mathbb{P}  \left( (\mathbf{T}, \mathbf{M}): 
\left|  \dynamicempiriciallongitude{q} 
- F_{\cotmeas_{_{1-t-s}}}  \! \left( \frac{1}{N} u_{\water+1}^{j^q(\mathbf{T})} \right) \right| 
> \varepsilon \right) = \\
\mathbb{P}  \left( (\mathbf{T}, \mathbf{M}):
\left| \empiricallongitude^q \! \left(\frac{1}{N} \content^{j^q(\mathbf{T})}_{\water+1} \right)
- F_{\cotmeas_{_{1-t-s}}}  \! \left( \frac{1}{N} u_{\water+1}^{j^q(\mathbf{T})} \right) \right|
> \varepsilon \right).
\end{multline*}
The latter is bounded from above by 
\begin{equation}
    \label{eq:my-little-pony}
    \mathbb{P}   \left( (\mathbf{T}, \mathbf{M}):
\sup_{x \in \R} \left| \empiricallongitude^q (x)
- F_{\cotmeas_{_{1-t-s}}}  ( x ) \right|
> \varepsilon \right).
\end{equation}
The random event in \eqref{eq:my-little-pony} 
is expressed purely in terms of the random tableau~$j^q(\mathbf{M})$ 
and does not involve the random tableau $\mathbf{T}$. For this
reason the probability in \eqref{eq:my-little-pony}, 
due to the condition \ref{item:41b} on page
\pageref{item:41b}
and the bijectivity of the \jdtt $\jdtcomplete$, 
is equal to
\[
\mathbb{P}_{\widetildetableaux_{\sq_N}^{(w-q+1,w-q+k)}} 
\! \left(\multi_N: 
\sup_{x \in \R} \left| \empiricallongitude (x)
- F_{\cotmeas_{_{1-t-s}}}  ( x ) \right|
> \varepsilon \right) 
+ O\!\left(\frac{k^2}{N}\right)
\]
with the constant in the $O$-notation depending only on $t$. 
Since $\lim \frac{w-q}{N^2} = 1-t-s > 0$ we can apply 
\cref{prop:longitude-experimental} 
which completes the proof. 
\end{proof}

\subsection{Proof of the upper bound \eqref{eq:sink-B} in \cref{thm:all-the-same}}
\label{subsec:proof-theoretical-longitude-b}

Let $k = k(N)$ be such that $k(N) \to \infty$ 
and $k(N) = o(\sqrt{N})$, i.e.,
\eqref{eq:minimal-assumption-k} is satisfied. 
We will use the results from \cref{sec:propproof} to prove
the upper bound \eqref{eq:sink-B}. 
Let $q := \lfloor t_2 N^2 \rfloor - \lfloor t_1 N^2 \rfloor$. 
Recall that $w = \lceil (1-t_1) N^2 \rceil - 1$ 
(cf.~\cref{sec:propproof}).  

By \ref{item:41a} from \cref{sec:propproof}
we can translate the probability 
on the left-hand-side of \eqref{eq:sink-B}
into the setting of $\mathbf{T}$ in the following way
\begin{multline*}
\Pp \! \left(\tab \in \tableaux_{\sq_N}: \TheoreticalLongitude(\tim_2) - \TheoreticalLongitude(\tim_1) > \varepsilon \right)  =  \\
\mathbb{P} \! \left((\mathbf{T}, \mathbf{M}): F_{\cotmeas_{_{1-t_2}}} \! \left( \frac{1}{N} u_{\water+1}^{j^q \left(\mathbf{T} \right)} \right)
- F_{\cotmeas_{_{1-t_1}}} \! \left( \frac{1}{N} u_{\water+1}^{\mathbf{T}} \right)  > \varepsilon \right);
\end{multline*} 
note that the event on the right-hand side \emph{does not} involve the tableau $\mathbf{M}$.  
The latter probability can be estimated from above via the union bound  
by the sum of probabilities of the following three events:
\begin{itemize}
\item the fraction of the multisurfers in the final position (i.e., in time $t_2$) 
which are to the left of the surfer is `unusually small', that is 
\[
A := \left\{ (\mathbf{T}, \mathbf{M}): 
F_{\cotmeas_{_{1-t_2}}} \! \left( \frac{1}{N} u_{\water+1}^{j^q \left(\mathbf{T} \right)} \right) 
- \dynamicempiriciallongitude{q} > \frac{\varepsilon}{2} \right\};
\]
\item the number of the multisurfers which are to the left of the surfer 
increases over time, more precisely
\[
B := \left\{ (\mathbf{T}, \mathbf{M}): 
\dynamicempiriciallongitude{q} 
- \dynamicempiriciallongitude{0} >0  \right\};
\]
\item the fraction of the multisurfers in the initial position (i.e, in time $t_1$) 
which are to the left of the surfer is `unusually big', that is 
\[
C := \left\{ (\mathbf{T}, \mathbf{M}): 
\dynamicempiriciallongitude{0} 
- F_{\cotmeas_{_{1-t_1}}} \! \left( \frac{1}{N} u_{\water+1}^{\mathbf{T}} \right)  > \frac{\varepsilon}{2} \right\}.
\]
\end{itemize}

By \cref{prop:dynamicalempiricallongitude} the probabilities of the events $A$ and $C$ 
are of order $O\! \left(\frac{1}{k} + \frac{k^2}{N} \right)$. 
By \cref{lem:ghosts} the event $B = \emptyset$ is impossible. 
The choice of the sequence $k$ as in the beginning of this subsection
implies that the upper bound \eqref{eq:sink-B} holds.

\subsection{Proof of the lower bound \eqref{eq:sink-A} in \cref{thm:all-the-same}}
\label{subsec:proof-theoretical-longitude-a}

For any Young diagram $\lambda$ and tableau $\tab \in \tableaux_\lambda$
we will denote by $\lambda^{\on{tr}}$ and $\tab^{\on{tr}}$, respectively,  
the diagram and the tableau obtained by a \emph{transposition}
of the diagram $\lambda$ and the tableau $\tab$. 

The transposition of
tableaux gives a natural bijection 
between the sets $\tableaux_\lambda$ and
$\tableaux_{\lambda^{\on{tr}}}$ of standard tableaux, respectively, of shape
$\lambda$ and its transpose~$\lambda^{\on{tr}}$.
Moreover, under the transposition the $u$-coordinate of 
the given box in a standard tableau $\tab$
changes its sign, namely for any $n \in \{1, \dots, |T|\}$
\begin{equation}
\label{eq:u-coord-transposed}
u_n^\tab = - u_n^{\tab^{\on{tr}}}.
\end{equation}
In particular, the $u$-coordinate of the surfer $u(X_t)$ changes its sign 
under the transposition, i.e., $u(X_t) = - u(X_t^{\on{tr}})$
where $X_t^{\on{tr}}$ denotes the position of the surfer in the transposed tableau $\tab^{\on{tr}}$. 

Recall that for any $\alpha \in (0,1)$ by $\cotmeas_\alpha$ we denote the limit
measure found by Pittel and Romik on the circle of latitude $\alpha$, see
\cref{sec:random-position}. 
The pushforward of the measure $\cotmeas_\alpha$ under the
involution $\R \ni z \mapsto -z$  is a measure
$\tilde{\nu}_\alpha$ which fulfills the
following equality
\begin{equation}
\label{eq:limit-measure-nu}
\tilde{\nu}_\alpha \! \left((-\infty, u] \right) 
= \cotmeas_\alpha \! \left([-u, \infty) \right)
\quad \text{ for all $u \in \R$}.
\end{equation}
Notice that by \eqref{eq:u-coord-transposed} 
the measure $\tilde{\nu}_\alpha$ 
is the limit measure on the circle of latitude~$\alpha$ 
for the transposed sequence of Young diagrams $(\sq_N^{\on{tr}})$.

Let us denote the position of the surfer in the transposed tableau
by $X_t^*$, cf.~\eqref{eq:evacuation}, and the theoretical longitude
of the surfer in the transposed tableau by $\eta$, i.e.,
\[
\eta(t) := F_{\tilde{\nu}_{1-t}} \! \left( u (X_t^*) \right)
\quad 
\text{for $t \in [0,1]$}.  
\]
Observe that by \eqref{eq:limit-measure-nu}
and the continuity of $F_{\cotmeas_{1-t}}$
\begin{equation}
\label{eq:limit-measure-transposition-correspondence}
\eta(t) = 1 - F_{\cotmeas_{1-t}}\big( {- u(X_t^*)} \big).
\end{equation}

\cref{eq:sink-B} applied to 
the transposed tableaux gives
\[
\lim_{N\to\infty} \Pp \bigg\{ \tab \in \tableaux_{\sq_N^{\on{tr}}}:  
    \eta(\tim_2) - \eta(\tim_1)  > \varepsilon \bigg\} =0. 
\]
On the other hand
by \eqref{eq:limit-measure-transposition-correspondence}
and then \eqref{eq:u-coord-transposed}
we get 
\begin{multline*}
\eta(\tim_2) - \eta(\tim_1) = 
\big[ 1 - F_{\cotmeas_{1-\tim_2}}\big( {-u(X_{\tim_2}^* )} \big) \big] 
- \big[ 1 - F_{\cotmeas_{1-\tim_1}}\big( {-u(X_{\tim_1}^* )} \big) \big]  = \\[2ex]
F_{\cotmeas_{1-\tim_1}}\Big( u \big( (X_{\tim_1}^*)^{\on{tr}} \big) \Big) 
- F_{\cotmeas_{1-\tim_2}}\Big( u \big( (X_{\tim_2}^*)^{\on{tr}} \big) \Big).
\end{multline*}
Since $(X_{\tim}^*)^{\on{tr}}$ reflects the position of the surfer 
in the original (not-transposed) tableau we have
\[
\eta(\tim_2) - \eta(\tim_1)
= \TheoreticalLongitude(\tim_1) - \TheoreticalLongitude(\tim_2).
\]
Since we consider the uniform distribution 
on the set of tableaux,
this ends the proof of the lower bound \eqref{eq:sink-A} 
and completes the proof of \cref{thm:all-the-same}.

	\section{Proof of \cref{thm:evacuation}}
	\label{sec:proof-of-evacuation-for-squares}

		\subsection{Plan for the proof of \cref{thm:evacuation}}
		\label{subsec:plan-proof-thm-evacuation}

		We will make the following steps in the~proof of \cref{thm:evacuation}:
		\begin{enumerate}[label=(S\arabic*)]
		\item \label{step:1}
        Pick a~candidate for the~random variable $\Longitude_N\colon \tableaux_{\sq_N} \to [0,1]$.
		\item \label{step:2}
        Prove a~\emph{pointwise version of \cref{thm:evacuation}:} 
		with the~help of \cref{lem:true-position} and \cref{thm:all-the-same}
		we will show that the~chosen candidate 
		gives a~good approximation of surfer's position 
		for an arbitrary $t \in (0,1)$,
		i.e.,
		\begin{equation}
		\label{eq:pointwise-thm-evacuation}
		\bigforall_{t \in (0,1)}\ 
			\lim_{N \to \infty} 
			\Pp \left(\tab_N \in \tableaux_{\sq_N}: 
			\left| \mypoint_\tim(\tab_N) -
			\point_{1-\tim, \Longitude_N(\tab_N)}  \right|
			> \varepsilon \right) 
			= 0.
		\end{equation}
        The proof will be given in \cref{subsec:pointwise-thm-version-proof}.
    
		\item \label{step:3} Prove \emph{the full (i.e., the original) version of
    \cref{thm:evacuation}}, i.e.,
		\[
			\lim_{N \to \infty} 
			\Pp \left(\tab_N \in \tableaux_{\sq_N}: 
			\sup_{\tim \in [0,1]} 
			\left| \mypoint_\tim(\tab_N) -
			\point_{1-\tim, \Longitude_N(\tab_N)}  \right|
             > \varepsilon \right) 
			= 0.
		\] 
		We will start with a finite $\varepsilon$-net of the~family of
		level curves $\{h_\alpha : \alpha \in [0,1]\}$ parameterized by \mbox{$0 = \alpha_0 <
		\alpha_1 < \dots < \alpha_p < \alpha_{p+1} = 1$}. By the~previous point,
		\eqref{eq:pointwise-thm-evacuation} holds uniformly for 
		$t \in \{\alpha_0, \dots, \alpha_{p+1} \}$ in a finite set. Then for the intermediate moments of time
		\mbox{$\alpha_i < t < \alpha_{i+1}$} we will use the~monotonicity of the \jdtp and in this way
		justify that \eqref{eq:pointwise-thm-evacuation} holds \emph{uniformly} over $t \in
		(0,1)$. 
		This proof will be given in \cref{sec:proof-step3}. 
		\end{enumerate}
		
		In the very end, in \cref{subsec:limit-distribution}, we will show that 
		the probability distribution of the~random variable $\Longitude_N$ converges 
		to the uniform distribution on the unit interval $[0,1]$.

		\subsection{Auxiliary notation}
		\label{subsec:auxiliary-notation}
		
		In the~proof we will switch 
		between the~position of the~surfer 
		in the $XY$ and the $UV$ coordinate systems as well as in the geographic coordinates.
		For any $t \in (0,1)$ let us introduce 
		the notation for the true and the theoretical 
		$u-$ and $v-$coordinate; 
		namely we define the following 
		functions $\tableaux_{\sq_N} \to \R$ by 
		\[
		\begin{array}{l}
		u_t := u\! \left(X_{t} \right), 
		\\[5pt]
		v_t := v\! \left( X_{t} \right)
		\end{array}
		\qquad	\text{and} \qquad
		\begin{array}{l}
		u^{\theor}_t:= \left(F_{\cotmeas_{_{1-t}}}\right)^{-1} \! \left(\TheoreticalLongitude(t) \right), 
		\\[5pt]
		v^{\theor}_t := h_{1-t} \! \left( u^{\theor}_{t}  \right)
		\end{array}
		\]
		(recall that $\TheoreticalLongitude(t) = F_{\cotmeas_{_{1-t}}}(u_t)$,
		cf.~\eqref{eq:estimate-longitude},
		and see \cref{sec:levelcurves}
		for the definition of $h_t$).
		Denote additionally for any $t \in [0,1]$ 
		\[
		\left( x_t, y_t \right) := X_{t} 
		\text{\quad and \quad}
		\left( x^{\theor}_t, y^{\theor}_t \right) 
		:= \left( \frac{v^{\theor}_t 
		- u^{\theor}_t}{2}, \frac{u^{\theor}_t + v^{\theor}_t}{2}  \right).
		\]

\subsection{The~surfer's position can be asymptotically recovered from the theoretical longitude}
		
		We start with the result which shows that, in principle, it is possible to
recover the the~true position of the~surfer $X_t$ in the time $t\in (0,1)$ from
its theoretical longitude $\TheoreticalLongitude(\tim)$.
		
		\begin{lemma}
		\label{lem:true-position}
		Let $0 < t < 1$. 
		For any $\varepsilon>0$
		\begin{align}
		\label{eq:u-th-limit}
		\lim_{N \to \infty} 
		\Pp \! \left(\tab_N \in \tableaux_{\sq_N}: \ 
		u_t(\tab_N) \neq u^{\theor}_t(\tab_N)
		\right)
		&= 0
		\intertext{and} 
		\label{eq:v-th-limit}
		\lim_{N \to \infty} 
		\Pp \! \left(\tab_N \in \tableaux_{\sq_N}: \ 
		\left| v_t(\tab_N)
		- v^{\theor}_t(\tab_N)  \right| 
		> \varepsilon \right) 
		&= 0.
		\end{align}
		\end{lemma}

		\begin{proof}
		Let $0 < t < 1$ and $\varepsilon > 0$.
		We start with the proof of \eqref{eq:u-th-limit}.
		Recall  that the~position of the~surfer~$X_t$
		corresponds to the~position of the~box with the~number \mbox{$\lceil (1-t)N^2 \rceil$} 
		in the tableau $(\jdtcomplete)^{\lfloor t N^2 \rfloor} (\tab_N)$, cf.~\eqref{eq:evacuation},
		and so by \cite[Theorem~2]{Pittel2007} the random variable $u_t$ converges in distribution 
		to the measure $\cotmeas_{1-t}$ which has no atoms
		and has a compact connected support 
		\[
		\on{supp}(\cotmeas_{1-t}) = \left[ -2\sqrt{t(1-t)},\; 2\sqrt{t(1-t)} \right].
		\]
		
        \medskip
        
		Observe that since $\cotmeas_{1-t}$ has no atoms, 
		whenever $u_t \in \on{supp}(\cotmeas_{1-t})$ 
		then by the definition
		\begin{equation}
		\label{eq:ut-in-supp}
		u^{\theor}_t 
		=  F_{\cotmeas_{_{1-t}}}^{-1} \! \left(\TheoreticalLongitude(t) \right) 
		= u_t.
		\end{equation}
		On the other hand 
		\begin{equation}
		\label{eq:ut-not-in-supp}
		\Pp \! \left(\tab_N \in \tableaux_{\sq_N}: \  
		u_t(\tab_N) \notin \on{supp}(\cotmeas_{1-t}) \right) \xrightarrow{N \to \infty} 0.
		\end{equation}
		Indeed, for any $\varepsilon > 0$ 
		consider the function $f_\varepsilon : \R \to \R$
		given by 
		\[
		f_\varepsilon(x) :=
		1 - \min\left\{ 1, \frac{1}{\varepsilon} \on{dist}\left(x, \R \setminus \on{supp}(\cotmeas_{1-t}) \right) \right\},
		\quad
		x \in \R,
		\]
		where $\on{dist}(x, A) := \inf\{ d(x,y): y \in A \}$ 
		is the Hausdorff distance
		of the point $x$ from a set $A$.
		Clearly $f_\varepsilon$ is continuous and bounded.
		
		The distribution of the random variable~$u_t$
		is a pushforward of the measure~$\Pp$
		under the mapping $\tab_N \mapsto u_t(\tab_N)$.
		Therefore by the convergence in distribution of 
		the random variable $u_t$ to 
		the distribution given by the measure $\cotmeas_{1-t}$
		we get for any \mbox{$\varepsilon > 0$}
		\begin{multline}
		\label{eq:indicator-estimation}
		\Pp \! \left(\tab_N \in \tableaux_{\sq_N}: \  
		u_t(\tab_N) \notin \on{supp}(\cotmeas_{1-t}) \right) 
		\leq \\
		\int_{\tableaux_{\sq_N}} f_\varepsilon(u_t) \diff \Pp 
		\xrightarrow{N \to \infty} \int_\R f_\varepsilon \diff \cotmeas_{1-t}.
		\end{multline}
		Clearly, $f_\varepsilon$ converges pointwise as $\varepsilon \to 0$
		to the indicator function $\mathds{1}_{\overline{\R \setminus \on{supp}(\cotmeas_{1-t})}}$,
		therefore by the Lebesgue dominated convergence theorem 
		\[
		\lim_{\varepsilon \to 0} \int_\R f_\varepsilon \diff \cotmeas_{1-t}
		= \cotmeas_{1-t}\left( \overline{\R \setminus \on{supp}(\cotmeas_{1-t})} \right)
		= \cotmeas_{1-t}\Big( \partial \big(\on{supp} (\cotmeas_{1-t}) \big) \Big)
		= 0.
		\]
		This together with \eqref{eq:indicator-estimation} 
		implies \eqref{eq:ut-not-in-supp}. 
        
		A conjunction of \eqref{eq:ut-in-supp} 
		and \eqref{eq:ut-not-in-supp} proves \eqref{eq:u-th-limit}.

		\medskip
        
		\cref{eq:v-th-limit} follows from the result of Biane \cite[Theorem 1.5.1]{Biane1998}.
		We will shortly describe it here, but 
		the more developed discussion and precise statements 
		formulated using our notation
		are placed in \cref{sec:generalization}.

		The boundary of a Young diagram~$\lambda$ seen in the $(u,v)$-coordinate system 
		can be viewed as a non-negative $1$-Lipschitz function $\omega_\lambda$,
		see \cref{fig:continuous-diagram}.
		The function $\omega_\lambda$ is initially defined on the interval $I$ 
		given by the range of the $u$-coordinates of~$\lambda$,
		but it can be extended to a function defined on the real line $\R$
		by gluing $\omega_\lambda|_{I}$ with the modulus function $x \mapsto |x|$.
		This extended function $\omega_{\lambda}\colon \R\to\R_+$ is called \emph{the profile of $\lambda$},
		cf.~\cref{subsec:continuous-diagrams}.
		
		\begin{figure}[tbp]
\centering
\begin{tikzpicture}[scale=0.8]

\begin{scope}[scale=1]
 
        \begin{scope}[draw=gray,rotate=45,scale=sqrt(2)]
          \fill[fill=blue!10] (4,0) -- (4,1) -- (3,1) -- (3,2) -- (1,2) -- (1,3) -- (0,3) -- (0,0) -- cycle ;
        \end{scope}

       \begin{scope}
          \clip (-4.5,0) rectangle (5.5,5.5);
          \draw[thin, dotted] (-6,0) grid (6,6);
          \begin{scope}[rotate=45,draw=gray,scale=sqrt(2)]
              \clip (0,0) rectangle (4.5,5.5);
              \draw[thin, dotted] (0,0) grid (6,6);
          \end{scope}
      \end{scope}

      \draw[->,thick] (-6,0) -- (6,0) node[anchor=west]{$u$};
      \foreach \z in {-4, -3, -2, -1, 1, 2, 3, 4, 5}
            { \draw (\z, -2pt) node[anchor=north] {\tiny{$\z$}} -- (\z, 2pt); }

      \draw[->,thick] (0,-0.4) -- (0,6) node[anchor=south]{$v$};
      \foreach \t in {1, 2, 3, 4, 5}
            { \draw (-2pt,\t) node[anchor=east] {\tiny{$\t$}} -- (2pt,\t); }

  \begin{scope}[draw=gray,rotate=45,scale=sqrt(2)]

          \draw[->,thick] (0,0) -- (6,0) node[anchor=west,rotate=45] {\textcolor{gray}{{$x$}}};
          \foreach \x in {1, 2, 3, 4, 5}
              { \draw (\x, -2pt) node[anchor=north,rotate=45] {\textcolor{gray}{\tiny{$\x$}}} -- (\x, 2pt); }

          \draw[->,thick] (0,0) -- (0,5) node[anchor=south,rotate=45] {\textcolor{gray}{{$y$}}};
          \foreach \y in {1, 2, 3, 4}
              { \draw (-2pt,\y) node[anchor=east,rotate=45] {\textcolor{gray}{\tiny{$\y$}}} -- (2pt,\y); }

          \draw[ultra thick,draw=blue] (5.5,0) -- (4,0) -- (4,1) -- (3,1) -- (3,2) -- (1,2) -- (1,3) -- (0,3) -- (0,4.5) ;

  \end{scope}

\end{scope}
\end{tikzpicture}

\caption{A Young diagram $\lambda=(4,3,1)$ shown in the Russian convention. 
The blue solid line represents its profile $\omega_\lambda$. 
The $(u,v)$-coordinate system corresponding to the Russian convention 
and the $XY$- coordinate system corresponding to the French convention are shown.}
\label{fig:continuous-diagram}
\end{figure}

		The restriction $\tab|_{\leq \lceil (1-t) N^2 \rceil}$ of the random tableau $\tab$
		has a (random) shape 
		\[ \lambda_{1-t} := \sh \tab|_{\leq \lceil (1-t) N^2 \rceil}\] 
		whose (random) profile $\omega_{\lambda_{1-t}}$ 
		is such that the following equality 
		in terms of the $u$- and $v$-coordinates of 
		the surfer (in time $t$) holds true:
		\[
		v_t (\tab_N) = \omega_{\lambda_{1-t}}\! \left( u_t(\tab_N) \right).
		\]
			On the other hand if $u_t(\tab_N) = u^{\theor}_t(\tab_N)$ then
		\[
		v_t^{\on{th}}(\tab_N) 
		= h_{1-t}(u_t^{\on{th}}(\tab_N)) 
		= h_{1-t} (u_t(\tab_N)).
		\]
			
			We proved in \eqref{eq:u-th-limit} that  
		the set of tableaux $\tab_N \in \tableaux_{\sq_N}$
		for which $u_t(\tab_N) = u^{\theor}_t(\tab_N)$ has asymptotically full probability.
		Therefore with asymptotically full probability 
		the following inequality holds
		\begin{equation}
    \label{eq:supremum-ok}
             \left| v_t (\tab_N) - v_t^{\on{th}}(\tab_N) \right| \leq \sup_{x \in \R} 
        \left|  \omega_{\lambda_{1-t}}(x) - h_{1-t}(x) \right|. 
		\end{equation}
		
		By the aforementioned result \cite[Theorem 1.5.1]{Biane1998} 
		(see \cref{lem:limit-shape} for the precise general statement)
	the right-hand-side of \eqref{eq:supremum-ok}
	converges in probability to $0$.
		This completes the proof of \eqref{eq:v-th-limit}
		and \cref{lem:true-position}.

		\end{proof}

		\subsection{A candidate for the random variable $\Longitude_N(\tab_N)$, 
		development of \ref{step:1}}
		\label{subsec:random-variable}
		
		Pick any $t_0 \in (0,1)$ and define 
		for $\tab_N \in \tableaux_{\sq_N}$
		\[
		\Longitude_N(\tab_N) := \TheoreticalLongitude(t_0)(\tab_N).
		\]
    We will show that the random variable $\Longitude_N$ 
		has the desired properties from 
		\cref{thm:evacuation}.
		
		For any $t\in[0,1]$ we define the \emph{approximated $u$- and $v$-coordinate} of the surfer as  
		\[
		\begin{array}{l}
		u^{\Longitude}_t:= \left(F_{\cotmeas_{_{1-t}}}\right)^{-1} \! \left(\Longitude_N \right), 
		 \\[5pt]
		v^{\Longitude}_t := h_{1-t} \! \left( u^{\Longitude}_{t}  \right).
		\end{array}
		\]
		These definitions were chosen in such a way that that the \emph{approximated position
of the surfer} in the $XY$ coordinate system
		\[
		\left( x^{\Longitude}_t, y^{\Longitude}_t \right) 
		:= \left( \frac{v^{\Longitude}_t 
		- u^{\Longitude}_t}{2}, \frac{u^{\Longitude}_t + v^{\Longitude}_t}{2}  \right)
		= \point_{1-\tim, \Longitude_N(\tab_N)}
		\]
is the point which appears in the statement of \cref{thm:evacuation} and
Eq.~\eqref{eq:pointwise-thm-evacuation}.

		With this notation the expression in the modulus in the event 
		in \eqref{eq:pointwise-thm-evacuation} takes the form
		\[
		\left| \mypoint_\tim(\tab_N) -
			\point_{1-\tim, \Longitude_N(\tab_N)}  \right|
		= \frac{1}{\sqrt{2}} \left| \left(u_{t}, v_{t} \right) 
		- \left(u^{\Longitude}_{t}, v^{\Longitude}_{t} \right) \right|.
		\]

		\subsection{The proof of \ref{step:2} -- the pointwise version of \cref{thm:evacuation}}
		\label{subsec:pointwise-thm-version-proof}
		
		Let \mbox{$\varepsilon > 0$} and $t \in (0,1)$.
		By the triangle inequality, 
		\begin{align}
		\left| u_t(\tab_N)  - u_t^{\Longitude}(\tab_N) \right|
		&\leq 
		\left| u_t(\tab_N)  - u_t^{\theor}(\tab_N) \right|
		+ 
		\left| u_t^{\theor}(\tab_N)  - u_t^{\Longitude}(\tab_N) \right|
		\label{eq:first-inequality}
		\intertext{and}
		\left| v_t(\tab_N)  - v_t^{\theor}(\tab_N) \right|
		&\leq 	\left| v_t(\tab_N) - v_t^{\theor}(\tab_N) \right| 
		+ 		\left| v_t^{\theor}(\tab_N) - v_t^{\Longitude}(\tab_N)\right|.
		\label{eq:second-inequality}
		\end{align}
		
		By \cref{lem:true-position},
		in each of the above two inequalities
		the first summand on the right-hand side 
		converges in probability to $0$, that is,  
		\[
		\left| u_t(\tab_N)  - u_t^{\theor}(\tab_N) \right| \xrightarrow{P} 0
		\quad \text{and} \quad 
		\left| v_t(\tab_N)  - v_t^{\theor}(\tab_N) \right| \xrightarrow{P} 0.
		\]
		
The second summands on the right hand side of 
\eqref{eq:first-inequality} and \eqref{eq:second-inequality} 
are the distances between the values of uniformly continuous functions, respectively, 
$\psi \mapsto u_t^\psi$ and the composition $\psi \mapsto h_{1-t}(u_t^\psi)$
evaluated at the arguments 
\[
\TheoreticalLongitude(t)
\quad \textrm{ and } \quad 
\TheoreticalLongitude(t_0) = \Longitude_N.
\]
By \cref{thm:all-the-same}
with $t_1 = \min(t,t_0)$ and $t_2 = \max(t,t_0)$ 
the distance between these two arguments 
converges in probability to $0$, i.e., 
		\[ 
		\TheoreticalLongitude(t)
		- \Longitude_N  \xrightarrow{P} 0.
		\]	
	Therefore we get 
		\[
		\left| u_t^{\theor}(\tab_N)  - u_t^{\Longitude}(\tab_N) \right| \xrightarrow{P} 0 
		\quad \text{and} \quad 
		\left| v_t^{\theor}(\tab_N) - v_t^{\Longitude}(\tab_N)\right|  \xrightarrow{P} 0.
		\]
		
		As the result
		\[
		\left| u_t(\tab_N)  - u_t^{\Longitude}(\tab_N) \right| \xrightarrow{P} 0 
		\quad \text{and} \quad 
		\left| v_t(\tab_N) - v_t^{\Longitude}(\tab_N)\right|  \xrightarrow{P} 0
		\]
		which completes the proof \eqref{eq:pointwise-thm-evacuation} which is the
pointwise version of \cref{thm:evacuation}.

		\subsection{The proof of \ref{step:3} -- the full version of \cref{thm:evacuation}}
		\label{sec:proof-step3}

		\subsubsection{Uniform continuity of the geographic coordinate system}
		
		We start with showing that the geographic coordinate system
		on the square (recall \cref{sec:geographic})
		is uniformly continuous. 
				
		\begin{lemma}
		\label{lem:uniform-continuity}
		The function
		\begin{align*}
		[0,1] \times [0,1] \ni 
		(\alpha, \psi) \mapsto P_{\alpha, \psi} 
		= \left( x_\alpha^\psi, y_\alpha^\psi \right) 
		\end{align*}
		is uniformly continuous.
		\end{lemma}
		
		\begin{proof}
		Since the mapping 
		$\left( x_\alpha^\psi, y_\alpha^\psi \right)  \mapsto 
		\frac{1}{\sqrt{2}} \left( u_\alpha^\psi, v_\alpha^\psi \right)$
		is an isometry (as a rotation in $\R^2$)
		it is enough to show that each coordinate of the function  
		\[
		[0,1] \times [0,1] \ni 
		(\alpha, \psi) \mapsto 
		\left( u_\alpha^\psi, v_\alpha^\psi \right)
		= \left( u_\alpha^\psi, h_{1-\alpha}\! \left(u_\alpha^\psi \right) \right)
		\]
		is uniformly continuous.
		
		Recall that the limit distribution $\cotmeas_\alpha$
		has density \eqref{eq:density}.
		It is easy to show that for any $\alpha \in (0,1)$ 
		the cumulative distribution function $F_{\cotmeas_\alpha}$ of $\cotmeas_\alpha$ 
		fulfills the equation
		\[
		F_{\cotmeas_\alpha}(u) = F_{\operatorname{SC}}\! \left(\frac{u}{2\sqrt{\alpha(1-\alpha)}} \right) 
		\quad \text{for all $|u| \leq 2\sqrt{\alpha(1-\alpha)}$}
		\]
		where $F_{\operatorname{SC}}$ denotes 
		the cumulative distribution function of 
		the standard semicircle distribution 
		with the density \eqref{eq:SC}.
		This implies that for any $\psi \in [0,1]$ and $\alpha \in (0,1)$ we get
		\[
		\psi = F_{\cotmeas_{1-\alpha}}\! \left( u^{\psi}_{\alpha} \right) 
		= F_{\operatorname{SC}}\! \left( \frac{u^{\psi}_{\alpha}}{2\sqrt{\alpha(1-\alpha)}} \right),
		\]
		which after applying $F_{\operatorname{SC}}^{-1}$ gives 
		\begin{equation}
		\label{eq:square-alpha-quantile}
		u^{\psi}_{\alpha} = 2\sqrt{\alpha(1-\alpha)} \cdot F_{\operatorname{SC}}^{-1}(\psi).
		\end{equation}
		Moreover, by the definition $P_{0, \psi} = (0,0) \in \R^2$ and $P_{1, \psi} = (1,1) \in \R^2$
		(cf.~\cref{sec:geographic}),
		hence \eqref{eq:square-alpha-quantile} holds for all $\psi \in [0,1]$ and $\alpha \in [0,1]$.
    
		Therefore the mapping $(\alpha, \psi) \mapsto u_\alpha^\psi$ 
		is uniformly continuous
		since $\psi \mapsto F_{\operatorname{SC}}^{-1}(\psi)$ is uniformly continuous.
		Indeed, $F_{\operatorname{SC}}^{-1}$ is the inverse function of $F_{\operatorname{SC}}$ 
		which is injective on $\left[-2\sqrt{\alpha(1-\alpha)}, 2\sqrt{\alpha(1-\alpha)}\right]$,
		continuous and has compact domain and range.  
		
		\medskip
		
		The function $(\alpha, u) \mapsto h_{\alpha}(u)$
		is uniformly continuous
		on the domain 
		\[
		\Diamond := \left\{ (\alpha, u):
		\alpha \in [0,1]
		\text{ and }
		|u| \leq 2\sqrt{\alpha(1-\alpha)}
		\right\}
		\]
		(as a continuous mapping on the compact set $\Diamond$).
		Hence the function
		$(\alpha, \psi) \mapsto v_\alpha^\psi$
		is uniformly continuous as the composition of two uniformly continuous functions:
		\[
		\Diamond \ni (\alpha, u) \mapsto h_{\alpha}(u)
		\quad \text{and} \quad
		(\alpha, \psi) \mapsto (1-\alpha, u_\alpha^\psi) \in \Diamond. \qedhere
    \] 
		\end{proof}

		\subsubsection{The proof of \ref{step:3} -- the full version of \cref{thm:evacuation}}
		\label{subsec:polishing-step3}
		
		Let $\varepsilon > 0$.
		By \cref{lem:uniform-continuity} the function
		\[
		[0,1] \times [0,1] \ni 
		(\alpha, \psi) \mapsto P_{\alpha, \psi} 
		\]
		is uniformly continuous, so there exists $\delta>0$ such that
		\begin{equation}
		\label{eq:unif-cont-delta-net}
		\bigforall_{s,t \in [0,1]}\ 
		|s-t| < \delta \implies 
		\bigforall_{\psi \in [0,1]}\ 
		\left| 
		P_{s, \psi} 
		- P_{t, \psi} 
		\right| 
		< \varepsilon.
		\end{equation}
		
		Let us take a finite $\delta$-net $0 = \alpha_1 < \dots < \alpha_n = 1$ of the
interval $[0,1]$. By the pointwise version \ref{step:2} of \cref{thm:evacuation},
which we proved in \cref{subsec:pointwise-thm-version-proof},
		\begin{equation}
		\label{eq:poinwise-result-use}
		\left| \mypoint_{\alpha_i}(\tab_N) -
			\point_{1-\alpha_i, \Longitude_N(\tab_N)} \right|
			 \xrightarrow{P} 0 
			\quad 
			\text{for \; $i \in \{1, \dots, n \}$}
		\end{equation}
		(the latter holds for $i = 1$ and $i=n$ 
		by the definition of $P_{\alpha, \psi}$, cf.~\cref{sec:geographic}).
		Therefore	
		there exists a subset $\tableaux^*_N\subseteq \tableaux_{\sq_N}$ 
		of asymptotically full measure
        which consists of tableaux $\tab_N$
		with the property that 
		for each $i \in \{1, \dots, n \}$
		\begin{equation}
		\label{eq:surfer-trajectory}
		\left| x_{\alpha_i}(\tab_N) - x^{\Longitude}_{\alpha_i}(\tab_N) \right| < \varepsilon
		\text{\quad and \quad}
		\left|y_{\alpha_i}(\tab_N) - y^{\Longitude}_{\alpha_i}(\tab_N) \right| < \varepsilon.
		\end{equation}
		
		By the monotonicity of the \jdtp
		for any $i \in \{1, \dots, n-1 \}$
		and $t \in [\alpha_i, \alpha_{i+1}]$ 
		\begin{equation}
		\label{eq:jdt-monotonicity}
		x_{\alpha_{i+1}} \leq x_t \leq x_{\alpha_i} 
		\text{\quad and \quad}
		y_{\alpha_{i+1}} \leq y_t \leq y_{\alpha_i}.
		\end{equation}
		
		By \eqref{eq:surfer-trajectory} and \eqref{eq:jdt-monotonicity}  
		for any $\tab_N \in \tableaux^*_N$
		and any $i \in \{1, \dots, n-1 \}$
		and $t \in [\alpha_i, \alpha_{i+1}]$
		the following system of inequalities is satisfied:
        \begin{empheq}[left=\empheqlbrace]{equation}
            \notag
            \begin{aligned}
	- \varepsilon + \left( x^{\Longitude}_{\alpha_{i+1}} - x^{\Longitude}_t \right) 
		 & \leq x_t - x^{\Longitude}_t 
		\leq \left(x^{\Longitude}_{\alpha_i} - x^{\Longitude}_t \right) + \varepsilon; 
		\\
		- \varepsilon + \left( y^{\Longitude}_{\alpha_{i+1}} - y^{\Longitude}_t \right) 
		& \leq y_t - y^{\Longitude}_t 
		\leq \left(y^{\Longitude}_{\alpha_i} - y^{\Longitude}_t \right) + \varepsilon.
            \end{aligned}
        \end{empheq}
		
		Since for any $t \in [0,1]$
		\[
		X_t = (x_t, y_t) 
		\quad \text{and} \quad 
		\point_{1-t, \Longitude_N(\tab_N)} = \left(x^{\Longitude}_t , y^{\Longitude}_t \right)
		\]
		and by \eqref{eq:unif-cont-delta-net}
		for any $i \in \{ 1, \dots, n-1 \}$
		and $t \in [\alpha_i, \alpha_{i+1}]$
		\[
		\max\left\{ 
		\left| x^{\Longitude}_{\alpha_i} - x^{\Longitude}_t  \right|, \
		\left| x^{\Longitude}_{\alpha_{i+1}} - x^{\Longitude}_t  \right|, \
		\left| y^{\Longitude}_{\alpha_i} - y^{\Longitude}_t  \right|, \
		\left| y^{\Longitude}_{\alpha_{i+1}} - y^{\Longitude}_t  \right|
				\right\}
				< \varepsilon
		\]
		we infer that
		for $\tab_N \in \tableaux^*_N$
		\[
		\sup_{t \in [0,1]} \left| \mypoint_{t}(\tab_N) 
		- \point_{1-t, \Longitude_N(\tab_N)}  \right| < 2 \varepsilon.
		\]
		This completes the proof of \ref{step:3} since 
		$\tableaux^*_N$ has asymptotically full probability.

\subsection{Limit distribution of the random variable~$\Longitude_N$}
\label{subsec:limit-distribution}
		       
		We will show the second component of \cref{thm:evacuation}, namely 
		that the random variable~$\Longitude_N$ converges in distribution
		to the~uniform distribution on the~unit interval~$[0,1]$.
		
		Let $G_N$ denote the cumulative distribution function of the random variable 
		$\Longitude_N : \tableaux_{\sq_N} \to [0,1]$.
		For any $z \in [0,1]$ we have (recall \eqref{eq:estimate-longitude})
		\begin{multline}
            \label{eq:cdf}
		G_N(z) 
		= \Pp \! \left( \tab_N: \ 
		F_{\cotmeas_{_{1-t_0}}} \! \left( u_{t_0}(\tab_N) \right) \leq z 
		\right) = \\
        \Pp \! \left(\tab_N \in \tableaux_{\sq_N}: \ 
		u_{t_0}(\tab_N) \leq \left(F_{\cotmeas_{_{1-t_0}}} \right)^{-1}(z) \right).
	\end{multline}

		By \cite[Theorem~2]{Pittel2007}, the distribution of the random variable
$u_{t_0}$ converges weakly (as $N\to\infty$) to the measure
$\cotmeas_{_{1-t_0}}$ which has no atoms, so
the right-hand side of \eqref{eq:cdf} converges to 
\[ F_{\cotmeas_{_{1-t_0}}}\left(   F_{\cotmeas_{_{1-t_0}}}^{-1}(z) \right)=z\]
which is the cumulative distribution function of the uniform measure $U(0,1)$. 
		This completes the proof of 
		the second component of \cref{thm:evacuation},
		and hence the proof of \cref{thm:evacuation}.

		\section{The correspondence between evacuation and \jdtps}
\label{sec:equivalence-of-problems}

The results which we consider in this section hold for general, 
not necessarily square tableaux. 
For a tableau $\tab\in\tableaux_{\lambda}$ with $n=|\lambda|$ boxes we 
denote by 
\[    \revevac(T)=
\Big( 
\pos_n\! \big( j^{n-1}(\tab) \big),\ 
\dots,\
\pos_n\! \big( j^{1}(\tab) \big),\
\pos_n(\tab)   
\Big)
\]
the evacuation path \eqref{eq:evacuationtrajectory}
written in the reverse order.

The following result shows an intimate relationship between the \jdtps and the
evacuation paths and, in particular, implies equivalence of
\cref{thm:evacuation,thm:jdt} 
(see \cref{sec:equivalence-of-main-problems} for
the proof).

\begin{proposition}
    \label{thm:equivalence} 
    
    Let $\lambda$ be a fixed Young diagram and let $\tab \in \tableaux_{\lambda}$ 
		be a random standard Young tableau sampled according
    to the uniform measure on $\tableaux_\lambda$. Then the~probability
    distributions of the lazy \jdtp $\pat(\tab)$ and the evacuation
    path $\revevac(\tab)$ coincide.
\end{proposition}
\begin{proof}
    We will construct a certain bijection $\varepsilon^*\colon \tableaux_{\lambda}
    \to\tableaux_{\lambda}$ on the~set of tableaux of shape $\lambda$. 
		Clearly, the random tableaux $\tab$ and $\varepsilon^*(\tab)$ have the~same
    distribution. 
		An application of
    \cref{prop:intertwine} below completes the proof.
\end{proof}

In the remaining part of this section we will present the details of the~map~$\varepsilon^*$ and we will prove that \cref{prop:intertwine} indeed holds true.

\subsection{Dual evacuation}

The dual evacuation has a beautiful algorithmic
description in terms of the manipulations of the boxes of $\tab$,
cf.~\cite[Definition 2.10]{Pon2011}, however we will
not make use of it. For our purposes it is more convenient to define the dual
evacuation $\varepsilon^*$ implicitly by Robinson--Schensted--Knuth correspondence 
as follows. For a permutation $\sigma=(\sigma_1,\dots,\sigma_n)\in\Sym_{n}$ we denote
\[\sigma^\sharp:=(n+1-\sigma_n, \ldots, n+1-\sigma_1) \in \Sym_n.\]
If $\sigma$ corresponds to a pair $(P, Q)$ under RSK, then $\sigma^\sharp$
corresponds to $\big( \varepsilon^*(P), \varepsilon^*(Q)\big)$ under RSK, see
\cite[A1.2.10]{Stanley1999}.

\medskip

We will use the following fact (see \cite[Proposition~3.9.3]{Sagan2001} for the proof).
\begin{fact}[{\cite{Schuetzenberger1963}}]
    For any $\sigma \in \Sym_n$, the following identity holds up 
		to renumbering of the boxes on the left-hand side, 
		so that the resulting tableau becomes standard
    \begin{equation}\label{shift}
        \left(j \circ Q\right)(\sigma) 
				= \left(Q \circ \; s\right) (\sigma) 
				\qquad \text{for }\sigma = (\sigma_1, \ldots, \sigma_n) \in \Sym_n
    \end{equation} 
    where $s(\sigma) = s(\sigma_1, \sigma_2, \ldots, \sigma_n) := (\sigma_2, \ldots, \sigma_n)$ is a shift.
\end{fact}

\begin{proposition}
    \label{prop:intertwine}
    \[\pat(\tab)= \revevac\big(\varepsilon^*(T)\big).\]
\end{proposition}

\begin{proof}
    For any tableaux $R,S$ we will use a shorthand notation
    \[ R/ S = \sh R /\sh S \]
    for the skew diagram obtained by subtracting their shapes. In all examples
    below this skew diagram $ R/ S =\left\{ \Box\right\}$ consists of
    a single box; we will write shortly $\Box = R/ S$.
    
    Let $\tab = Q(\sigma)$ be a recording tableau of some permutation
    $\sigma$; with these notations $\varepsilon^*(\tab)=Q\left( \sigma^\sharp
    \right)$. 
    
    \medskip
    
    By \eqref{shift}, the lazy \jdtp fulfills for $i\in[n]$
    \begin{multline}
        \label{eq:evacuation-vs-jdt-A}
        \pat_{i}(\tab) = 
        Q\left(\sigma_1,  \ldots, \sigma_{i} \right) /
        j\big(  Q\left(\sigma_1, \ldots, \sigma_{i} \right) \big) =
        \\
        Q\left(\sigma_1, \sigma_2, \ldots, \sigma_{i} \right) / Q\left(\sigma_2, \ldots, \sigma_{i} \right).
    \end{multline}
    
    \medskip
    
    On the other hand, by \eqref{shift}, applying \jdtincomplete $n-i$ times to
    $\varepsilon^*(\tab)=Q\left( \sigma^\sharp \right)$, leads to the tableau (up
    to renumbering boxes on the left-hand side so that the tableau becomes
    standard)
    \[j^{n-i}\left( \varepsilon^*(\tab) \right) = Q\left(n + 1 -\sigma_{i}, \ldots,n + 1 -\sigma_2, n + 1-\sigma_{1}\right).\]
    The~position of~the~box with the maximal entry in a recording tableau can
    be found by comparing this tableau to the recording tableau of a truncated
    sequence; it follows that
    \begin{multline}
        \label{eq:evacuation-vs-jdt-B}
        \pos_n j^{n-i}\left( \varepsilon^*(\tab) \right) =
        Q\left(n + 1 -\sigma_{i}, \ldots,n + 1 -\sigma_2, n + 1-\sigma_{1}\right) / \\
        Q\left(n + 1 -\sigma_{i}, \ldots,n + 1 -\sigma_2 \right).        
    \end{multline}
    By the result of Schensted \cite[Lemma~7]{Schensted1961} and Greene theorem \cite[Theorem~3.1]{Greene1974} 
		the shapes of the tableaux which contribute to the right-hand sides of \eqref{eq:evacuation-vs-jdt-A} 
		and~\eqref{eq:evacuation-vs-jdt-B} are equal which concludes the proof.
\end{proof}

\subsection{Proof of \cref{thm:jdt}}
\label{sec:equivalence-of-main-problems}

\begin{proof}[Proof of \cref{thm:jdt}]
Let $N \in \N$. Let  $\Longitude_N \colon \tableaux_{\sq_N} \to [0,1]$ be the
random variable which is given by  \cref{thm:evacuation}.
We define the random variable $\LongitudeTilde_N \colon \tableaux_{\sq_N} \to [0,1]$ by 
\[\LongitudeTilde_N(T_N):= \Longitude_N(\varepsilon^*(T_N)).\]
Since  $\varepsilon^*$ is a bijection,
		\begin{multline*}
		\Pp \Big\{ \tab_N \in \tableaux_{\sq_N}: 
	\sup_{\tim \in [0,1)} 
	\left| \mypoint_\tim(\tab_N) -
	\point_{1-\tim, \Longitude_N(\tab_N)}  \right|
							 > \varepsilon \Big\} = \\
		\Pp \Big\{ \tab_N \in \tableaux_{\sq_N}: 
\sup_{\tim \in [0,1)} 
\left| \mypoint_\tim(\varepsilon^*(\tab_N)) -
\point_{1-\tim, \Longitude_N(\varepsilon^*(\tab_N))}  \right|
> \varepsilon \Big\} = \\
		\Pp \Big\{ \tab_N \in \tableaux_{\sq_N}: 
	\sup_{\tim \in [0,1)} 
	\left| \frac{1}{N} \pat_{\lceil (1-\tim) N^2\rceil}(\tab_N) -
	\point_{1-\tim, \LongitudeTilde_N(\tab_N)}  \right|
							 > \varepsilon \Big\} = \\
		\Pp \Big\{ \tab_N \in \tableaux_{\sq_N}: 
	\sup_{\tim \in (0,1]} 
	\left| \frac{1}{N} \pat_{\lceil \tim N^2\rceil}(\tab_N) -
	\point_{\tim, \Longitude_N(\tab_N)}  \right|
							 > \varepsilon \Big\},
	\end{multline*}
where the second equality is a consequence of \cref{prop:intertwine}. By
\cref{thm:evacuation} the left-hand side converges to $0$ in the limit $N \to
\infty$; on the other hand the right-hand side is the probability which appears
in \cref{thm:jdt}.
\end{proof}

\section{Generalizations of the main results \\ for non-square tableaux}
\label{sec:generalization}

\subsection{Continuous diagrams}
\label{subsec:continuous-diagrams}

We call a~function $\omega: \R \to \R$ a~\emph{continuous diagram} 
\cite{Kerov1993a,Kerov1998} if
\begin{itemize}
    \item $\omega$ is a $1$-Lipschitz function, i.e., 
    \[|\omega(u_1) - \omega(u_2)| \leq |u_1 - u_2| \qquad \text{for all $u_1, u_2 \in \R$;} \]
    \item $\omega(u) = |u|$ for sufficiently large $|u|$.
\end{itemize}
We will denote the set of continuous diagrams by $\mathcal{C Y}$;
we endow this set with the $L^\infty$-metric. 
(Our definition is more specific than the one of Kerov \cite{Kerov1993a}
who allows to additionally 
translate our \emph{centered} continuous diagrams
along the real line.)

Any (usual) Young diagram~$\lambda$ seen in 
the $(u,v)$-coordinate system is a $1$-Lipschitz function 
defined on some interval 
(given by the range of the $u$-coordinates of $\lambda$) 
and has slopes equal to $\pm 1$.
It can be extended outside its initial domain 
by a modulus function $x \mapsto |x|$.
In this way we obtain a continuous diagram $\omega_\lambda$
to which we will refer as \emph{the profile of $\lambda$}, 
see~\cref{fig:continuous-diagram}. 

Let $s > 0$. For a continuous diagram $\omega$ 
we define \emph{the scaling of $\omega$ by $s$}
as the following continuous diagram denoted by $s \omega$
\[
s \omega \colon \ 
\R \ni u \mapsto s \cdot \omega \! \left( s^{-1} u \right).
\] 

Let $(\lambda_N)$ be a sequence of Young diagrams
with the property that the sequence of 
the corresponding rescaled profiles
\[ \left(\frac{1}{\sqrt{|\lambda_N|}} \omega_{\lambda_N} \right)\] converges to 
a continuous diagram $\Lambda$
in the $L^\infty$-metric,
that is,
\[
\sup_{u \in \R} \left| 
\frac{1}{\sqrt{ |\lambda_N| }} \,
\omega_{\lambda_N} \! \left( \sqrt{ |\lambda_N| }\ u  \right)
- \Lambda(u)   \right| 
\xrightarrow{N \to \infty} 0.
\]
We will call $\Lambda$ \emph{the limit shape} for 
the sequence of Young diagrams $(\lambda_N)$ 
and denote such convergence by 
$\frac{1}{\sqrt{ \left|\lambda_N \right| }} 
\lambda_N
\to
\Lambda$.

\subsection{The asymptotic setup}
\label{sec:asymptotic-setup}

Let $C\geq 1$ be a fixed constant. For each integer $N\geq 1$ let $\lambda_N$ be
a $C$-balanced Young diagram.
We assume that
\[ 
\lim_{N \to \infty} \left|\lambda_N\right| = \infty 
\] 
and that there exists a limit shape $\Lambda \in \mathcal{C Y}$ 
for the sequence $(\lambda_N)$,
i.e., that $\frac{1}{\sqrt{\left|\lambda_N\right|}} \lambda_N \to \Lambda$.

Our goal in \cref{sec:generalization} is to find counterparts of
\cref{thm:evacuation,thm:jdt} in which the sequence $(\sq_N)$
of square diagrams is
replaced by the sequence $(\lambda_N)$ of $C$-balanced Young diagrams.

\subsection{The limit curves}

The result of Pittel and Romik concerning the existence of the level curves
\cite[Theorem~1(i)]{Pittel2007}, cf.~\cref{sec:levelcurves}, is a special case of
a more general phenomenon. Using the results of Biane \cite[Theorem~1.2 and
Theorem~1.5.1]{Biane1998} one can show that under the assumptions from
\cref{sec:asymptotic-setup} there exists a family of level curves for a uniformly
random Young tableau of the shape $\lambda_N$ (in the limit as $N\to\infty$).
The following proposition describes precisely this result. 

		\begin{proposition}
    \label{lem:limit-shape}
		Let $(\lambda_N)$ be a sequence of $C$-balanced Young diagrams 
		such that $|\lambda_N| \to \infty$ and 
		$\frac{1}{\sqrt{ \left|\lambda_N \right| }} \lambda_N
		\to \Lambda$ 
		for some continuous diagram $\Lambda \in \mathcal{C Y}$
		(i.e., $(\lambda_N)$ fulfills 
		the assumptions in \cref{sec:asymptotic-setup}).
		For any $\alpha \in [0, 1]$ 
    there exists a continuous diagram 
		$\Lambda_\alpha \in \mathcal{C Y}$
        such that
    \[ 
		\frac{1}{\sqrt{|\lambda_N|}} 
		\sh \! \left( \tab_N \big|_{\leq \lfloor \alpha \cdot |\lambda_N| \rfloor} \right) 
		\xrightarrow{P} 
		\Lambda_{\alpha}
		\]
		where $\tab_N$ is a uniformly random element of $\tableaux_{\lambda_N}$.
		\end{proposition}

		We will say that $\Lambda_\alpha$ is
		\emph{the $\alpha$-level curve} for the sequence $(\lambda_N)$
		or, shortly, the $\alpha$-level curve for the diagram $\Lambda$. 
		Note that $\Lambda_0$ is the empty diagram
		and $\Lambda_1 = \Lambda$.
		
		For example, in the case when $\lambda_N = \sq_N$ is a square Young diagram,
		the $\alpha$-level curve for $(\lambda_N)$ is 
		the curve $h_\alpha$,
		cf.~\cref{sec:levelcurves}.

		\begin{remark}
		Biane proved his results \cite[Theorem~1.2 and Theorem~1.5.1]{Biane1998} with
the tools of the free probability theory \cite{Mingo2017}, but
\cref{lem:limit-shape} can be also showed using the beam models \cite{Sun2018}
or by solving a gradient variational problem \cite{Kenyon2021}.
		\end{remark}

\subsection{The limit measures on the level curves}
\label{sec:limit-measure}

The second result of Pittel and Romik which gives explicitly the limit measure on
the $\alpha$-level curve for the square diagram \cite[Theorem~2]{Pittel2007},
cf.~\cref{sec:random-position}, is a special case of another general result. It
turns out that to every continuous diagram we can associate two natural measures
-- the \emph{transition measure} and the \emph{cotransition measure} -- which
have very natural interpretations in the case of the usual Young diagrams. 

		\subsubsection{Transition measure of a continuous diagram}
		\label{subsec:transition-measure}
		To any continuous diagram $\omega \in \mathcal{C Y}$, 
		one can associate a~probability measure $\tmeas_\omega$, 
		called \emph{the transition measure of~$\omega$} \cite{Kerov1993,Biane1998}, 
		as the unique compactly supported measure on~$\R$ 
		such that its Cauchy transform
		\[
		G_{\tmeas_{\omega}}(z) := \int_\R \frac{1}{z-x} \diff \tmeas_\omega (x)
		\]
		is given by the~equation
		\[
		G_{\tmeas_\omega}(z) 
		= \frac{1}{z} \exp \int_\R \frac{1}{x-z} \sigma'(x) \diff x 
		 = \frac{1}{z} \exp \int_\R \frac{1}{(x-z)^2} \sigma(x) \diff x 
		\]
		where $\sigma(u) := (\omega(u) - |u|)/2$. 
		
		The motivations for this notion are related to random walks on the set of Young diagrams:
		the atoms of the transition measure $\tmeas_{\omega_\lambda}$ of a usual Young diagram $\lambda$
		correspond to the Markov's transition probabilities in the Plancherel growth process starting in $\lambda$
        \cite[Section~3.2]{Kerov1993}.

		The mapping which to a continuous diagram $\omega$ 
		assigns the transition measure $\tmeas_\omega$
		is a \emph{homeomorphism} \cite[Section 2.3]{Kerov1993}. 
		Moreover, a continuous diagram is uniquely determined
		by its transition measure.

		\subsubsection{Cotransition measure of a continuous diagram}
		\label{subsec:cotransition-measure}
		
		For a continuous diagram $\omega \in \mathcal{C Y}$ we define its 
		\emph{area} as the area of the region between the profile and the $x$- and the $y$-axis: 
		\[
		A(\omega) = \int_\R \left( \omega(x) - |x| \right) \diff x.
		\]
		\emph{The cotransition measure} $\cotmeas_\omega$ of $\omega$ is defined as 
		the unique probability measure with the Cauchy transform $G_{\cotmeas_\omega}$
		given by the following equation \cite[Equation~(8)]{Romik2004}:
		\begin{equation}
		\label{eq:relation-measures}
		\frac{A(\omega)}{2} G_{\cotmeas_\omega}(x) = x - \frac{1}{G_{\tmeas_\omega}(x)}.
		\end{equation}
        By convention, the cotransition measure of the empty diagram 
        $\cotmeas_{\omega_\emptyset} = \delta_0$
        is defined to be the measure concentrated in $0$.    

		The cotransition measure $\cotmeas_{\omega_\lambda}$ 
		of a usual Young diagram $\lambda$
		is the distribution of the $u$-coordinate of the box with the maximal entry $|\lambda|$
		in a uniformly random standard tableau of shape $\lambda$,
		cf.~\cite[page 628 
		and the comment below Eq.~(6)]{Romik2004}.
		In particular, the measure $\cotmeas_\alpha$ 
		introduced in \cref{sec:random-position}
		to which we referred as \emph{the limit measure} 
		is the cotransition measure corresponding to 
		the continuous diagram $h_\alpha$
		(more precisely, to the proper extension 
		of the function $h_\alpha$ given by \eqref{height}
		by $x \mapsto |x|$ and $x \mapsto 2-|x|$). 
        Be advised that
		diagrams of different shape
		may have the same cotransition measure.
		For example, the square Young diagram $\sq_N$ 
		has the same cotransition measure $\cotmeas_{\omega_{\sq_N}} = \delta_0$
		concentrated at the point $u=0$, no matter which size of the square $N\in\N$ we choose.
		However, if we restrict our considerations
		to the set of (centered) continuous diagrams of fixed positive area,
		then any diagram is uniquely determined by its cotransition measure
		and this correspondence is a homeomorphism, 
		\cite[Theorem~6]{Romik2004}. 

		\subsubsection{The limit measures on the level curves of continuous diagram}
		\label{subsec:limit-measure-on-level-curves}

		We will refer to the cotransition measure 
		$\cotmeas_{\Lambda_\alpha}$ 
		corresponding to the $\alpha$-level curve $\Lambda_\alpha$ 
		for a continuous diagram $\Lambda$
		as \emph{the limit} (or \emph{cotransition}) \emph{measure} on the level curve $\Lambda_\alpha$. 
		
		The cumulative distribution function of the limit measure $\cotmeas_{\Lambda_\alpha}$	
		will be denoted by $F_{\Lambda_\alpha}$, i.e., 
		\[
		F_{\Lambda_\alpha}(u) := \cotmeas_{\Lambda_\alpha}\! \left( (-\infty, u] \right)
		\text{ for each $u \in \R$}.
		\]
		The density of the limit measure $\cotmeas_{\Lambda_\alpha}$	
		will be denoted by $f_{\Lambda_\alpha}$, whenever this density exists.

\subsection{The geographic coordinate system}
\label{subsec:geographic-coord-system2}

		We will endow a continuous diagram $\Lambda$ 
		with the system of geographic coordinates, cf.~\cref{sec:geographic}.
		For this purpose we view the shape $\Lambda$ 
		as the following compact subset of the $(x,y)$-Cartesian plane
		\begin{multline*}
		\Lambda^{\on{Cart}} := \overline{ \left\{ (x, y) \in \R^2 : \; |x-y| < x+y < \Lambda(x-y) \right\}   } = \\
		\overline{ \left\{ (x, y) \in [0,\infty)^2: \; x+y < \Lambda(x-y) \right\}   }.
		\end{multline*}
		(Recall that $ u = x-y $ and $v = x+y$ are the $(u,v)$ coordinates,
		cf.~\cref{sec:levelcurves}.)

	For any $\alpha \in [0,1]$ we define \emph{the quantile function} 
	for the limit measure $\cotmeas_{\Lambda_\alpha}$ 
	by the formula 
	\begin{equation}
	\label{eq:generalized-quantile}
	Q_{\Lambda_\alpha}(\psi) := \inf \! \left\{ 
	u \in \on{\supp}(\cotmeas_{\Lambda_\alpha}): \  
	F_{\Lambda_\alpha}(u) \geq \psi \right\}
	\quad \text{for $\psi \in [0,1]$}
	\end{equation}
	where $\on{supp}(\mu)$ denotes 
	\emph{the support} of the measure $\mu$.
	In particular, $Q_{\Lambda_0} \equiv 0$.

	For given $\alpha \in [0,1]$ and $\longitude\in [0,1]$ 
	there is exactly one point 
	$p = (x,y) \in \Lambda^{\on{Cart}}$
	such that 
	\begin{itemize}
	\item $p$ lies on the level curve $\Lambda_\alpha$
	seen in the $XY$-coordinates system, i.e., 
	\[ x+y = \Lambda_\alpha(x-y) ;\]
	\item the $u$-coordinate of $p$ is given by the quantile function:
    \[u(p) = x-y = Q_{\Lambda_\alpha}(\longitude) .\]
    \end{itemize}
    We will denote this point by 
	$\point_{\alpha, \longitude} := 
	(x_{\alpha}^\longitude, y_{\alpha}^\longitude) \in \Lambda^{\on{Cart}}$.
	In particular, by the definition of $\nu_{\Lambda_0}=\delta_0$ 
	we have $\point_{0, \psi} = (0,0) \in \R^2$
	for any $\psi \in [0,1]$. 
	Additionally we denote by 
\[
u_{\alpha}^{\longitude} := x_\alpha^\longitude - y_\alpha^\longitude
\quad \text{ and } \quad
v_{\alpha}^{\longitude} := x_\alpha^\longitude + y_\alpha^\longitude
\] 
	the~$u$- and $v$-coordinate of the point~$\point_{\alpha, \longitude}$. 
	
	We will refer to the mapping 
	$[0,1]^2 \ni (\alpha, \longitude) \mapsto \point_{\alpha, \longitude} \in \Lambda^{\on{Cart}}$
	as to \emph{the geographic coordinates system}.
	
	\begin{remark}
	There may be some problem with defining 
	the counterparts of the longitude and 
	the latitude (the geographic coordinates)
	of the point $p \in \Lambda^{\on{Cart}}$
	like we did in \cref{sec:geographic}.
    We defined \emph{the latitude} 
	as \emph{the unique} \mbox{$\alpha \in [0,1]$} for which
	$p$ lies on the level curve $h_{\alpha}$
	(more precisely, on the restriction of the level curve $h_\alpha$
	to the support of the corresponding cotransition measure~$\cotmeas_{\alpha}$).
	We are not sure if such uniqueness holds 
	in a general case
	when $\Lambda$ is an arbitrary continuous diagram.  
	To make things worse, 
	the definition of \emph{the longitude} 
	depends on the limit measure $\cotmeas_{\Lambda_{\alpha(p)}}$
	corresponding to the $\alpha(p)$-level curve 
	(circle of latitude with the latitude $\alpha(p)$, 
	cf.~\cref{sec:levelcurves,sec:geographic}).
	In the worst scenario, 
	not only we shall pick some latitude $\alpha(p)$,
	but also the limit measure $\cotmeas_{\Lambda_{\alpha(p)}}$ may have atoms. 
	In particular, an attempt of using these direct counterparts of 
	the definitions from \cref{sec:geographic} in the more general context
	may lead to the situation in which several points 
	have the same geographic coordinates
	or some geographic coordinates $(\alpha, \psi)$ are not used.
	The case of the $L$-shape diagram, see~\cref{fig:L-shape},
	is an example of the first problem.
	\end{remark}

	\begin{problem}
	Let $\Lambda$ be a continuous diagram. 
	Show that for any point 
	\[
	p \in \Lambda^{\on{UV}} := 
	\overline{ \left\{ (u, v) : \ u \in \R \ 
	\text{ and } \ \; |u| < v < \Lambda(u) \right\}   }
	\]
	there is a unique $\alpha \in [0,1]$
	with the property that 
	\[
	p \in \left\{(u, \Lambda_\alpha(u)): \ 
	u \in \on{supp}(\cotmeas_{\Lambda_\alpha}) \right\}.
	\]
	In other words, show that for any point $p \in \Lambda^{\on{UV}}$
	there exists a unique level curve $\Lambda_\alpha$ 
	(for some $\alpha \in [0,1]$) 
	which restricted to the support of
	the corresponding cotransition measure $\cotmeas_{\Lambda_\alpha}$ 
	contains $p$.
    \end{problem}

\subsection{Extension of the main results}

\begin{theorem}
    \label{thm:extension}
		
		Let $(\lambda_N)$ fulfill the assumptions in \cref{sec:asymptotic-setup})  
		and assume that 
        \begin{enumerate}[label=($\star$)]
            \item \label{eq:uniform-continuity-assumption} the geographic
coordinates system is continuous, i.e., the map 
$[0,1]^2 \ni (\alpha, \psi) \mapsto P_{\alpha, \psi}$, 
is continuous.
         \end{enumerate}

		Then the  analogues of \cref{thm:evacuation,thm:jdt} hold true, i.e., if
$\tab_N$ is a uniformly random element of $\tableaux_{\lambda_N}$ then:
		\begin{itemize}
		\item there exists a family of random variables 
		$\Psi_N: \tableaux_{\lambda_N} \to [0,1]$ 
		indexed by $N \in \N$
		such that 
	\begin{equation}
	\label{eq:generalized-main-theorem}
		\sup_{\tim \in [0,1]} 
					\left| \mypoint_\tim(\tab_N) -
					\point_{1-\tim, \Longitude_N(\tab_N)}  \right|
					\xrightarrow{P} 0,
	\end{equation}
		\item there exists a family of random variables 
		$\LongitudeTilde_N: \tableaux_{\lambda_N} \to [0,1]$ 
		indexed by $N \in \N$
		such that 
	\[
		\sup_{\tim \in [0,1]} 
				\left| 
				 \frac{1}{N} \pat_{\lceil \tim |\lambda_N| \rceil}(\tab_N)  -
				P_{\tim,\LongitudeTilde_N(\tab_N)} \right| 
				\xrightarrow{P} 0.
	\]
	\end{itemize}
	
	The probability distribution
	of the random variable~$\Longitude_N$ (respectively, $\LongitudeTilde_N$) 
	converges, as \mbox{$N\to\infty$}, to
the~uniform distribution on the~unit interval~$[0,1]$.

\end{theorem}

\begin{proof}
The proof of \cref{thm:evacuation} 
is applicable in this more general case --
one shall replace all occurrences of the $\sq_N$
by $\lambda_N$ and change the references to the limit measures. 
We enumerate the properties 
of square Young diagrams 
which played the crucial role 
in the proof of \cref{thm:evacuation}.
\begin{itemize}
\item $\sq_N$ is a $1$-balanced Young diagram
(we used this property in \cref{prop:single-and-multi,prop:longitude-experimental});
\item the cotransition measure $\cotmeas_\alpha$ 
corresponding to the level curve $h_\alpha$ 
has no atoms and has a connected support
(we used this property in \cref{thm:all-the-same,prop:single-and-multi});
\item the distribution of $u_\alpha$ -- the $u$-coordinate of the surfer 
in time $\alpha$ -- converges to the cotransition measure $\cotmeas_{1-\alpha}$
(we used this property in \cref{lem:true-position});
\item the uniform continuity of the geographic coordinates system, 
i.e., the uniform continuity of the mapping $(\alpha, \psi) \mapsto P_{\alpha, \psi}$,
cf.~\cref{lem:uniform-continuity}.
\end{itemize}
Notice that the counterparts of all these properties 
are present in our new setting.
In particular, the distribution of the $u$-coordinate of the surfer 
converges to the proper limit measure since the correspondence 
between the (centered) continuous diagrams of fixed positive area
and the cotransition measures is a homeomorphism, cf.~\cref{subsec:cotransition-measure}.
Moreover, the assumption~\ref{eq:uniform-continuity-assumption}
on the continuity of the geographic coordinates system
assures that for any $\alpha \in [0,1]$
the support of the limit measure $\cotmeas_{\Lambda_\alpha}$ is connected.
\end{proof}

\subsection{Example: random rectangular tableaux}
\label{subsec:rectangular}

		For any real numbers \mbox{$a,b > 0$} let $\sq_{a \times b}$ denote 
		the rectangle with the left bottom corner positioned in $(0,0) \in \R^2$
		and with sides $a$ and $b$
		(on $X$ and $Y$ axis, respectively, when seen in the $(x,y)$-coordinates system).
		
    Let $(M_i)$ and $(N_i)$ be two sequences of positive integers which fulfill
    the conditions from \cref{sec:setup}, i.e., 
		$M_i \to \infty$ and $N_i \to \infty$, and 
		there is some (shape parameter) $\theta > 0$ 
		such that
		\[ 
		\lim_{i\to\infty} \frac{M_i}{N_i} = \theta.
		\]
		We define for $i \in \N$
    \[ 
		\lambda_i := \sq_{M_i \times N_i} 
		\]
    to be the rectangular Young diagram which has $M_i$ rows and $N_i$ columns.
		
		The following theorem generalizes \cref{thm:jdt}
		(which is a special case for $M_i=N_i=i$).
		
		\begin{corollary}
		\label{thm:rectangle-evacuation}
		For $i \in \N$ let $\tab_i$ be a uniformly random tableau of 
		the shape $\sq_{M_i \times N_i}$. 
		Then the rescaled lazy sliding path 
		$\frac{1}{\sqrt{M_i N_i}} \pat(\tab_i)$
		with respect to the supremum norm converges in probability
		to the random function
		\begin{equation}
		\label{eq:u-coord-rectangle}
		\tim \mapsto \traj_S(\tim)
		= 2 \sqrt{\tim(1-\tim)}\ S + \frac{\theta-1}{\sqrt{\theta}}\ \tim
		\end{equation}
		where $S$ denotes the random variable 
		with the standard semicircular distribution, cf.~\eqref{eq:SC}.
	\end{corollary}        
        
		\begin{proof}
		We will check that the assumptions of \cref{thm:extension}
		are fulfilled and find the explicit formula 
		for the geographic coordinates system.
		
		Clearly, the sequence 
		$\left( \frac{1}{\sqrt{M_i N_i}} \sq_{M_i \times N_i} \right)$
		converge to the limit shape $\Lambda = \sq_{\sqrt{\theta} \times 1/\sqrt{\theta}}$.
		We apply \cref{lem:rectangle} below with 
		$a := 1/\sqrt{\theta}$ and $b:= \sqrt{\theta}$ 
		to see that 
		\begin{itemize}
		\item the geographic coordinates system on 
		$\sq_{\sqrt{\theta} \times 1/\sqrt{\theta}}$ 
		is continuous,
		\item for any $\psi \in [0,1]$ and $\alpha \in [0,1]$
		the $\psi$-th quantile $u_\psi$ of 
		the cotransition measure $\cotmeas_{\Lambda_\alpha}$ 
		is given by 
		\[
		u_\psi = 2\sqrt{\alpha(1-\alpha)} F_{\operatorname{SC}}^{-1}(\psi) 
		+ \frac{\theta - 1}{\sqrt{\theta}} \alpha
		\]
		where $F_{\operatorname{SC}}$ denotes 
		the cumulative distribution function 
		of the standard semicircle distribution, cf.~\eqref{eq:SC}.
\end{itemize}
		
		Notice that if $U$ is a random variable 
		with the uniform $U(0,1)$ distribution
		then $F_{\operatorname{SC}}^{-1}(U)$
		has the standard semicircle distribution.
		
		By \cref{thm:extension} (its second and third part)
		the (rescaled) lazy sliding path 
		$\frac{1}{\sqrt{M_i N_i}} \pat(\tab_i)$
		with respect to the supremum norm
		converges in probability to 
		the random function \eqref{eq:u-coord-rectangle}. 
		\end{proof}

		\begin{lemma}
		\label{lem:rectangle}
		Let $a, b > 0$ and $\Lambda := \omega_{\sq_{a \times b}}$ 
		be the profile of\/ $\sq_{a \times b}$.
		Then for any $\alpha \in (0,1)$ 
		the cotransition measure $\cotmeas_{\Lambda_\alpha}$ 
		corresponding to the $\alpha$-level curve~$\Lambda_\alpha$ 
		has the density
		\begin{equation}
		\label{eq:rectangle-density}
		f_{\cotmeas_{_{\Lambda_\alpha}}} \! (x) 
		= 
		\frac{1}{2 \sqrt{ab \cdot \alpha(1-\alpha)}} \,
		f_{\operatorname{SC}} 
		\! \left( \frac{x + \alpha(a-b)}{2 \sqrt{ab \cdot \alpha(1-\alpha)}}  \right)
		\end{equation}
		where $f_{\operatorname{SC}}$ is the density 
		of the standard semicircular distribution, cf.~\eqref{eq:SC}.
		\end{lemma}
		
		\begin{proof}
		The following calculations use \eqref{eq:relation-measures}
		and some relations
		between the $R$-transform and the Cauchy transform of 
		the appropriate transition measures, 
		see \cite[Section~3]{Mingo2017} for the theory. 
		
		The Cauchy transform of the transition measure $\tmeas_\Lambda$ 
        of the rectangular diagram
		is given by \cite[Equation~(2)]{Romik2004}
		\begin{equation}
		\label{eq:G-rectangular}
		G(z) := G_{\tmeas_\Lambda} (z) = \frac{ z - (b-a) }{(z+a) (z-b)}.
		\end{equation}
		It is an analytic function and in some neighborhood of $\infty$
		it is invertible \cite[Section~3, Theorem~17(i)]{Mingo2017}. 
		Moreover, there exists a neighborhood $U \subset \C$ of~$0$
		for which 
		$G|_{G^{-1}(U)}$ is invertible \cite[Section~3, Theorem~17(ii)]{Mingo2017}
		and we can calculate the $R$-transform of the measure $\tmeas_\Lambda$ 
		with the formula \cite[Section~3, Theorem~17(iii)]{Mingo2017}
		\[
		R(z) := R_{\tmeas_\Lambda}(z) = G^{-1}(z) - \frac{1}{z}
		\text{\quad for $z \in U \setminus\{0\}$}.
		\]
		Substituting in the latter $z$ with $G(z')$ 
		(for some $z' \in G^{-1}(U \setminus\{0\})$) 
		we get the relation
		\begin{equation}
		\label{eq:G-relation}
		z = R\big(G(z)\big) + \frac{1}{G(z)}
		\text{\quad for $z \in G^{-1}(U \setminus\{0\})$}.
		\end{equation}
		We use this equality in order to substitute each occurrence of the variable $z$
on the right-hand side of \eqref{eq:G-rectangular}; by clearing of the
denominator we obtain
		\begin{equation}
		\label{eq:G-R-equation}
		R(G) + \frac{1}{G} - (b-a) 
		= G\cdot \! \left( R(G) + \frac{1}{G} + a \right)\! \left( R(G) + \frac{1}{G} - b \right), 
		\end{equation}
        where we used the shorthand notation $G=G(z)$.
		By the choice of $U$ as the neighborhood of $0$ 
		we get that \eqref{eq:G-R-equation} is fulfilled with $G$ 
		replaced by any complex number $z \in U \setminus \{0\}$,
		i.e., the $R$-transform fulfills the following quadratic equation 
        for any $z \in U \setminus \{0\}$
		\begin{equation}
		\label{eq:R-equation}
		R(z) + \frac{1}{z} - (b-a)
		= z \cdot \! \left( R(z) + \frac{1}{z} + a \right)\! \left( R(z) + \frac{1}{z} - b \right). 
		\end{equation}	
		
        \medskip
		
		Now, we will calculate the Cauchy transform 
		of the transition measure~$\tmeas_{\Lambda_\alpha}$ 
		on the level curve $\Lambda_\alpha$.
		Denote by $R_\alpha$ and $G_\alpha$, respectively,
		the $R$-transform and the $G$-transform 
		of~$\tmeas_{\Lambda_\alpha}$.
		The $R$-transforms of the transition measures 
		$\tmeas_\Lambda$ and $\tmeas_{\Lambda_\alpha}$ 
		are related by the following correspondence
		\cite[Theorem~1.2]{Biane1998}
		\[
		R_\alpha(z) := R_{\tmeas_{\Lambda_\alpha}} \! (z) = R( \alpha \cdot z )
		\text{\quad for $z \in U \setminus \{0\}$}.
		\]
		Let us put $\alpha \cdot z$ instead of $z$ in \eqref{eq:R-equation}
		(this substitution is legal since 
		the neighborhood~$U$ can be taken to be a convex set).
		Equation \eqref{eq:R-equation} implies therefore that $R_\alpha(z)$ is a solution to the following quadratic equation
		\begin{multline}
		\label{eq:R-alpha-equation}
		R_\alpha(z) + \frac{1}{\alpha z} - (b-a) = \\
		\alpha z \left( R_\alpha(z) + \frac{1}{\alpha z} + a \right)
		\! \left( R_\alpha(z)+ \frac{1}{\alpha z} - b \right)
		\end{multline}
		for any $z \in U \setminus \{0\}$.

        In \eqref{eq:R-alpha-equation} we substitute each occurrence of the
variable $z$ by $G_\alpha(z)$; this substitution is valid as long as
$|z|$ is big enough so that $G_\alpha(z)\in U \setminus \{0\}$. Let us
denote additionally $H_\alpha = \frac{1}{G_\alpha}$; then we use the
relation \eqref{eq:G-relation} and substitute each occurrence of
$R_\alpha\big(G_\alpha(z)\big)$ by $z-H_\alpha(z)$. Then
\eqref{eq:R-alpha-equation} takes the form
		\begin{multline}
		\label{eq:G-H-level-curve-equation}
		H_\alpha(z) \! \left(z - (b-a) + \left(\frac{1}{\alpha} - 1 \right) H_\alpha(z) \right)
		= \\
		\alpha
		\left(z + a + \left(\frac{1}{\alpha} - 1 \right) H_\alpha(z) \right)
		\! \left(z - b + \left(\frac{1}{\alpha} - 1 \right) H_\alpha(z) \right)
		\end{multline}
		which holds if $|z|$ is big enough. %
For any $z$ big enough, the latter is a quadratic equation in $H_\alpha(z)$ 
which has two solutions given by
explicit (but complicated, so we omit writing them here) formulas. 
These solutions come from two branches of the complex square root. 
The function $H_\alpha$ must be analytic in some neighborhood of $\infty$ 
and therefore can be given by only one (family) of these solutions. 
Moreover, $H_\alpha$ must have a proper asymptotics, more precisely 
\cite[Section~3.1, Lemma~3]{Mingo2017}
\[
\lim_{y \to \infty} \frac{H_\alpha(iy)}{y} = i,
\]
which allows us to choose the proper solution. 

The formula for $H_\alpha$ gives also an explicit formula for 
the Cauchy transform of
the cotransition measure $\cotmeas_{\Lambda_\alpha}$ on the level curve
$\Lambda_\alpha$ as (cf.~\eqref{eq:relation-measures})
		\[
		G_{\cotmeas_{_{\Lambda_\alpha}}} \! (z) 
		= \frac{1}{\alpha ab} \! \left( z - H_\alpha(z) \right).
		\]
The function $G_{\cotmeas_{_{\Lambda_\alpha}}}$ is analytic (since $H_\alpha$ is analytic);
we will use its analytic continuation to the upper halfplane $\C^+$.
		
		With this (complicated) formula for $G_{\cotmeas_{_{\Lambda_\alpha}}}$ we can now recover 
		the density of the cotransition measure $\cotmeas_{\Lambda_\alpha}$
		using the Stieltjes inversion formula \cite[Section~3, Theorem~6]{Mingo2017}.
		One can easily show that this density is the properly rescaled 
		and translated standard semicircle distribution, cf.~\eqref{eq:SC}.
		\end{proof}

\subsection{What if the geographic coordinates system is not uniformly continuous?}

		The situation when the geographic coordinates system
		$(\alpha, \psi) \mapsto P_{\alpha, \psi} \in \Lambda^{\on{Cart}}$ 
		(recall~\cref{subsec:geographic-coord-system2})
		is not uniformly continuous is not rare. 
		One among many examples is 
		the limit shape $\Lambda$
		which is the $L$-shape, cf.~\cref{fig:L-shape}. 
		In this case there are some $\alpha$-level curves 
		(for $\alpha$ big enough)
		for which the corresponding 
		cotransition measure is supported 
		on two disjoint intervals. 
		This forces 
		the function $\psi \mapsto P_{\alpha, \psi}$
		to have a discontinuity related to the hole between the intervals.
		Therefore in such situation \cref{thm:extension} does not apply
		since the assumption~\ref{eq:uniform-continuity-assumption} is not fulfilled.      

		\begin{problem}
		Find a counterpart of the assumption~\ref{eq:uniform-continuity-assumption} 
		for which the conclusion of \cref{thm:extension} holds true.
		\end{problem}

\section{The correspondence between Young tableaux and 
particle systems. Proof of \cref{thm:second-class-B}}
\label{sec:particles}

We will first describe a bijection
which links 
a recording tableau 
with a unique history of 
a particular Totally Asymmetric Simple Exclusion Process,
recall \cref{sec:scp}.
We will base on the articles of Rost \cite{Rost1981}, as well as 
Romik and the second named author \cite[Section~7]{Romik2015a}. 
Then, in \cref{subsec:proof-second-class-particle}
we will prove \cref{thm:second-class-B}.

\subsection{The correspondence between a Young diagram 
with a distinguished corner
and a configuration of particles -- Rost's mapping}
\label{sec:Rost-mapping}

Let $\lambda$ be a non-empty Young diagram 
and $\Box$ be one of its inner corners,
i.e., $\Box$ is a cell of $\lambda$ 
such that the shape $\lambda \setminus \Box$
is still a Young diagram, 
see \cref{fig:Rost-mapping}.
Following \cite[Section~7.1]{Romik2015a}, 
we will present the two-step algorithm in which 
to the pair $(\lambda, \Box)$ 
we assign a configuration of holes and particles 
with exactly one second class particle. 
To our best knowledge the foundations for this mapping 
were first laid in \cite[Remark~1]{Rost1981},
and therefore we will call it \emph{Rost's mapping}.

\begin{figure}[tbp]
    \begin{center}
        \begin{tikzpicture}[scale=0.7]
            \draw[->,very thick] (-7,0)  -- (8,0) node[anchor=west]{{{$u$}}};
            \draw[->,very thick] (-7,-2) -- (8,-2) node[anchor=west]{{{$u$}}};
            \clip (-7,-2.5) rectangle  (7-0.25,7-0.25) ;
            \draw[fill=blue!10](0,0) -- (4,4) -- (2,6) -- (0,4) -- (-1,5) -- (-3,3) -- cycle;
            \draw[thin, dotted] (-6,-2) grid[xstep=1,ystep=100] (20,20);
            \foreach \x in {-6, -5, -4, -3, -2, -1, 0, 1, 2, 3, 4, 5, 6} { \draw[thick] (\x, -4pt)  -- (\x, 4pt); } 
            \foreach \x in {-5, -4, -3, -2, -1, 0, 1, 2, 3, 4, 5} {
                \draw[thick] (\x, -2) +(0,-0.4)  -- +(0,0.4);
                \node[fill=white] at (\x,-1) {\tiny $\x$}; }
            \draw[dashed](0,0) -- (8,8);
            \draw[dashed](-1,1) -- (7,9);
            \draw[dashed](-2,2) -- (6,10);
            \draw[dashed](-3,3) -- (5,11);
            \draw[dashed](-4,4) -- (4,12);
            \draw[dashed](-5,5) -- (3,13);
            \draw[dashed](-6,6) -- (2,14);
            \draw[dashed](-7,7) -- (1,15);
            \draw[dashed](0,0) -- (-8,8);
            \draw[dashed](1,1) -- (-7,9);
            \draw[dashed](2,2) -- (-6,10);
            \draw[dashed](3,3) -- (-5,11);
            \draw[dashed](4,4) -- (-4,12);
            \draw[dashed](5,5) -- (-3,13);
            \draw[dashed](6,6) -- (-2,14);
            \draw[dashed](7,7) -- (-1,15);
            \fill[pattern=north west lines, pattern color=red](2,4) -- (1,5) -- (2,6) -- (3,5) -- cycle;
            \draw[](8,8) -- (0,0) -- (-8,8);
            \fill[pattern=north west lines, pattern color=red,ultra thick] (1.1,-0.5) rectangle (2.9,0.5);
            \draw[ultra thick,black](7,7) -- (4,4) -- (2,6) -- (0,4) -- (-1,5) -- (-3,3)  -- (-7,7);
            \fill[fill=blue!50,draw=black] (3.5,4.5) circle(0.3);
            \fill[fill=blue!50,draw=black] (2.5,5.5) circle(0.3);
            \fill[fill=blue!50,draw=black] (-0.5,4.5) circle(0.3);
            \fill[fill=blue!50,draw=black] (-3.5,3.5) circle(0.3);
            \fill[fill=blue!50,draw=black] (-4.5,4.5) circle(0.3);
            \fill[fill=blue!50,draw=black] (-5.5,5.5) circle(0.3);
            \fill[fill=blue!50,draw=black] (3.5,0) circle(0.3);
            \fill[fill=blue!50,draw=black] (2.5,0) circle(0.3);
            \fill[fill=blue!50,draw=black] (-0.5,0) circle(0.3);
            \fill[fill=blue!50,draw=black] (-3.5,0) circle(0.3);
            \fill[fill=blue!50,draw=black] (-4.5,0) circle(0.3);
            \fill[fill=blue!50,draw=black] (-5.5,0) circle(0.3);
            \fill[fill=white,draw=black] (-2.5,3.5) circle(0.3);
            \fill[fill=white,draw=black] (-1.5,4.5) circle(0.3);
            \fill[fill=white,draw=black] (0.5,4.5) circle(0.3);
            \fill[fill=white,draw=black] (1.5,5.5) circle(0.3);
            \fill[fill=white,draw=black] (4.5,4.5) circle(0.3);
            \fill[fill=white,draw=black] (5.5,5.5) circle(0.3);
            \fill[fill=white,draw=black] (-2.5,0) circle(0.3);
            \fill[fill=white,draw=black] (-1.5,0) circle(0.3);
            \fill[fill=white,draw=black] (0.5,0) circle(0.3);
            \fill[fill=white,draw=black] (1.5,0) circle(0.3);
            \fill[fill=white,draw=black] (4.5,0) circle(0.3);
            \fill[fill=white,draw=black] (5.5,0) circle(0.3);
            \fill[fill=white,draw=black] (-2,-2) circle(0.3);
            \fill[fill=white,draw=black] (-1,-2) circle(0.3);
            \fill[fill=white,draw=black] (1,-2) circle(0.3);
            \fill[fill=white,draw=black] (4,-2) circle(0.3);
            \fill[fill=white,draw=black] (5,-2) circle(0.3);
            \fill[fill=blue!50,draw=black] (3,-2) circle(0.3);
            \draw[preaction={fill,white},pattern=north west lines, pattern color=red] (2,-2) circle (0.3); 
            \fill[fill=blue!50,draw=black] (-0,-2) circle(0.3);
            \fill[fill=blue!50,draw=black] (-3,-2) circle(0.3);
            \fill[fill=blue!50,draw=black] (-4,-2) circle(0.3);
            \fill[fill=blue!50,draw=black] (-5,-2) circle(0.3);
            \draw[red,ultra thick] (1.1,-0.5) rectangle (2.9,0.5);

        \end{tikzpicture}
    \end{center}
    \caption{Above:~the Young diagram $(4,4,2)$ and its marked inner
        corner in the second row. In the middle: the corresponding configuration of
        particles on the shifted lattice $\Z'$ via Rost's mapping. The marked inner
        corner corresponds to the pair of nodes in the rectangle.  Below:~the
        corresponding configuration of particles (including the second class particle) on
        the lattice $\Z$.} 
		\label{fig:Rost-mapping}
\end{figure}

In the first step of the Rost's mapping, given a Young diagram $\lambda$ we draw 
its profile $\omega_\lambda$ in the Russian coordinates system
(see \cref{fig:continuous-diagram} 
and \cref{subsec:continuous-diagrams} for 
the definition of the profile).
To the profile $\omega_\lambda$ there corresponds 
a unique configuration of holes and (first class) particles 
on $\Z' := \Z + \frac{1}{2}$ 
which appears in the following way.
For each $m \in \Z$ exactly one of the following two cases holds true:
\begin{itemize}
\item in the case when the slope of the profile $\omega_\lambda$ on the interval $[m, m+1]$ 
is equal to $-1$ then we put a (first class) particle 
at the site $m+\frac{1}{2}$ of the lattice $\Z' := \Z + \frac{1}{2}$; 
\item in the case when the slope of the profile $\omega_\lambda$ on the interval $[m, m+1]$ 
is equal to $+1$ then we put a hole at the site $m+\frac{1}{2}$
(in other words, the site $m+\frac{1}{2}$ is vacant). 
\end{itemize}
The distinguished inner corner~$\Box$ 
corresponds in the above particle system 
to a \emph{hole--particle pair}.
We outline this hole--particle pair with a red rectangle, 
see~the top and the middle part of \cref{fig:Rost-mapping}.

In the second step of Rost's mapping 
we define a particle configuration on~$\Z$.
We start with merging the hole--particle pair outlined in the red rectangle
into the single particle which we will call \emph{the~second class particle}.
We put it in the middle of the initial interval containing the hole-particle pair. 
Then we translate
all holes and particles which are placed to the left of the second class particle 
by $+\frac{1}{2}$
and we translate
all holes and particles which are placed to the right of the second class particle
by $-\frac{1}{2}$. 
These steps are illustrated in the middle and bottom part of \cref{fig:Rost-mapping}.
In this way we end up with the configuration of particles on $\Z$
with a single second class particle.

\subsection{The correspondence between the standard tableaux and the histories of TASEP}
\label{sec:enhanced-Rost-mapping}

Recall that for a standard tableau $\tab$ and 
a positive integer \mbox{$p \leq |T|$} 
we define the restricted tableau $\tab|_{\leq p}$ 
to be the tableau which consists of only these boxes of $\tab$
which have entries $\leq p$. 

For a standard tableau $\tab$ with $n \geq 1$ boxes 
and any $p \in \{1, \dots, n\}$ let us define 
\begin{align*}
\lambda^{(p)} & := \lambda^{(p)}(\tab) := \on{sh} \! \left( \tab|_{\leq p} \right),
\\
\Box^{(p)}    & := \Box^{(p)}(\tab) := \pat_{p}(\tab),
\end{align*} 
that is, $\big(\lambda^{(p)}, \Box^{(p)} \big)$ is 
the pair which consists of 
the Young diagram $\lambda^{(p)}(\tab)$ 
which is the shape of the restricted tableau $\tab|_{\leq p}$
and the last box along the sliding path in $\tab$ 
which contains a number $\leq p$, cf.~\cref{sec:jdt}. 

Let $\tab$ be a standard tableau with $n \geq 1$ boxes.
For any $t \in \{1, \dots, n\}$ let us consider 
the system of particles $P_t := P_t(\tab)$ 
which corresponds to the pair $(\lambda^{(t)}, \Box^{(t)})$
via Rost's mapping defined in \cref{sec:Rost-mapping}.
Notice that the initial configuration $P_1$ 
is such that 
the second class particle is located at the site $u=0$,
all negative nodes are occupied by the first class particles
and all positive nodes are occupied by holes.
Such configuration is called \emph{the Dirac sea}. 
Observe that for any $t \in \{1, \dots, n-1\}$ 
the neighboring states $P_t$ and $P_{t+1}$ 
differ by 
one of the three transitions  
described in \cref{sec:setup-scp} (cf.~\cref{fig:transition},
see \cite[Sections~7.2 and 7.3]{Romik2015a} for a step-by-step proof). 
Therefore the family $(P_t)_{t \in \{1, \dots, n\}}$ is 
a history of the particle system starting from the Dirac sea.

Given the above observation, it is easy to see that the mapping
which to the standard tableau $\tab$ with $n$ boxes 
associates the history $\left( P_t(\tab) \right)_{t \in \{1, \dots, n\}}$
of the particle system is 
\emph{a bijection} between the set of standard tableaux with $n$ boxes
and the possible $n$-step histories of particle systems 
starting from the Dirac sea.
Moreover, for a tableau $\tab$
and any $t \in \{1, \dots, |\tab|\}$
the $u$-coordinate 
of the box $\pat_t(\tab)$ in the lazy \jdtp
is the position of the 
second class particle in time $t$ 
in the corresponding TASEP,
see \cite[Proposition~7.1]{Romik2015a}.

\subsection{Proof of \cref{thm:second-class-B}}
\label{subsec:proof-second-class-particle}

\begin{proof}[Proof of \cref{thm:second-class-B}]
Let us denote by $\tableaux_{M \times N}$ the set of standard tableaux 
of $M \times N$ rectangular shape (where $M$ is a number of rows
and $N$ is a number of columns).

In the following we discard the particles and holes which do not occupy the nodes
\eqref{eq:set-of-nodes}. In this way the correspondence defined in
\cref{sec:enhanced-Rost-mapping} gives a bijection between $\tableaux_{M_i \times
    N_i}$ and the set of histories of the particle system considered in
\cref{sec:setup-scp}.
In this correspondence for any $t \in \{1, \dots, M_i N_i \}$ the position of the
second class particle in time $t$ in the TASEP corresponds to the $u$-coordinate
of the box $\pat_t(\tab)$ in the lazy \jdtp. In the special case $M_i=N_i=i$ when
$M_i \times N_i = \sq_i$ an application of \cref{thm:jdt} completes the proof. In
the general case we apply \cref{thm:rectangle-evacuation} instead.
\end{proof}

\section*{Acknowledgments} 
Research is supported by \emph{Narodowe Centrum Nauki}, grant number 
\linebreak
2017/26/A/ST1/00189. 
We thank Dan Romik and Serban Belinschi
for valuable discussions. 

\printbibliography

\end{document}